\documentclass[12pt,a4paper,reqno]{amsart}
\usepackage[a4paper,top=3cm,bottom=3cm,inner=3cm,outer=3cm,includehead,includefoot]{geometry}
\usepackage{amsmath,amssymb,amsthm,mathtools,wasysym,calc,verbatim,enumitem,tikz,pgfplots,url,hyperref,mathrsfs,cite,fullpage,mathtools}
\usetikzlibrary{shapes.misc,calc,intersections,patterns,decorations.pathreplacing,backgrounds}
\hypersetup{colorlinks=true,citecolor=blue,linktoc=page}
\usepackage{setspace}

\newtheorem{theorem}{Theorem}

\newtheorem{lemma}[theorem]{Lemma}
\newtheorem{prop}[theorem]{Proposition}
\newtheorem{claim}[theorem]{Claim}

\theoremstyle{definition}
\newtheorem{definition}[theorem]{Definition}
\newtheorem{remark}[theorem]{Remark}

\newenvironment{clmproof}[1]{\begin{proof}[Proof of Claim~\ref{#1}]\let\qednow\qedsymbol\renewcommand{\qedsymbol}{}}{\; \qednow \end{proof}}

\numberwithin{theorem}{section}
\mathtoolsset{showonlyrefs}
\setlist[itemize]{leftmargin=1cm}
\setlist[enumerate]{leftmargin=1.3cm}
\newlist{itemize2}{itemize}{1}
\setlist[itemize2]{labelindent=0.6cm, labelwidth=\widthof{(M2)}, labelsep=0.3cm, leftmargin=!}

\everydisplay{\thickmuskip=7mu}
\renewcommand{\leq}{\leqslant}
\renewcommand{\geq}{\geqslant}
\renewcommand{\le}{\leqslant}
\renewcommand{\ge}{\geqslant}
\renewcommand{\to}{\rightarrow}
\newcommand{\wdiam}{\operatorname{wd}}
\newcommand{\diam}{\operatorname{diam}}
\newcommand{\ice}{\operatorname{Ice}^\times}

\newcommand{\interior}{\operatorname{int}}
\newcommand{\conv}{\operatorname{conv}}
\newcommand{\sgn}{\operatorname{sgn}}

\let\eps\varepsilon

\def\B{\mathcal{B}}
\def\BB{\mathbf{B}}
\def\C{\mathcal{C}}

\def\Ex{\mathbb{E}}

\def\H{\mathbb{H}}
\def\HH{\mathbf{H}}
\def\cH{\mathcal{H}}

\def\J{\mathcal{J}}
\def\K{\mathcal{K}}

\def\N{\mathbb{N}}
\def\cN{\mathcal{N}}

\def\P{\mathbb{P}}
\def\Pr{\mathbb{P}}

\def\R{\mathbb{R}}

\def\S{\mathcal{S}}
\def\SS{\mathbf{S}}
\def\T{\mathcal{T}}
\def\U{\mathcal{U}}

\def\X{\mathcal{X}}
\def\Y{\mathcal{Y}}
\def\Z{\mathbb{Z}}

\def\<{\langle}
\def\>{\rangle}

\def\0{\mathbf{0}}
\def\ih{\mathrm{IH}}

\title[Universality for monotone cellular automata]{Universality for Monotone Cellular Automata}

\author[P. Balister \and B. Bollob\'as \and R. Morris \and P.J. Smith]{Paul Balister \and B\'ela Bollob\'as \and Robert Morris \and Paul Smith}

\address{Mathematical Institute, University of Oxford, Radcliffe Observatory Quarter, Woodstock Road, Oxford, OX2 6GG, UK}\email{Paul.Balister@maths.ox.ac.uk}

\address{Department of Pure Mathematics and Mathematical Statistics, Wilberforce Road, Cambridge, CB3 0WA, UK, and Department of Mathematical Sciences, University of Memphis, Memphis, TN 38152, USA}
\email{b.bollobas@dpmms.cam.ac.uk}
\address{IMPA, Estrada Dona Castorina 110, Jardim Bot\^anico, Rio de Janeiro, 22460-320, Brazil}
\email{rob@impa.br}
\address{Clerkenwell, London}
\email{paulsmith@cantab.net}
\thanks{P.B.\ and B.B.\ were partially supported by NSF grant DMS~1855745, R.M.\ by FAPERJ (Proc.~E-26/202.993/2017) and CNPq (Proc.~303681/2020-9), and by ERC Starting Grant 680275 MALIG, and P.S. by ERC Starting Grant 676970, by ISF grants 1147/14 and 1207/15, and by a CNPq bolsa PDJ} 

\date{\today}

\begin{document}

\begin{abstract}
In this paper we study monotone cellular automata in $d$ dimensions. We develop a general method for bounding the growth of the infected set when the initial configuration is chosen randomly, and then use this method to prove a lower bound on the critical probability for percolation that is sharp up to a constant factor in the exponent for every `critical' model. This is one of three papers that together confirm the Universality Conjecture of Bollob\'as, Duminil-Copin, Morris and Smith. 
\end{abstract}

\maketitle

\tableofcontents

\section{Introduction}

The Universality Conjecture of Bollob\'as, Duminil-Copin, Morris and Smith~\cite{BDMS} states that every $d$-dimensional monotone cellular automaton (see Definition~\ref{def:MCA}) is a member of one of $d+1$ `universality classes', which are characterized by their behaviour on sparse random sets. More precisely, it states that the critical probability for percolation in $\Z_n^d$ of each such automaton is either bounded away from zero, or is of order $\big( \log_{(r-1)} n \big)^{-\Theta(1)}$ for some $r \in \{1,\ldots,d\}$, where $\log_{(r)}$ denotes an $r$-times iterated logarithm. The automata with non-trivial phase transitions were characterized in ~\cite{BBMSupper,BBMSsub}, and for those with trivial phase transitions an upper bound on the critical probability was proved in terms of an explicit parameter $r = r(\U)$. In this paper we 
prove a matching lower bound, and hence complete the proof of the Universality Conjecture. 

\pagebreak

The mathematical study of general cellular automata was initiated over 70 years ago,\footnote{However, according to von Neumann~\cite{vN51}, it was already the case in 1951 that `Automata have been playing a continuously increasing, and have by now attained a very considerable, role in the natural sciences. This is a process that has been going on for several decades.'} by von Neumann (see~\cite{vN51,vN66,S58,U58}) and Ulam~\cite{U50}, who were inspired by the then-incipient theory of computing, and by the model of neural networks that had been introduced a few years earlier by McCulloch and Pitts~\cite{MP43}. They can be used to simulate the behaviour of a wide variety of biological, chemical and physical processes, and have become important objects of study in many areas of applied mathematics. For example, some of the most extensively studied models in statistical physics are (probabilistic) cellular automata, including the (kinetic) Ising model of ferromagnetism~\cite{Hugo20,FV17,M99} and kinetically constrained models of the liquid-glass transition~\cite{CMRT2,GST,T22} (see Sections~\ref{KCM:sec} and~\ref{Ising:sec}). Perhaps surprisingly, an important first step in understanding the metastable behaviour of these models is often to first study the corresponding monotone model, in which sites can only change state in one direction (see, for example,~\cite{CMnuc2,FSS,MGlauber,HMT1,HMT2,MMT}). 

In this paper we consider general monotone cellular automata in $d$ dimensions. The study of such processes originated in the work of Chalupa, Leath and Reich~\cite{CLR} in 1979, who introduced the family of models now known as `$r$-neighbour bootstrap percolation' (see below), and many other special cases have since been studied~\cite{DE,DH,vEF,GG93,HLR,Mount}. The following definition, however, is due to Bollob\'as, Smith and Uzzell~\cite{BSU}, who were the first to study completely general models in two dimensions.\smallskip

\begin{definition}\label{def:MCA}
Let\/ 
$\U = \{ X_1,\ldots,X_m \}$ be a 
finite collection of finite, non-empty subsets of $\Z^d \setminus \{\0\}$. Let $A \subset \Z_n^d$ be a set of initially \emph{infected} sites, set $A_0 = A$, and define
$$A_{t+1} = A_t \cup \big\{ x \in \Z_n^d \,:\, x + X \subset A_t \text{ for some } X \in \U \big\}$$
for each $t \ge 0$. The \emph{$\U$-closure} of $A$ is the set $[A]_\U := \bigcup_{t\geq 0} A_t$ of eventually infected sites.
\end{definition}

We call this model \emph{\,$\U$-bootstrap percolation}, and we are interested in understanding the behaviour of this process when the initial set $A$ is chosen randomly. To be precise, let us say that $A$ is a \emph{$p$-random} subset of $\Z_n^d$ if each vertex of $\Z_n^d$ is included in $A$ independently with probability $p$, and write $\P_p$ for the associated product probability measure. Our main aim is to determine the behaviour (as $n \to \infty$) of the \emph{critical probability} 
\begin{equation}\label{def:pc}
p_c\big( \Z_n^d,\U \big) := \inf \Big\{ p \in (0,1] \,:\, \P_p\big( [A]_\U = \Z_n^d \big) \ge 1/2 \Big\}
\end{equation}
for each $d$-dimensional \emph{update family} $\U$. 

The most extensively-studied of these models is the classical $r$-neighbour bootstrap process, in which vertices are infected at time $t$ if they have at least $r$ already-infected neighbours. The update family $\cN_r^d$ of this model consists of the ${2d \choose r}$ subsets of size $r$ of the $2d$ nearest neighbours of the origin.  Over the past 40 years a great deal of work has been done on this model,  
and its behaviour is now well-understood, both in two dimensions~\cite{AL,HM,Hol} and in higher dimensions~\cite{BBM3d,BBDM,CC,CM,Sch1}. In particular, 
it was proved by Balogh, Bollob\'as, Duminil-Copin and Morris~\cite{BBM3d,BBDM}, that 
\begin{equation}\label{eq:r:neighbour:threshold}
p_c\big( \Z_n^d, \cN_r^d \big) \, = \, \bigg( \frac{\lambda(d,r) + o(1)}{\log_{(r-1)} n} \bigg)^{d-r+1}
\end{equation}
as $n \to \infty$ for each $2 \le r \le d$, where $\lambda(d,r) > 0$ is an explicit constant. If $r > d$, on the other hand, then it is not hard to see that $p_c( \Z_n^d, \cN_r^d ) \to 1$ as $n \to \infty$, since an uninfected copy of $\{0,1\}^d$ cannot be infected in the process. Similar results have been proved for some other families of automata in two dimensions~\cite{BDMSDuarte,DE,DH,DEH}, and weaker bounds are known for a small number of specific three-dimensional update families~\cite{Blanq19,vEF}.

For general $\U$-bootstrap processes the challenges are rather different, and the first bounds for arbitrary update families\footnote{Slightly earlier, Duminil-Copin and Holroyd~\cite{DH} had obtained very precise results for a large family of `balanced' and `symmetric' two-dimensional monotone cellular automata.} were obtained by Bollob\'as, Smith and Uzzell~\cite{BSU}, who proposed a partition of the family of all two-dimensional models into three classes, which they called `subcritical', `critical' and `supercritical'. They proved that $p_c(\Z_n^2,\U)$ is polynomial in $n$ for supercritical families, and poly-logarithmic in $n$ for critical families. Balister, Bollob\'as, Przykucki and Smith~\cite{BBPS} then proved that the critical probability for subcritical two-dimensional families is bounded away from zero, confirming a conjecture from~\cite{BSU}. For critical two-dimensional families, much more precise bounds were obtained by Bollob\'as, Duminil-Copin, Morris and Smith~\cite{BDMS}, who determined the critical probability up to a constant factor. 

One of the key insights of Bollob\'as, Smith and Uzzell~\cite{BSU} was that, at least in two dimensions, the typical behaviour of the $\U$-bootstrap process 
should, roughly speaking, be determined by its action on discrete half-spaces; specifically, whether or not discrete half-spaces are `stable' under the action of $\U$. Given a $d$-dimensional update family $\U$, we define the set of \emph{stable directions} of $\,\U$ to be
\begin{equation}\label{def:S}
\S(\U) := \big\{ u \in \SS^{d-1} :\, [ \H^d_u ]_\U = \H^d_u \big\},
\end{equation}
where, for each $u \in \SS^{d-1}$, we write $\H^d_u := \{ x \in \Z^d : \< x,u \> < 0 \}$ for the discrete half-space with normal $u$. Note that $u \in \S(\U)$ if and only if $X \not\subset \H^d_u$ for every $X \in \U$. 

The following classification of $d$-dimensional update families was proposed by Bollob\'as, Smith and Uzzell~\cite{BSU} in the case $d = 2$, by Balister, Bollob\'as, Przykucki and Smith~\cite{BBPS} for subcritical families, and by Bollob\'as, Duminil-Copin, Morris and Smith~\cite{BDMS} in general. Given a sphere $\SS \subset \R^d$ of arbitrary dimension and a set $\T \subset \R^d$, we write $\interior_\SS(\T)$ for the interior of $\T$ in $\SS$ with respect to the topology induced by geodesic distance. 

\begin{definition}\label{def:tri}
A $d$-dimensional update family $\U$ with stable set $\S = \S(\U)$ is:
\begin{itemize}
\item \emph{supercritical} if $H \cap \S = \emptyset$ for some open hemisphere $H \subset \SS^{d-1}$;\smallskip
\item \emph{critical} if there exists a hemisphere $H \subset \SS^{d-1}$ such that $\interior_{\SS^{d-1}}(H \cap \S) = \emptyset$ and if $H \cap \S \neq \emptyset$ for every open hemisphere $H \subset \SS^{d-1}$;\smallskip
\item \emph{subcritical} if $\interior_{\SS^{d-1}}(H \cap \S) \neq \emptyset$ for every hemisphere $H \subset \SS^{d-1}$.
\end{itemize}
\end{definition}

It was conjectured by Balister, Bollob\'as, Przykucki and Smith~\cite{BBPS}, and proved by the authors in~\cite{BBMSsub,BBMSupper}, that $p_c(\Z^d,\U) > 0$ if and only if\/ $\U$ is subcritical.\footnote{$p_c(\Z^d,\U)$ is defined as in~\eqref{def:pc}, except the $\U$-bootstrap process is run on the infinite graph $\Z^d$ instead of the finite torus $\Z_n^d$. An alternative proof of the fact that $p_c(\Z^d,\U) > 0$ for all subcritical update families $\U$ was found recently by Hartarsky and Szabó~\cite{HS22}.} 
Moreover, in~\cite{BBMSupper} we proved upper bounds on $p_c(\Z_n^d,\U)$ for all critical and supercritical models. In this paper we prove matching lower bounds for all critical models in $d$ dimensions, and hence complete the proof of the following `universality' theorem, which confirms a conjecture of Bollob\'as, Duminil-Copin, Morris and Smith~\cite{BDMS}, and implies that \emph{every}\/ $\U$-bootstrap process `resembles' (in its large-scale behaviour) one of the $r$-neighbour models.

\begin{theorem}\label{conj:universality}
Let\/ $\U$ be a $d$-dimensional update family.
\begin{itemize}
\item[$(a)$] If\/ $\U$ is supercritical then $p_c\big( \Z_n^d,\U \big) = n^{-\Theta(1)}$.\smallskip
\item[$(b)$] If\/ $\U$ is critical then there exists $r \in \{2,\dots,d\}$ such that
$$p_c\big( \Z_n^d,\U \big) = \bigg(\frac{1}{\log_{(r-1)} n}\bigg)^{\Theta(1)}.$$
\item[$(c)$] If\/ $\U$ is subcritical then $p_c\big( \Z^d,\U \big) > 0$.
\end{itemize}
\end{theorem}

We shall define the quantity $r = r(\U)$ explicitly in Section~\ref{sec:r}, and use it to state a more precise version of the theorem above (see Theorem~\ref{thm:universality}). We shall then be able to state the main result of this paper (Theorem~\ref{thm:lower}), which implies the lower bound in Theorem~\ref{conj:universality} for critical update families. It moreover does so with a universal exponent, depending only on $d$, which is sharp for the $2$-neighbour bootstrap process. It is proved in~\cite{BBMSuncomp}, however, that for critical models with $2 \le r < d$ and for supercritical models with $d \ge 2$, the exponent in the critical probability (if it exists) is in general uncomputable.\footnote{More precisely, it is shown in~\cite{BBMSuncomp} that for each Turing machine program $P$ and each $1 \leq r < d$, there exists a $d$-dimensional update family $\U$ with $r(\U) = r$, such that if $P$ halts then the exponent of $\U$ is at most $2/3$ and if $P$ does not halt then the exponent is at least $1$. If the exponents were always computable, then this would violate the undecidability of the halting problem. The exponent in the case $d = r = 2$ (which was determined in~\cite{BDMS}) \emph{is} computable, but computing it is NP-hard, see~\cite{HMez}.}

All previous lower bounds for critical models with $r \ge 3$ (see ~\cite{BBM3d,BBDM,Blanq22,CC,CM,vEF,HMod}) have used the method introduced by Cerf and Cirillo~\cite{CC} in their groundbreaking work on the $3$-neighbour bootstrap process. In their proof, they coupled this model with a sequence of adjacent two-dimensional models, using the fact that the most `difficult' directions for growth come in pairs of the form $\{u,-u\}$. One of the main problems that we face is that our models do not (in general) have any form of symmetry, and we have therefore had to develop a new and completely different method for controlling the growth. This new method involves a delicate induction on the `size' of an `iceberg' (see Sections~\ref{sec:outline} and~\ref{one:step:sec}--\ref{sec:ind}), and was inspired by Holroyd's method of hierarchies~\cite{Hol}. 

In Sections~\ref{sec:r} and~\ref{sec:outline} we shall state our main theorem precisely, and give a detailed sketch of its proof. The proof itself is given in Sections~\ref{sec:trees}--\ref{sec:ind}. Before embarking on this journey, however, let us provide some further motivation for Theorem~\ref{conj:universality} by briefly discussing some potential applications of our method to two (non-monotone) spin models that have been extensively studied in the statistical physics literature: the Ising model of ferromagnetism, and kinetically constrained models of the liquid-glass transition. 


\subsection{Kinetically constrained models}\label{KCM:sec}

A glass is a disordered material that nevertheless behaves mechanically like a solid, and is formed by rapidly cooling a viscous liquid. Understanding this liquid-glass transition is an important open problem in condensed matter physics (see for example~\cite{ALBB,DS01}). Kinetically constrained models were introduced in the 1980s (see~\cite{FA}, or the reviews~\cite{GST,RS03}) in order to model the liquid-glass transition, and exhibit several key properties of cooled liquids near the glass transition point.

Perhaps the simplest way to understand a kinetically constrained model is as a biased random walk on the family of percolating sets\footnote{More precisely, this is true if the initial set of empty sites percolates for the $\U$-bootstrap process. If the empty sites at time zero are chosen to be $p$-random, then this holds almost surely for every $p > 0$ if $\U$ is not subcritical, by Theorem~\ref{conj:universality} (specifically, by the upper bounds proved in~\cite{BBMSupper}).} in $\U$-bootstrap percolation on $\Z^d$. More specifically, the state (either `empty' or `occupied') of each site $x \in \Z^d$ updates at rate $1$ as long as the set $x + X$ is entirely empty for some $X \in \U$; when it updates, it becomes empty with probability $p$, and occupied with probability $1-p$, independently of its current state (and all other events). Observe that if the initial set $A_0$ of empty sites is chosen to be $p$-random, then this process $(A_t)_{t \ge 0}$ is stationary. Well-studied examples of kinetically constrained models include the East model (see, e.g.,~\cite{CFM}), whose update family 
consists of the single set $\{-1\}$, and the $r$-facilitated Friedrickson--Andersen model, introduced in~\cite{FA}, which corresponds to $r$-neighbour bootstrap percolation.

Over the past few years, there have been some dramatic advances in our understanding of general kinetically constrained models in two dimensions, mirroring those in the study of\, $\U$-bootstrap percolation. One of the main quantities of interest for a kinetically constrained model (starting from equilibrium) is the infection time of the origin\footnote{The mean of this quantity also provides a lower bound on the relaxation time; see~\cite[Section~2.2]{MMT}.}
$$\tau(\Z^d,\U) \, := \, \inf\big\{ t \ge 0 : \0 \in A_t \big\},$$
which, by the fundamental work of Cancrini, Martinelli, Roberto and Toninelli~\cite{CMRT}, is almost surely finite whenever $p > p_c(\Z^d,\U)$. Moreover, a lower bound on the mean infection time is given (up to a constant factor) by the median infection time of the origin in the corresponding $\U$-bootstrap process (see~\cite[Lemma~4.3]{MT}). This bound is sometimes sharp (for example, for the $2$-neighbour model, see~\cite{HMT3}), but for many other update families (for example, for the Duarte model, see~\cite{MaMaT}) the more complex behaviour of kinetically constrained models (in particular, the existence of `energy barriers') causes $\tau(\Z^d,\U)$ to be significantly larger than the bound given by the $\U$-bootstrap process. 

Motivated by the connection with $\U$-bootstrap percolation, a number of conjectures were made in~\cite{M17} regarding the rate of growth of the mean infection time as $p \to 0$ for critical update families in two dimensions, and for supercritical update families in $d$ dimensions. These have now all been either proved or disproved (see~\cite{HMT1,HMT2,MaMaT,MMT}), and in two dimensions the situation is extremely well-understood, with the full universality picture determined (see~\cite{Har,HarMa,HMT3} for the most recent developments). More precisely, both the critical and supercritical families need to be partitioned into two different universality classes~\cite{HMT1,HMT2}, and then further refined for logarithmic corrections, see~\cite{Har,HarMa}. 

In higher dimensions, corresponding results are known only for the $r$-neighbour model, the main obstruction to studying more general models being the lack of tools for controlling the corresponding $\U$-bootstrap models. Having developed such tools in this paper and in~\cite{BBMSupper}, it now seems feasible that substantial progress could be made on this problem. For example, the following lower bound is a straightforward consequence of the results of this paper. The statement depends on the parameter $r(\U)$, which will be defined explicitly in Section~\ref{sec:r} (see Definition~\ref{def:r}). We write $\exp_{(r)}$ for an $r$-times iterated exponential, so $\exp_{(0)} n = n$ and $\exp_{(r)}n = \exp\big( \exp_{(r-1)} n \big)$ for each~$r \ge 1$.

\begin{theorem}\label{thm:KCM}
Let\/ $\U$ be a critical $d$-dimensional update family. Then
$$\Ex\big[ \tau(\Z^d,\U) \big] \ge \exp_{(r-1)}\left( p^{-c} \right),$$
for some constant $c = c(d) > 0$, where $r = r(\U) \in \{2,\ldots,d\}$.
\end{theorem}

The deduction of Theorem~\ref{thm:KCM} from our results is standard (see, e.g.,~\cite{MT}), and we therefore leave the details to the interested reader. It seems reasonable to conjecture that the bound in Theorem~\ref{thm:KCM} is sharp up to the value of the constant $c$ for all critical families (even though this is not true for supercritical families, see~\cite{MaMaT}), and it may perhaps even be possible to prove a matching upper bound using the techniques of~\cite{CMRT,MT,MMT}, together with those introduced in~\cite{BBMSupper}. However, there are likely to be very significant technical challenges involved in making such an approach rigorous.

\subsection{The Ising model of ferromagnetism}\label{Ising:sec}

Introduced over 100 years ago, in an attempt to model Curie's temperature,\footnote{This is the temperature above which a magnetic substance loses its ferromagnetic properties.} the Ising model has become one of the most extensively-studied models in statistical physics. The equilibrium (Gibbs) measures of the model are now relatively well-understood, and a wide array of powerful tools have been developed for their study; we refer the reader to the excellent surveys by Duminil-Copin~\cite{Hugo16,Hugo20}, and the recent papers~\cite{AD,ADS,ADTW,CDHKS,CS,DGR,S10}, for more details. 

The kinetic (or stochastic) Ising model, on the other hand, which was one of the original motivations for the introduction of the $r$-neighbour bootstrap process in~\cite{CLR}, and also the early rigorous work in~\cite{AL}, is rather less well-understood, despite the recent breakthroughs~\cite{LS13,LS16} on the mixing time in the high temperature regime. For example, consider the case of the zero-temperature kinetic Ising model on $\Z^d$ with zero external field. In this setting, the state ($+$ or $-$) of each site of $\Z^d$ updates at random times to agree with the majority of its (nearest) neighbours, with ties broken randomly. It is a longstanding (possibly folklore) conjecture that, for all $d \ge 2$, if the initial states are chosen independently with probability $p > 1/2$ of being $+$, then the system almost surely fixates at $+$. However, this has only been proved when $p$ is sufficiently close to~$1$, by Fontes, Schonmann and Sidoravicius~\cite{FSS}, and when $p \to 1/2$ (slowly) as $d \to \infty$, by Morris~\cite{MGlauber}. We remark that results on the $r$-neighbour bootstrap process (from~\cite{AL} and~\cite{BBMmaj}, respectively) play a crucial role in the proofs in~\cite{FSS} and~\cite{MGlauber}. 

The $r$-neighbour bootstrap process has also found applications in the study of the metastable behaviour of the stochastic Ising model at small positive temperatures and with non-zero external field. 
For example, the relaxation of the system, starting with all states equal, was studied by Schonmann~\cite{Sch4} (who confirmed a conjecture of Aizenman and Lebowitz~\cite{AL})  when the temperature is small but fixed and the external field tends to zero, and by Dehghanpour and Schonmann~\cite{DS1} and Cerf and Manzo~\cite{CMnuc2} when the external field is fixed and the temperature tends to zero. 


We expect that the techniques introduced in this paper will similarly find applications in the study of the stochastic Ising model with more general (finite-range) interactions. For example, it was conjectured in~\cite{M17} that for any critical $d$-dimensional update family $\U$, fixation occurs almost surely in the $\U$-Ising dynamics\footnote{In these dynamics, which were introduced in~\cite{M17}, the state of a site $x$ updates to $+$ at rate $1$ if there exists $X \in \U$ such that $x + X$ is entirely $+$, and similarly for $-$. When $\U = \cN_d^d$, this is equivalent to the zero-temperature nearest-neighbour Ising model, and therefore in this case the conjecture follows from the theorem of Fontes, Schonmann and Sidoravicius~\cite{FSS}.} 
if $p$ is sufficiently close to $1$. There are two main obstructions to using the method of~\cite{FSS} to prove the conjecture: the first, covering the sites in state $-$ with well-separated $\U$-closed droplets, is dealt with by the method of this paper; the second, showing that a $(-)$-droplet disappears in polynomial time, seems to require additional new ideas (see~\cite{Blanq20} for more details). Similarly,~it~seems plausible that the results of~\cite{CMnuc2} could be extended to models with finite-range interactions, and we would expect the techniques introduced in this paper to play an important role in the proof of any such generalization.

\section{The resistance of an update family}\label{sec:r}

In this section we shall define the parameter $r = r(\U)$ of a $d$-dimensional update family (see Definition~\ref{def:r}), which we call the \emph{resistance} of\, $\U$. This will then allow us to state a refined version of Theorem~\ref{conj:universality} (see Theorem~\ref{thm:universality}). To finish the section, we will describe some of the key new ideas and techniques that will be introduced in this paper. 

The idea behind the definition of the resistance $r(\U)$ is quite simple: we will define the `induced' resistance of a point $u \in \SS^{d-1}$ (with respect to the stable set $\S(\U)$) inductively, by considering the resistance of a small copy of $\SS^{d-2}$ centred at $u$. The resistance of $\S(\U)$ will then be defined using the hemisphere in which the maximum resistance is smallest. The details are somewhat technical, however, so the reader may find it helpful to have in mind pictures of the stable sets of the 2- and 3-neighbour models in three dimensions. For $\cN_2^3$, the stable set is the six points $\{\pm e_1, \pm e_2, \pm e_3\}$, and for $\cN_3^3$ the stable set is the union of the three great circles orthogonal to the standard basis vectors. 


First we introduce a family of objects called `$\SS$-stable sets', which are subsets of the sphere $\SS$, and generalize the notion of a stable set of an update family. 

\begin{definition}\label{def:happy}
Let $\T\subset\SS^{d-1}$ and let $\SS$ be a sphere (of arbitrary size and dimension) embedded in $\SS^{d-1}$.
We say that $\T$ is \emph{$\SS$-stable} if there exists a finite collection $\cH_1,\ldots,\cH_m$ of finite families of closed hemispheres of $\SS$, such that
\begin{equation}\label{eq:happy}
\T \cap \SS = \bigcap_{i=1}^m \bigcup_{H \in \cH_i} H.
\end{equation}
\end{definition}


In particular, the stable set of a $d$-dimensional update family is $\SS^{d-1}$-stable.

\begin{lemma}\label{lem:hemispheres}
If\/ $\U$ is a $d$-dimensional update family, 
then
$$\S(\U) = \bigcap_{X \in \, \U} \bigcup_{x \in X} \big\{ u\in \SS^{d-1} : \< x,u \> \geq 0 \big\}.$$
In particular, $\S(\U)$ is $\SS^{d-1}$-stable.
\end{lemma}

\begin{proof}
This is nothing more than the observation that a direction $u \in \SS^{d-1}$ is unstable for $\U$ if and only if there exists $X \in \U$ such that $\< x, u \> < 0$ for every $x \in X$.
\end{proof}

Now, given a sphere $\SS \subset \SS^{d-1}$, $u \in \SS$ and $\eta > 0$, define the sub-sphere\footnote{Here, and throughout the paper, we write $\| \cdot \|$ for the Euclidean norm on $\R^d$.}
\[
S_\eta(\SS,u) := \big\{ v\in \SS : \|u - v\| = \eta \big\}.
\]
We say that two sets $A,B \subset \R^d$ are \emph{equivalent}, and write $A\equiv B$, if one can be obtained from the other by 
a composition of translations, dilations and rotations. The following simple lemma, which was proved in~\cite{BBMSupper}, says that the intersections of a stable set with any two sufficiently small spheres centred at a given point are equivalent. 

\begin{lemma}[Lemma~2.3 of~\cite{BBMSupper}]\label{lem:inducedwelldef}
Let $\SS \subset \SS^{d-1}$ be a sphere, let $\T \subset \SS^{d-1}$ be $\SS$-stable, and let $u \in \SS$. Then there exists $\eta_0 = \eta_0(\SS,\T,u) > 0$ such that, for all $0 < \eta < \eta_0$, 
$$\T \cap S_\eta(\SS,u) \equiv \T \cap S_{\eta_0}(\SS,u)$$
and $\T$ is $S_\eta(\SS,u)$-stable.
\end{lemma}

We shall often use Lemma~\ref{lem:inducedwelldef} implicitly by writing $S_\eta(\SS,u)$ without specifying $\eta$. This should always be taken to mean that $\eta < \eta_0$, where $\eta_0 = \eta_0(\SS,\T,u)$ (for the $\T$ currently under consideration) is given by Lemma~\ref{lem:inducedwelldef}.  

We are now ready to define the resistance of a $d$-dimensional update family. For each $k \in \{0,\ldots,d - 1\}$, let us write $\C^k$ for the set of all $k$-dimensional spheres embedded in $\SS^{d-1}$, and for each $\SS \in \C^k$, let us write $\B(\SS)$ for the set of all $\SS$-stable sets. 

\begin{definition}\label{def:r}
For each $1 \leq k \leq d$ and each $\SS \in \C^{k-1}$, we define two functions
\begin{align*}
\rho^{k-1}(\SS; \, \cdot, \, \cdot) &\colon \B(\SS) \times \SS \to \{0,1,\dots,k\}\\
r^k(\SS; \, \cdot) &\colon \B(\SS) \to \{1,\dots,k+1\}
\end{align*}
inductively as follows. Set $r^0 \equiv 1$, and let $1 \leq k \leq d$, $\SS \in \C^{k-1}$, $\T \in \B(\SS)$ and $u \in \SS$. We define the \emph{induced resistance of $u$ with respect to $\T$ in $\SS$} to be
\begin{equation}\label{eq:rho}
\rho^{k-1}(\SS; \T, u) := 
\begin{cases} 
r^{k-1}\big( S_\eta(\SS,u); \T \big) \quad & \text{if } u \in \T, \\ 
0 & \text{otherwise}, 
\end{cases}
\end{equation}
and the \emph{resistance of $\T$ in $\SS$} to be
\begin{equation}\label{eq:r}
r^k(\SS; \T) := \min_H \, \max _{u \in H} \, \rho^{k-1}\big( \SS; \T, u \big) + 1,
\end{equation}
where the minimum is taken over all open hemispheres $H \subset \SS$.

Now, given a $d$-dimensional update family $\U$, define the \emph{resistance} of $\U$ to be
$$r(\U) := r^d\big( \SS^{d-1}; \S(\U) \big).$$
\end{definition}

In words, the resistance of $\U$ is determined by the largest induced resistance of a direction in the `easiest' hemisphere of $\SS^{d-1}$, and the induced resistance of a stable direction $u$ is given by the resistance of a small sphere centred at $u$. It was shown in~\cite[Observation~2.6]{BBMSupper} that, for each $1 \le r \le d+1$, the update family $\cN_r^d$ has resistance~$r$. 

\medskip
\pagebreak

The following lemma, which was proved in~\cite{BBMSupper}, shows that $r(\U)$ refines Definition~\ref{def:tri}. 

\begin{lemma}[Lemma~2.7 of~\cite{BBMSupper}]\label{lem:tri}
If\/ $\U$ is a $d$-dimensional update family, then:
\begin{itemize}
\item[$(a)$] $\U$ is supercritical if and only if\/ $r(\U) = 1$.\smallskip
\item[$(b)$] $\U$ is critical if and only if\/ $r(\U) \in \{2,\dots,d\}$.\smallskip
\item[$(c)$] $\U$ is subcritical if and only if\/ $r(\U) = d + 1$.\smallskip
\end{itemize}
\end{lemma}

Lemma~\ref{lem:tri} follows in a straightforward way from Definition~\ref{def:r} (see~\cite[Section~2]{BBMSupper} for the details). Having defined $r(\U)$, we can now state a more precise universality theorem. Note in particular that, by Lemma~\ref{lem:tri}, it implies Theorem~\ref{conj:universality}. 

\begin{theorem}\label{thm:universality}
Let $\U$ be a $d$-dimensional update family.
\begin{itemize}
\item[$(a)$] If\/ $\U$ is supercritical then $p_c(\Z_n^d,\U) = n^{-\Theta(1)}$.\smallskip
\item[$(b)$] If\/ $\U$ is critical then 
$$p_c(\Z_n^d,\U) = \bigg(\frac{1}{\log_{(r-1)} n}\bigg)^{\Theta(1)},$$
where $r = r(\U)$.\smallskip
\item[$(c)$] If\/ $\U$ is subcritical then $p_c(\Z^d,\U) > 0$.
\end{itemize}
\end{theorem}

The upper bounds in $(a)$ and $(b)$ were proved in~\cite{BBMSupper}, and part $(c)$ was proved in~\cite{BBMSsub,HS22}. In the next section we shall state the precise lower bound that will be proved in this paper, and provide a detailed outline of our approach. First, however, let us briefly describe a few of the key new ideas and techniques that will be introduced during the proof.

To begin, observe that, by~\eqref{eq:r}, there exists a finite set 
\begin{equation}\label{eq:choosing:T}
\T \subset \big\{ u \in \SS^{d-1} :\, \rho^{d-1}( \SS^{d-1}; \S, u) \ge r - 1 \big\}
\end{equation}
intersecting every open hemisphere of $\SS^{d-1}$. Our main task (see Proposition~\ref{prop:critdrop}) will be to bound the probability that a `$\T$-droplet' $D$ (see Definition~\ref{def:droplet}) is `internally spanned' (see Definition~\ref{def:intspan}). Roughly speaking, this is the event that the $\U$-closure of $D \cap A$ (where $A$ is the set of initially-infected vertices) `connects' the faces of~$D$. Intuitively, the most likely way in which this can occur is via the growth of a droplet inside $D$. It will therefore be necessary to control the growth of the infected set on the faces of this smaller droplet (and on the faces of droplets growing on these faces, and so on). In order to do so, we first construct (in Section~\ref{sec:trees}) a tree of stable directions: for each $u \in \T$, there exists a set of directions (cf.~\eqref{eq:choosing:T}), each with induced resistance at least $r - 2$, that intersects every open hemisphere of $S_\eta(\SS^{d-1},u)$, and so on. There is, however, a significant subtlety in doing this, which we discuss in Section~\ref{trees:sketch:sec}.

Our main tools for controlling growth on the faces of droplets (that are themselves growing on the faces of droplets) are `icebergs', which we introduce in this paper, and which are analogues of droplets for growth that is assisted by a collection of half-spaces. These objects play a central role in our proof (see, for example, our main induction hypothesis, Definition~\ref{def:ih}), and in order to work with them we will need to develop an entirely new set of tools (see Sections~\ref{sec:icebergs} and~\ref{sec:subadd}). Developing these tools required a number of significant innovations; we mention in particular the `iceberg containers' introduced in Section~\ref{sec:icebergs}, and the sub-additivity lemma for iceberg containers proved in Section~\ref{sec:subadd}. These will be described in more detail in Section~\ref{iceberg:chat}. 

We would like to emphasize that, although our icebergs are named after their (much simpler) two-dimensional counterparts from~\cite{BDMS}, the difficulties that we encounter in dealing with the $d$-dimensional versions are completely different, and have no analogue in two dimensions. In particular, it turns out that the usual sub-additivity lemma for droplets in two dimensions \emph{is not true} in higher dimensions, and recovering this situation (which, at first sight, appears to be catastrophic) requires a series of substantial conceptual innovations. To be more precise, we actually prove a sub-additivity lemma for iceberg containers using a subtle and somewhat surprising notion of size (see Definition~\ref{def:diamstar}). The definition of iceberg containers (see Definition~\ref{def:iceberg-container}) is moreover delicate and (necessarily) somewhat complicated, and was designed very carefully in order to allow us to prove this sub-additivity lemma. 

In order to bound the probability that a droplet is internally spanned, we need to bound the probability of all possible ways in which growth can occur, not only the simplest. As mentioned above, in all previous work on models with $r \ge 3$ this was done using the Cerf--Cirillo method, which is not available to us due to the lack of symmetry. We have therefore developed, in Sections~\ref{one:step:sec} and~\ref{sec:ind}, a new method, which was inspired by Holroyd's method of hierarchies~\cite{Hol} (see Lemmas~\ref{lem:goodsat} and~\ref{lem:sideways-iceberg}), but is used here in a very different way, and for a very different purpose. To be more precise, Holroyd used his method to group the possible paths of a growing droplet in the $2$-neighbour model into a bounded collection of `types', which then allowed him to use a simple union bound. Our hierarchies, on the other hand, are used to handle droplets at \emph{very} different scales, and to deal with the lack of symmetry. This is also the first time that hierarchies have been used in a setting where $r(\U) \ge 3$.\footnote{While hierarchies were used in~\cite{BBM3d,BBDM} to study the $r$-neighbour model, the authors of those papers first used the symmetry of the model to reduce to a setting whose behaviour resembled that of the two-neighbour model, and then used hierarchies in a similar way to Holroyd's original paper.} We remark that we also do not use a union bound over hierarchies, but instead a more subtle inductive argument (see Lemma~\ref{lem:indstep}).

\section{The main theorem, and an outline of the proof}\label{sec:outline}

Our task in the remainder of the paper is to prove the following theorem, which implies the lower bound in Theorem~\ref{thm:universality} for critical update families.

\begin{theorem}\label{thm:lower}
Let\/ $\U$ be a critical $d$-dimensional update family. Then
$$p_c\big( \Z_n^d,\U \big) \geq \bigg( \frac{1}{\log_{(r-1)} n} \bigg)^{d - 1 + o(1)},$$
where $r = r(\U)$ is the resistance of\/ $\U$. 
\end{theorem}

The constant $d - 1$ in the exponent is sharp for the $2$-neighbour model, for which the threshold (see~\eqref{eq:r:neighbour:threshold}) was first determined by Aizenman and Lebowitz~\cite{AL}. Recall also that, by the results of~\cite{BBMSuncomp}, the best possible exponent is uncomputable in general if $r < d$. Determining the correct exponent for all $d$-dimensional update families with $r(\U) = d$ also seems likely to be extremely challenging, and would require significant new ideas. 

The rest of this section will be taken up by a detailed outline of the proof of Theorem~\ref{thm:lower}. In particular, we shall introduce (informally) the main new tools in the proof: bounding trees, icebergs, and one-step hierarchies. First, however, we shall introduce the fundamental notion of an `internally spanned $\T$-droplet', which will then allow us to state the key bound that we shall prove on the probability that a $\T$-droplet is internally spanned (see Proposition~\ref{prop:critdrop}). To avoid repetition, let us fix for the rest of the paper a critical $d$-dimensional update family $\U = \{X_1,\ldots,X_m\}$ with $r(\U) = r$. 


\subsection{Internally spanned droplets} 

We begin by defining droplets, which are just discrete polytopes. Note that 
we do not insist that droplets are finite. 



\begin{definition}\label{def:droplet}
Given a finite set $\T \subset \SS^{d-1}$, a \emph{$\T$-droplet} is a non-empty set of the form
\[
D = \bigcap_{u \in \T} \big\{ x \in \Z^d : \< x,u \> \le a_u \big\},
\]
for some collection $\{ a_u \in \R : u \in \T \}$. 
\end{definition}

Note that if $\T \subset \S(\U)$, then a $\T$-droplet is closed under the $\U$-bootstrap process; 
in particular, this holds if $\T$ is chosen as in~\eqref{eq:choosing:T}. We can therefore use such droplets to bound the growth of the process. In order to relate percolation to droplets, we shall use an idea from~\cite{AL,CC}, which we generalize to our setting as follows.   

Suppose that $[A]_\U = \Z_n^d$, and run the $\U$-bootstrap process with initial set $A$, except infecting only one new site in each step. In each step, the diameter\footnote{The diameter $\diam(\X)$ of a set $\X \subset \R^d$ is the supremum of $\|x - y\|$ over all $x,y \in \X$.} of the largest `component' $K$ of infected sites at most triples (once it is larger than some constant), and so at some point in the process $K$ will have diameter $\Theta(\log n)$, and moreover (if we choose the definition of component correctly) will be infected using only sites of $K \cap A$. In particular, if $D$ is the minimal $\T$-droplet containing $K$, then $K \subset [D \cap A]_\U$.

We shall, roughly speaking, follow this strategy (see Lemma~\ref{lem:AL-for-torus}), using the following notion of connectivity. Let us define $R = R(\U)$, the \emph{range} of $\U$, to be 
\begin{equation}\label{eq:R}
R := 2 \cdot \max_{x \in X \in \,\U} \, \|x\|.
\end{equation}

\begin{definition}\label{def:strconn}
We say that a set $K \subset \R^d$ is \emph{strongly connected} if it is connected in the graph $G$ with vertex set $\R^d$, and edge set 
$$E(G) = \big\{  xy : \|x-y\| \le R \big\}.$$
We also say that two sets $K,K' \subset \R^d$ are \emph{strongly connected to each other} if there exist 
$x \in K$ and $y \in K'$ with $\|x - y\| \le R$. 
\end{definition}

Now, observe that if we run the $\U$-bootstrap process on a set $A$, except infecting only one vertex in each step, and if $K$ is a strongly connected component of infected sites at time $t$, then $K \subset [K \cap A]_\U$. To see this, simply note that if $X \in \U$ and $x \in \Z^d$, then $x$ is strongly connected to each vertex of $x + X$. If $\T$ is a set of stable directions, then it follows that the minimal $\T$-droplet $D$ containing $K$ has the following key property. 

\begin{definition}\label{def:intspan}
A $\T$-droplet $D$ is \emph{internally $\T$-spanned} if there exists a set $K \subset D \cap A$ such that $[K]_\U$ is strongly connected and $D$ is the minimal $\T$-droplet containing $[K]_\U$.
\end{definition}

Note that the event that $D$ is internally $\T$-spanned depends on the set $A$, and also on the update family $\U$ (which we fixed above). 
We can now state our key result about droplets, which (by the argument above, together with a simple union bound over droplets) implies a slightly weaker version of Theorem~\ref{thm:lower}, with a worse exponent. 

\begin{prop}\label{prop:critdrop}
There exist $\eps > 0$ and $\delta > 0$, and a finite set $\T \subset \S(\U)$ intersecting every open hemisphere of\/ $\SS^{d-1}$, such that the following holds. Let $D$ be a $\T$-droplet with
$$\diam(D) \leq \exp_{(r-2)}\big(p^{-\eps}\big).$$
Then $D$ is internally $\T$-spanned with probability at most $\exp\big(- \delta \cdot \diam(D)\big)$.
\end{prop}

We remark that in the proof of Theorem~\ref{thm:lower} we shall actually use Proposition~\ref{prop:critdrop-v2}, which is a slightly stronger and more technical variant of this result.

Before describing how we shall prove Proposition~\ref{prop:critdrop}, let us begin by explaining why the methods used in previous works fail for general models. As mentioned above, all previous lower bounds for models with $r \ge 3$ (see~\cite{BBM3d,BBDM,CC,CM,vEF,HMod}) used a technique that was introduced by Cerf and Cirillo~\cite{CC} during their work on the $r$-neighbour process. Note that in that setting one can take $\T = \{\pm e_1, \dots, \pm e_d\}$, and thus $\T$-droplets are cuboids. To bound the probability that a cuboid $D$ is internally $\T$-spanned, they partitioned $D$ into $(d-1)$-dimensional slabs of thickness 2, coupled the $r$-neighbour process in $D$ with independent $(r-1)$-neighbour processes in these slabs, noting that each vertex has only one neighbour in $D$ outside its own slab, and used induction on $r$ to control the process inside the slabs. Observe that symmetry plays a crucial role in the coupling.

In the definition of $\,\U$ we made no assumptions of symmetry, and one of our main challenges in proving Proposition~\ref{prop:critdrop} will be to deal with the lack of such an assumption. Indeed, the strategy described above of partitioning into lower-dimensional slabs would not work for us, since we cannot expect to have two opposite stable directions, and therefore growth inside slabs can be assisted by $\T$-droplets growing outside those slabs and `breaking in' -- such $\T$-droplets are sometimes called `savers'. 

\subsection{One-step hierarchies}\label{one:step:chat}

In order to deal with the lack of symmetry, we introduce a novel `one-step' variant of the `hierarchies' technique that was pioneered by Holroyd~\cite{Hol} in his work on the two-neighbour model, and which has since become a fundamental tool in the study of bootstrap percolation (see for example~\cite{BBDM,BDMS,DE}). Hierarchies are normally used in order to carefully analyse the possible growth paths of droplets, and thus to prove sharp threshold results for symmetric models with $r(\U) = 2$. An important innovation of the paper~\cite{BDMS} was to apply them for a different purpose: to deal with the lack of symmetry in general two-dimensional models.\footnote{To be precise, the authors of~\cite{BDMS} only used hierarchies for unbalanced models, whose growth is asymptotically one-dimensional.} Here we again use them for this purpose, but the difficulties involved in doing so are very different, and we shall need to modify the method much more dramatically. In particular, in all previous applications of the method of hierarchies the droplets being studied had size $p^{-\Theta(1)}$. Here, on the other hand, we will need to `track' the growth of droplets (and, more generally, icebergs) 
from size $p^{-O(1)}$ all the way up to the `critical' scale, $\exp_{(r-2)}(p^{- \Theta(1)})$. 

{\setstretch{1.14}

In the usual hierarchies method (including in~\cite{BDMS}) one uses a union bound over hierarchies, but when $r(\U) \ge 3$ there would be far too many hierarchies for this to work. In order to overcome this problem we introduce a simpler family of objects, which we call `one-step hierarchies' (see Section~\ref{one:step:sec}). Roughly speaking, these consist either of two icebergs that are sufficiently large and disjointly `iceberg spanned'\footnote{The definition of an iceberg being `iceberg spanned' is introduced in Section~\ref{sec:icebergs} (see Definitions~\ref{def:iceberg-container} and~\ref{def:iceberg spanned}). It is somewhat more complicated than the corresponding notion for droplets.}, or of two icebergs, one inside the other, where the smaller iceberg is iceberg spanned and also a certain `crossing event' holds (see Lemma~\ref{lem:goodsat}). The construction of these objects relies heavily on the theory of icebergs that we develop in Sections~\ref{sec:icebergs} and~\ref{sec:subadd} (see also Section~\ref{iceberg:chat}), but otherwise follows the familiar strategy that was introduced in~\cite{Hol}. 

The second main innovation of our new technique is the way in which we use these one-step hierarchies to bound the probability that an iceberg is spanned. We do so by replacing the simple union bound used in previous applications by a carefully designed induction (see Definition~\ref{def:ih} and Lemma~\ref{lem:indstep}). Intuitively, this allows us to `pay' for the number of choices for the hierarchies by allowing our bound on the probability that an iceberg is spanned to increase gradually with the diameter. More precisely, 
our bound will increase (roughly speaking) from $p^x$ when the diameter $x$ is small, to $e^{-x}$ at the critical scale $x = \exp_{(r-2)}(p^{- \Theta(1)})$, see~\eqref{eq:qx}. We will apply the induction hypothesis to the icebergs in our one-step hierarchy, and also to a smaller iceberg 
whose existence we deduce from the crossing event (see Lemma~\ref{lem:sideways-iceberg}). 

The reader may be wondering why we worked with icebergs in the discussion above, rather than droplets. The reason is that when $r \ge 4$, we shall need to apply the strategy outlined above for droplets growing `on the sides' of larger droplets, and for droplets growing on the sides of those droplets, and so on. In this more complex setting the role of droplets is played by icebergs, and in order to carry out our approach using one-step hierarchies we will need to develop a new general method for using icebergs to bound growth in the $\U$-bootstrap process `assisted' by certain unions of half-spaces. 

}

\subsection{Icebergs}\label{iceberg:chat}

{\setstretch{1.14}

Consider a sequence of droplets, each growing on the sides of the droplets that came before it in the sequence. In order to study the growth in such situations, it will be convenient to replace each droplet in the sequence by a half-space $\H^d_{u_i}$, where $u_i$ is the normal to the side on which we are growing, and consider the process `assisted' by the set $\H^d_{u_1} \cup \cdots \cup \H^d_{u_k}$. We will choose the directions $u_1,\ldots,u_k$ carefully so that this set is closed under the $\U$-bootstrap process (see Proposition~\ref{lem:trees-exist}). In Sections~\ref{sec:icebergs} and~\ref{sec:subadd} we shall develop a theory of `icebergs', which play the role of droplets in this setting. 

As mentioned above, icebergs are named after some similar (but much simpler) two-dimensional objects that were introduced in~\cite{BDMS} in order to control the growth on the sides of droplets in models with `drift' (such as the Duarte model). In higher dimensions icebergs are significantly more complicated objects, and they will play a \emph{very} different (and \emph{much} more fundamental) role in our proof than they did in~\cite{BDMS}. Indeed, they are the central object of study in this paper: all of our results in Sections~\ref{sec:icebergs}--\ref{sec:ind}, including our main induction hypothesis (see Definition~\ref{def:ih}), are about properties of icebergs. 

}

Very roughly speaking, our strategy for controlling icebergs follows the approach that was first introduced in the groundbreaking work of Aizenman and Lebowitz~\cite{AL}, and that has been used in essentially all work on bootstrap percolation since then. In particular, we will need to design an `iceberg spanning algorithm' (see Definition~\ref{def:icespanalg}), and prove analogues for icebergs of the classical `extremal' and `Aizenman--Lebowitz' lemmas (see Lemmas~\ref{lem:AL} and~\ref{lem:extremal}). However, as always, the difficulty of carrying out this approach lies in the details, and here the obstructions to doing so are much greater than usual. 

In this subsection we will concentrate on the central obstruction that we shall need to overcome in the proof of Theorem~\ref{thm:lower}: designing and proving a suitable `sub-additivity' lemma for icebergs. For the $2$-neighbour model such a result is immediate: it simply says that if two (axis-parallel) rectangles $R_1$ and $R_2$ intersect, then the semi-perimeter of the smallest rectangle that contains both is at most the sum of the semi-perimeters of $R_1$ and $R_2$. A similar result also holds for arbitrary polygons in two dimensions, see~\cite[Lemma~6.4]{BDMS} for a simple proof. For polytopes in three or more dimensions, however, the corresponding sub-additivity lemma turns out to be false! 


Overcoming this obstacle is the main topic of Sections~\ref{sec:icebergs} and~\ref{sec:subadd}, and required a number of significant conceptual innovations, as well as the development of an entirely new and substantially more involved set of tools. The result we eventually manage to prove is somewhat complicated (see Lemma~\ref{lem:subadd}): 
it is a sub-additivity lemma for a certain measure (the `iceberg diameter') of the size of the `iceberg container' of the `$u$-closure' of a set of vertices $K$. To unpack this a little, let us take these definitions in reverse order. First, the $u$-closure of $K$ is simply the set of vertices infected in the $\U$-bootstrap process with the help of the set $\H^d_{u_1} \cup \cdots \cup \H^d_{u_k}$, starting from set $K$ (see Definition~\ref{def:u-closure}). Next, the iceberg container of a set (see Definition~\ref{def:iceberg-container}) is an iceberg that is chosen very carefully so that it is not too large, and has several useful properties (in particular, various kinds of monotonicity). Finally, the `iceberg diameter' (see Definition~\ref{def:diamstar}) of an iceberg $(J,u)$ is $\delta^k$ times the `weighted diameter' of a certain (again, carefully chosen) polytope that contains $J$, where $\delta > 0$ is a small constant and $k$ is a function of the `type' $u$ of the iceberg (roughly speaking, it is the number of half-spaces used when constructing $J$). The weighted diameter (see Definition~\ref{def:weighted-diam}), which plays a key role in the proof, measures the size of a polytope by taking a suitable linear sum of the distances of the faces to the origin. Moreover, and crucially, we can control its interaction with the definition of iceberg containers, and also with two particular families of polytopes (see Definition~\ref{def:FuGu} and several of the lemmas in Section~\ref{sec:subadd}, in particular Lemmas~\ref{lem:sub-add-droplets} and~\ref{lem:wd-permissible}). 

We remark that our sub-additivity lemma also plays a key role in the construction of our one-step hierarchies (see Lemma~\ref{lem:goodsat}), and that our induction hypothesis (Definition~\ref{def:ih}) will be defined using the iceberg diameter. We consider the techniques developed in the proof of Lemma~\ref{lem:subadd} to be of significant independent interest, and expect them to have many further applications in the study of high-dimensional cellular automata.

\subsection{Bounding trees}\label{trees:sketch:sec}

An important point, which we have so far avoided discussing, is that of exactly \emph{which} stable directions and half-spaces we use to define our droplets and icebergs. This question will be the subject of Section~\ref{sec:trees}, but let us briefly describe here the basic idea. First, choosing the set $\T$ in Proposition~\ref{prop:critdrop} is straightforward, since the definition of $r(\U)$ implies (see Definition~\ref{def:r}) that for every open hemisphere $H \subset \SS^{d-1}$, there exists a direction $u \in H$ such that 
$$\rho^{d-1}\big( \SS^{d-1}; \S(\U), u \big) \ge r - 1.$$
We can take $\T$ to be any finite collection of these directions that intersects every open hemisphere. Now, for each $u \in \T$, we need to choose a collection of directions that we can use to control growth on the $u$-side of a $\T$-droplet. We again use Definition~\ref{def:r}, which implies that for every open hemisphere $H \subset \SS_u := S_\eta(\SS^{d-1},u)$, we have 
\begin{equation}\label{eq:tree:sketch:u}
\rho^{d-2}\big( \SS_u; \S(\U), v \big) \ge r - 2
\end{equation}
for some direction $v \in H$. We choose a finite collection of these directions that intersects every open hemisphere of $S_\eta(\SS^{d-1},u)$, and use these directions (together with $-u$) to define the icebergs that we will use to control the $\U$-bootstrap process assisted by $\H_u^d$.

Now, however, comes a subtle and important point. If $r \ge 4$, then for each $v$ chosen above, we need to choose a collection of directions that we can use to control growth on the $v$-side of an iceberg that is itself growing on the $u$-side of a $\T$-droplet. Using Definition~\ref{def:r} in a naive way, we would obtain, for each open hemisphere $H \subset S_\eta(\SS_u,v)$, a direction $w \in H$ such that
$$\rho^{d-3}\big( S_\eta(\SS_u,v); \S(\U), w \big) \ge r - 3.$$
However, choosing a collection of such directions that intersects every open hemisphere of $S_\eta(\SS_u,v)$ is not sufficient. The issue is that $S_\eta(\SS_u,v)$ has too small a dimension, which means that droplets formed by these directions together with $-v$ may be infinite, and we need them to be finite (see Lemma~\ref{lem:icebergs-finite}). What we would like to do instead is to choose a direction in each hemisphere of $S_\eta(\SS^{d-1},v)$, but a priori we do not know that there exist directions with suitable resistances here. In order to resolve this issue, we use a property of stable sets (Lemma~\ref{lem:memphis}) that was proved in~\cite{BBMSupper}, and which plays a key role in the proof of the upper bounds in Theorem~\ref{thm:universality}. This lemma implies that if $v \in \SS_u$ satisfies~\eqref{eq:tree:sketch:u}, then
$$\rho^{d-1}\big( \SS^{d-1}; \S(\U), v \big) \ge r - 2,$$
and hence, for every open hemisphere $H \subset \SS_v := S_\eta(\SS^{d-1},v)$, there exists $w \in H$ with
$$\rho^{d-2}\big( \SS_v; \S(\U), w \big) \ge r - 3.$$
We choose a finite collection of these directions that intersects every open hemisphere of $\SS_v$, and use these directions (together with $-u$ and $-v$) to define the icebergs that control growth on the $v$-side of an iceberg that is growing on the $u$-side of a $\T$-droplet. 

Iterating this process gives a tree of stable directions that we use to define our icebergs; however, there are two other subtle points that we must consider when choosing this tree: the directions along each path from the root to a leaf must be linearly independent, and the union of any collection of half-spaces in these directions must be $\U$-closed. These conditions will allow us define icebergs (for this we will need the boundaries of any such collection of half-spaces to intersect), and to show that they are $\U$-closed, respectively. To deal with the requirement that the directions along a path are linearly independent (which only causes problems when $r \ge 5$), we shall prove a slightly technical lemma, which allows us to avoid the subspace spanned by the ancestors of a vertex in the tree when choosing its children (see Lemma~\ref{lem:memphis-subspace}). To prove that the union of half-spaces is closed is simpler: we just need to take $\eta$ sufficiently small (see Proposition~\ref{lem:trees-exist}). 

\subsection{The lattice $\Lambda$}

{\setstretch{1.12}

To finish this section, let us note a slightly subtle technical point: it will be convenient in the proof to treat the origin in $\R^d$ as different from the origin in $\Z^d$, since we shall need to consider shifted (continuous) half-spaces whose boundaries do not necessarily contain a lattice point (see Section~\ref{one:step:sec}). Equivalently, we shall think of the lattice $\Z^d$ as being translated by some unknown element of $\R^d$. The following lemma records the fact that shifting $\Z^d$ in this way does not affect the stable set $\S(\U)$. 

\begin{lemma}\label{lem:shift:Ssame}
Let $u \in \SS^{d-1}$ and $a \in \R$, and define\/ $\H_u^d(a) := \{ x \in \Z^d : \< x,u \> < a \}$. Then 
$$[\H_u^d(a)]_\U = \H_u^d(a) \qquad \Leftrightarrow \qquad [\H_u^d]_\U = \H_u^d \qquad \Leftrightarrow \qquad u \in \S(\U).$$
\end{lemma}

\begin{proof}
The second equivalence is the definition~\eqref{def:S} of $\S(\U)$, so it suffices to prove, for each $X \in \U$, that $X \subset \H_u^d$ if and only if $x + X \subset \H_u^d(a)$ for some $x \in \Z^d \setminus \H_u^d(a)$. One direction is immediate, since if $\< x + y,u\> < a$ and $\<x,u\> \ge a$, then $\< y, u\> < 0$. In the other direction, we are similarly done if there exists $x \in \Z^d$ with $\< x,u\> = a$, so assume not, and suppose that $X \subset \H_u^d$. Recalling that $X$ is finite, set 
$$\gamma := \max\big\{ \<y,u\> : y \in X \big\},$$
and note that $\gamma < 0$. Now, by increasing $a$ if necessary (without changing $\H_u^d(a)$) we may assume that there exists $x \in \Z^d \setminus \H_u^d(a)$ with $\<x,u\> < a - \gamma$. But now $\<x + y,u\> < a$ for every $y \in X$, and hence $x + X \subset \H_u^d(a)$, as required.
\end{proof}

Since the properties stated in Section~\ref{sec:r} only depended on the stable set $\S(\U)$, it follows from Lemma~\ref{lem:shift:Ssame} that they all hold equally well in the `shifted' setting. We shall therefore assume that the\/ $\U$-bootstrap process takes place on a lattice
$$\Lambda := y + \Z^d \subset \R^d$$ 
for some arbitrary $y \in \R^d$, and redefine $\H_u^d$ to be the discrete half-space in this lattice, i.e., $\H_u^d := \big\{ x \in \Lambda : \< x,u \> < 0 \big\}$. 
It will only be important in a few places~that we allow $y \not\in \Z^d$, and these will be pointed out to the reader when they occur. 
}


\section{Bounding trees}\label{sec:trees}

{\setstretch{1.12}

Recall that in Section~\ref{sec:outline} we fixed integers $d \ge r \ge 2$, and also a $d$-dimensional update family $\U = \{X_1,\ldots,X_m\}$ with $r(\U) = r$. The purpose of this section is to select the various finite subsets of $\S = \S(\U)$ whose elements will define the faces of the droplets and icebergs that we shall later use to bound growth under the action of the $\U$-bootstrap process. These sets of stable directions will be constructed in Proposition~\ref{lem:trees-exist} via a directed tree $T$ of depth $r - 1$, rooted at $\0$ (the origin) and with all other vertices elements of $\SS^{d-1}$. This tree will have the following properties: 
\begin{itemize}
\item[$(a)$] $\rho^{d-1}\big( \SS^{d-1}; \S, u \big) \ge r - i$ for every vertex $u \in V(T)$ at distance $i$ from the root;\smallskip
\item[$(b)$] if $\0 \to u_1 \to \cdots \to u_{r-1}$ is a directed path in $T$, then the set $\{u_1,\ldots,u_{r-1}\}$ is linearly independent;\smallskip 
\item[$(c)$] the children of the vertex $u \in V(T)$ (if they exist) all lie in a sufficiently small $(d - 1)$-sphere $S_\eta(\SS^{d-1},u) \subset \SS^{d-1}$ centred at $u$.
\end{itemize}
The vertices of $T$ at depth 1 from the root will form the set $\T$ in Proposition~\ref{prop:critdrop}, and for each path $\0 \to u_1 \to \cdots \to u_k$ in $T$ with $u_k$ not a leaf, the set of children of $u_k$, together with $\{-u_1,\ldots,-u_k\}$, will form the set of directions used to define `$u$-iceberg droplets' (see Definition~\ref{def:iceberg}). Icebergs will be the central subject of Sections~\ref{sec:icebergs} and~\ref{sec:subadd}.
}

To underline the importance of the tree $T$ defined in this section, let us remark that, once $T$ has been chosen, the only other property of $\U$ that will be used in the proof of Theorem~\ref{thm:lower} 
is its range $R$, which was defined in~\eqref{eq:R}.

We divide the (somewhat technical) definition of $T$ into two parts. In the first, we define the abstract concept of a `bounding tree' without any reference to stable sets. In the second, we introduce the connection between stable sets and bounding trees. The reader will see that we in fact define a pair $(T,\eta)$, where $T$ is a tree and $\eta$ plays the same role as the function whose existence (for a given stable set) is guaranteed by Lemma~\ref{lem:inducedwelldef}. Since we wish to define bounding trees without reference to stable sets, we allow an arbitrary function $\eta$ in Definition~\ref{def:tree}. We remark that 
the function $\eta$ will be discarded almost as soon as our final tree $T$ has been constructed (in Proposition~\ref{lem:trees-exist}), since all of the information that we will need in later sections will be encoded in the tree $T$. 

Given a sphere $\SS \subset \SS^{d-1}$ and a set $B \subset \SS$, we say that $B$ is an \emph{$\SS$-bounding set} if $B \cap H \neq \emptyset$ for every open hemisphere $H \subset \SS$. 
We also write $N_T^\to(v)$ to denote the out-neighbourhood of a vertex $v$ in a directed tree $T$. 

\begin{definition}\label{def:tree}
A \emph{bounding tree} is a pair $\B=(T,\eta)$, where $T$ is a finite, directed, rooted tree, with root vertex $\0\in\R^d$, all other vertices elements of $\SS^{d-1}$, and all edges directed away from the root, and the \emph{magnification function} $\eta \colon \{1,\ldots,s-1\} \to (0,\infty)$, for some $1 \le s \le d-1$, is a decreasing function  
such that the following hold:
\begin{itemize}
\item[$(a)$] Every leaf of $T$ is at (graph) distance $s$ from $\0$.\smallskip
\item[$(b)$] If $\0 \to u_1 \to \cdots \to u_s$ is a directed path in $T$, then
\begin{equation}\label{eq:tree-nbours}
u_{i+1} \in \SS_{u_i} := S_{\eta(i)} \big( \SS^{d-1}, u_i \big)
\end{equation}
for each $1 \le i \le s - 1$, and the set $\{u_1,\ldots,u_s\}$ is linearly independent.\smallskip
\item[$(c)$] If $u\in V(T)$ is not a leaf, then $N_T^\to(u)$ is an $\SS_u$-bounding set, where $\SS_\0 := \SS^{d-1}$.
\end{itemize}
The parameter $s$ is the \emph{depth} of $\B$.
\end{definition}

Thus, if $\B = (T,\eta)$ is a bounding tree, then the non-root vertices are elements of $\SS^{d-1}$, and the children of a vertex $u$ of depth $i$ in $T$ all lie in (and intersect every open hemisphere of) a sphere of radius $\eta(i)$ around $u$ in $\SS^{d-1}$. The reader should think of the function $\eta$ as decreasing sufficiently rapidly so that the children of vertices of greater depth lie much closer to their parents than vertices closer to the root.

\begin{definition}\label{def:tree-good}
A bounding tree $\B = (T,\eta)$ of depth $r - 1$ is \emph{$\S$-good} if the following holds. For every directed path $\0 \to u_1 \to \cdots \to u_{r-1}$ in $T$, 
\begin{equation}\label{eq:condition:d}
\rho^{d-1}(\SS^{d-1}; \S, u_i) \ge r - i
\end{equation}
for each $1 \le i \le r - 1$.
\end{definition}

Observe in particular that if $\B$ is $\S$-good then $\rho^{d-1}(\SS^{d-1}; \S, u) \ge 1$ for every vertex $u \in V(T)$, and hence every vertex of the tree $T$ is a stable direction of $\U$.

\subsection{Locally inherited resistance}\label{sec:memphis}

As explained at the end of Section~\ref{sec:outline}, we will use two technical lemmas involving stable sets in order to complete our construction of an $\S$-good bounding tree. The first of those lemmas was proved in~\cite{BBMSupper}. 

\begin{lemma}[Lemma~4.2 of~\cite{BBMSupper}]\label{lem:memphis}
If\/ $\T$ is an\/ $\SS^{d-1}$-stable set, then
\begin{equation}\label{eq:memphis}
\rho^{d-1} \big( \SS^{d-1}; \T, v \big) \ge \rho^{d-2} \big( S_\eta(\SS^{d-1},u); \T, v \big)
\end{equation}
for every $u \in \SS^{d-1}$ and $v \in S_\eta(\SS^{d-1},u)$.
\end{lemma}

The second of the two technical lemmas (Lemma~\ref{lem:memphis-subspace} below) will provide us with a lower bound on the right-hand side of~\eqref{eq:memphis}. For the statement, it will be convenient to work in a copy of the sphere $\SS^{d-2}$, rather than the sphere $\SS_u = S_\eta(\SS^{d-1},u)$. We shall later deduce a bound for $\SS_u$ using the homothety $\varphi_u$ that maps the centre of $\SS_u$ to the origin, and $\SS_u$ to the sphere $\SS^{d-1} \cap \{u\}^\perp$. 

Let us write $\SS^{d-2}$ to denote an arbitrary $(d-2)$-dimensional sphere in $\SS^{d-1}$, centred at the origin. Recall from~Definition~\ref{def:happy} that if $\T$ is an $\SS^{d-2}$-stable set, then
\begin{equation}\label{eq:happy:reminder}
\T \cap \SS^{d-2} = \bigcap_{i=1}^t \bigcup_{u \in Y_i} H_u
\end{equation}
for some finite collection $Y_1,\ldots,Y_t$ of finite families of vectors in $\SS^{d-2}$, where $H_u$ is the closed hemisphere of $\SS^{d-2}$ centred at $u$; that is, 
$$H_u := \big\{ v \in \SS^{d-2} : \< u,v \> \geq 0 \big\}.$$ 
For each $\SS^{d-2}$-stable set $\T$, let us choose a minimal representation as in~\eqref{eq:happy:reminder},\footnote{By `minimal' we mean that no proper subset of the collection of $H_u$ gives the same $\T \cap \SS^{d-2}$.} and define the \emph{set of centres of $\T$ with respect to $\SS^{d-2}$} to be \begin{equation}\label{def:CST}
C(\SS^{d-2};\T) := Y_1 \cup \cdots \cup Y_t.
\end{equation}
In particular, if $\T = \S$ then the set of centres can be taken to be (a subset of) the elements of the update sets (cf.~Lemma~\ref{lem:hemispheres}). Finally, recall that $B$ is an $\SS^{d-2}$-bounding set if $B \cap H \neq \emptyset$ for every open hemisphere $H \subset \SS^{d-2}$.\footnote{In particular, an $\SS^{d-2}$-bounding set can be infinite.}

The following second technical lemma will allow us to construct a bounding tree whose paths are linearly independent. 

\begin{lemma}\label{lem:memphis-subspace}
Let\/ $\T$ be an\/ $\SS^{d-2}$-stable set, and let\/ $W \subset \SS^{d-1}$ be such that $W^\perp \ne \{\0\}$ and $W^\perp \cap \SS^{d-1} \subset \SS^{d-2}$. If\/ $C(\SS^{d-2};\T) \subset W^\perp$, then there exists an $\SS^{d-2}$-bounding set $B$ with
\begin{equation}\label{eq:memphis-subspace-rho}
\rho^{d-2}\big( \SS^{d-2}; \T, v \big) \geq r^{d-1}\big( \SS^{d-2} ;\T \big) - 1
\end{equation}
for every $v \in B$, and such that $B \cap \< W \> = \emptyset$.
\end{lemma}

Lemma~\ref{lem:memphis-subspace} will be used to construct the neighbourhood of a vertex $u$ in our $\S$-good bounding tree. In our application we shall set $\T := \varphi_u( \S \cap \SS_u )$ and $W := \{ u_1,\ldots,u_s\}$, where $\varphi_u$ is the homothety that maps $\SS_u$ to $\SS^{d-2}$, and $\0 \to u_1 \to \cdots \to u_s = u$ is the directed path in $T$ from the root to $u$. We remark that the existence of an $\SS^{d-2}$-bounding set whose elements satisfy~\eqref{eq:memphis-subspace-rho} follows easily from Definition~\ref{def:r}, so the content of the lemma is that we can choose $B$ avoiding the span of $W$.  

In order to prove Lemma~\ref{lem:memphis-subspace}, we need another lemma from~\cite{BBMSupper}, whose statement requires a little extra notation. Given a set $W \subset \SS^{d-1}$, let us write 
$$\SS(W) := \SS^{d-1} \cap W^\perp,$$
and define the \emph{projection} $\pi(u,W^\perp)$ of a point $u \in \SS^{d-1} \setminus \< W \>$ onto $\SS(W)$ to be the unique element of $\SS(W)$ such that 
\begin{equation}\label{def:projection}
u = w + \lambda \cdot \pi(u,W^\perp)
\end{equation}
for some $w \in \< W \>$ and $\lambda > 0$. Note that $\varphi_u(x) = \pi(x,\{u\}^\perp)$, and that if $x \in W^\perp$, then 
\begin{equation}\label{obs:projection}
\sgn\big( \< x,u \> \big) = \sgn\big( \< x,\pi(u,W^\perp) \> \big),
\end{equation}
since $\< x,u \> = \lambda \cdot \< x, \pi(u,W^\perp) \>$. The following is the additional lemma that we need. 

\begin{lemma}[Lemma~4.5 of~\cite{BBMSupper}]\label{lem:memphis-vertical}
Let\/ $\T$ be an\/ $\SS^{d-2}$-stable set, and let $W \subset \SS^{d-1}$. If\/ $C(\SS^{d-2};\T) \subset W^\perp$, then
\begin{equation}\label{eq:memphis-vertical-rho}
\rho^{d-2}\big( \SS^{d-2};\T,u \big) = \rho^{d-2}\big( \SS^{d-2};\T,v \big),
\end{equation}
for every $u,v \in \SS^{d-2} \setminus \< W \>$ such that $\pi(u,W^\perp) = \pi(v,W^\perp)$. 
\end{lemma}

The proof of Lemma~\ref{lem:memphis-vertical} in~\cite{BBMSupper} is not too difficult: the key observation is that in a representation of $\T \cap \SS^{d-2}$ as in~\eqref{eq:happy:reminder}, the boundary of each hemisphere $H_u$ must either pass through both $u$ and $v$, or miss both (and hence miss an open ball around each).

With Lemma~\ref{lem:memphis-vertical} in hand, we are now ready to prove Lemma~\ref{lem:memphis-subspace}.

\begin{proof}[Proof of Lemma~\ref{lem:memphis-subspace}]
We begin by constructing an $\SS(W)$-bounding set $B_W$ such that~\eqref{eq:memphis-subspace-rho} holds for every $v \in B_W$. We shall later show (in~\eqref{eq:BW-to-B}) how to augment $B_W$ to obtain an $\SS^{d-2}$-bounding set $B$, recalling that $\SS(W) \subset \SS^{d-2}$ by assumption.

To construct $B_W$, observe first that, by~\eqref{eq:r}, for each $u \in \SS(W)$, there exists $w_u \in \SS^{d-2}$ such that $\< u,w_u \> > 0$ and 
\begin{equation}\label{eq:memphis-subspace-rho-repeated}
\rho^{d-2}\big( \SS^{d-2}; \T, w_u \big) \ge r^{d-1}\big( \SS^{d-2} ;\T \big) - 1.
\end{equation}
Note that $w_u \not\in \<W\>$, since $u \in \SS(W) \subset W^\perp$. We may therefore set $v_u := \pi(w_u,W^\perp)$, and define 
$$B_W := \big\{ v_u : u \in \SS(W) \big\}.$$
To see that $B_W$ is an $\SS(W)$-bounding set, observe that $\< u,v_u \> > 0$ for every $u \in \SS(W)$, since $\< u,w_u \> > 0$ by construction, and by~\eqref{obs:projection}. Moreover, by Lemma~\ref{lem:memphis-vertical} and~\eqref{eq:memphis-subspace-rho-repeated}, we have 
$$\rho^{d-2}\big( \SS^{d-2}; \T, v_u \big) = \rho^{d-2} \big( \SS^{d-2}; \T, w_u \big) \ge r^{d-1}\big( \SS^{d-2} ;\T \big) - 1,$$
since $v_u,w_u \not\in \<W\>$. It follows that~\eqref{eq:memphis-subspace-rho} holds for every $v \in B_W$, as claimed.

Now let $Y := \big\{ x \in \<W\> : \|x\| = 1/2 \big\}$, and for each $x \in Y$, define 
$$B_x := \big\{ v \in \SS^{d-2} : v - x \in W^\perp \,\text{ and }\, \pi(v,W^\perp) \in B_W \big\}.$$
We claim that the set
\begin{equation}\label{eq:BW-to-B}
B := B_W \cup \bigcup_{x \in Y} B_x
\end{equation}
is an $\SS^{d-2}$-bounding set, that $B \cap \<W\> = \emptyset$, and that~\eqref{eq:memphis-subspace-rho} holds for every $v \in B$. 
 
To show that $B \cap \<W\> = \emptyset$, recall that $B_W \subset W^\perp$ and that $x \in \<W\>$ and $\|x\| = 1/2$ for each $x \in Y$. Since $B_x \subset \SS^{d-2} \cap ( x + W^\perp )$, it follows that $B_x \cap \<W\> = \emptyset$. 

Next, to show that~\eqref{eq:memphis-subspace-rho} holds for $v \in B$, observe that, by Lemma~\ref{lem:memphis-vertical},  
$$\rho^{d-2}(\SS^{d-2};\T,v) = \rho^{d-2}\big( \SS^{d-2};\T,\pi(v,W^\perp) \big),$$
since $v \in B$ and $B \cap \<W\> = \emptyset$. Now, since $\pi(v,W^\perp) \in B_W$ for every $v \in B$, and since~\eqref{eq:memphis-subspace-rho} holds for every element of $B_W$, it follows that~\eqref{eq:memphis-subspace-rho} holds for $v$.  
 
Our only remaining task is therefore to verify that $B$ is an $\SS^{d-2}$-bounding set. To do so, we need to show that for each $u \in \SS^{d-2}$, there exists $v \in B$ such that  $\< u,v \> > 0$. If $u \in \<W\>$ then set $x := u/2$, and note that $x \in Y$. Now, observe that $B_x$ is non-empty, because $B_W$ is an $\SS(W)$-bounding set and $W^\perp \ne \{\0\}$, and so $B_W$ is non-empty. We may therefore choose an arbitrary element $v \in B_x$, and observe that $v = x + y$ for some $y \in W^\perp$, so $\< u,v \> = \< u,x\> = 1/2 > 0$, as required. On the other hand, if $u \not\in \<W\>$, then set $x := \pi(u,W^\perp) \in \SS(W)$. Since $B_W$ is an $\SS(W)$-bounding set, there exists $v \in B_W$ with $\< x,v \> > 0$, and since $v \in \SS(W)$, it follows by~\eqref{obs:projection} that $\< u,v \> > 0$. Thus $B$ is an $\SS^{d-2}$-bounding set, as claimed, and this completes the proof of the lemma.
\end{proof}

\subsection{Constructing the bounding tree} 

We are now ready to show that there exists an $\S$-good bounding tree. 
However, the existence of such trees will not be sufficient for our purposes; the tree also needs to be chosen so that the set $\H_{u_1}^d \cup \dots \cup \H_{u_{r-1}}^d$ is $\U$-closed for every path $\0 \to u_1 \to \cdots \to  u_{r-1}$ in $T$. We shall show in Proposition~\ref{lem:trees-exist} that a tree with this property exists.  
The slightly technical statement of Lemma~\ref{lem:trees-induction} is designed to facilitate the induction step in the proof of Proposition~\ref{lem:trees-exist}. In particular, the following variant of Definition~\ref{def:tree-good} will be required. 

\begin{definition}\label{def:tree-partially-good}
A bounding tree $\B = (T,\eta)$ of depth $s$ is \emph{partially $\S$-good} if, for every directed path $\0 \to u_1 \to \cdots \to u_s$ in $T$, we have 
\begin{equation}\label{eq:condition:d:repeat}
\rho^{d-1}(\SS^{d-1}; \S, u_i) \ge r - i
\end{equation}
for each $1 \le i \le s$, and moreover, for every $x \in X \in \U$ and $1 \le i < s$, we have either
\begin{equation}\label{eq:tree-partially-good-cond}
u_i \in \{x\}^\perp \qquad \text{or} \qquad \| u_i - \{x\}^\perp \| > r \cdot \eta(i).
\end{equation}
\end{definition}

Observe that a partially $\S$-good bounding tree of depth $r - 1$ is $\S$-good. Note also that condition~\eqref{eq:tree-partially-good-cond} is satisfied automatically as long as we choose $\eta(i)$ sufficiently small depending on $u_i$ (see Definition~\ref{def:tree}). 


We first show that there exists a partially $\S$-good bounding tree $\B = (T,\eta)$ of depth~1, which is the base case of our induction. 

\begin{lemma}\label{lem:trees-induction-basecase}
There exists a partially $\S$-good bounding tree of depth $1$.
\end{lemma}

\begin{proof}
To construct this tree, we simply take as our collection of neighbours of $\0$ a finite $\SS^{d-1}$-bounding set $B \subset \S$ such that $\rho^{d-1}(\SS^{d-1}; \S, u) \ge r - 1$ for all $u \in B$. The existence of such a set follows from our assumption that $r(\U) = r$; indeed, by~\eqref{eq:r} the set
$$\Big\{ \big\{ v \in \SS^{d-1} : \< u,v \> > 0 \big\} \,:\, u \in \SS^{d-1} \text{ and } \rho^{d-1}(\SS^{d-1}; \S, u) \ge r - 1 \Big\}$$
is an open cover of $\SS^{d-1}$, and we can let $B$ be the set of centres of the open hemispheres in a finite sub-cover. Note that the set $\{ u \}$ is automatically linearly independent for each $u \in B$, the function $\eta$ has empty domain, so does not need to be specified, and the condition~\eqref{eq:tree-partially-good-cond} is empty when $s = 1$. 
\end{proof}

We now use the technical lemmas from Section~\ref{sec:memphis} to prove the induction step.

\begin{lemma}\label{lem:trees-induction}
Every partially $\S$-good bounding tree $\B = (T,\eta)$ of depth $1 \le s \le r - 2$ can be extended to a partially $\S$-good bounding tree $\B' = (T',\eta')$ of depth $s + 1$ with $\eta'(s)$ arbitrarily small. In particular, there exists an $\S$-good bounding tree.
\end{lemma}

\begin{proof}
Fix $1 \le s \le r - 2$ and a partially $\S$-good bounding tree $\B = (T,\eta)$ of depth $s$. Set $\eta'(i) = \eta(i)$ for every $1 \le i \le s - 1$, and choose $\eta'(s)$ sufficiently small. To construct the bounding tree $\B' = (T',\eta')$, we need to choose a neighbourhood 
$$N_{T'}^\to(u) \subset \SS_u := S_{\eta'(s)}(\SS^{d-1}, u)$$
for each leaf $u$ of $T$. This neighbourhood must be an $\SS_u$-bounding set, and each of its elements must satisfy~\eqref{eq:condition:d:repeat} and must not lie in the span of the elements on the path in $T$ from $\0$ to $u$. As we noted above, the condition~\eqref{eq:tree-partially-good-cond} in Definition~\ref{def:tree-partially-good} is automatically satisfied as long as we choose $\eta'(s)$ sufficiently small (depending on $T$).

To construct these neighbourhoods, let $\0 \to u_1 \to \cdots \to u_s = u$ be a directed path in $T$, and observe that, by~\eqref{eq:rho} and since $\B$ is partially $\S$-good, we have  
\begin{equation}\label{eq:trees:induction:bounds}
r^{d-1}\big( \SS_u; \S ) = \rho^{d-1}(\SS^{d-1}; \S, u) \ge r - s.
\end{equation}
We shall use Lemma~\ref{lem:memphis-subspace} to prove the following claim.

\begin{claim}\label{claim:trees-induction}
There exists a finite $\SS_u$-bounding set $B_u$ with
\begin{equation}\label{eq:memphis-rho}
\rho^{d-2}\big( \SS_u; \S, v \big) \ge r^{d-1}\big( \SS_u; \S \big) - 1
\end{equation}
for every $v \in B_u$, and such that $B_u \cap \< W \> = \emptyset$, where $W := \{u_1,\dots,u_s\}$. 
\end{claim}

\begin{clmproof}{claim:trees-induction}
Set $\T := \varphi_u(\S \cap \SS_u)$, where $\varphi_u \colon \SS_u \to \SS(\{u\})$ is the homothety defined by $\varphi_u(x) = \pi(x,\{u\}^\perp)$. 
We will apply Lemma~\ref{lem:memphis-subspace} to the set $\T$, and then map the $\SS(\{u\})$-bounding set given by the lemma to a subset of $\SS_u$ using $\varphi_u^{-1}$.
 
Observe first that $\varphi_u(\SS_u) = \SS(\{u\})$ is a copy of $\SS^{d-2}$, that $\SS(W) \subset \SS(\{u\})$, since $u \in W$, and that $W^\perp \ne \{\0\}$, since $|W| = s < r \le d$. In order to apply Lemma~\ref{lem:memphis-subspace}, we therefore only need to show that $\T$ is an $\SS(\{u\})$-stable set, and that $C(\SS(\{u\});\T) \subset W^\perp$ for some valid choice of the set of centres of $\T$ with respect to $\SS(\{u\})$.

To choose the set $C(\SS(\{u\});\T)$, observe first that, by Lemma~\ref{lem:hemispheres}, we have 
\begin{equation}\label{eq:choice:of:C:proof}
\S \cap \SS_u = \bigcap_{X \in\, \U} \bigcup_{x \in X} \big( H_x \cap \SS_u \big),
\end{equation}
where $H_x = \{ v \in \SS^{d-1} : \< x,v\> \ge 0 \}$. Moreover, since $\U$ is a finite collection of finite sets and $\eta'(s)$ is sufficiently small, for each $x \in X \in \U$ we have either $u \in \{x\}^\perp$ or $\SS_u \cap \{x\}^\perp = \emptyset$. We claim that if $u \in \{x\}^\perp$ then $x \in W^\perp$. 
Indeed, if $u_i \not\in \{x\}^\perp$ for some $1 \le i < s$ then, since $\B$ is partially $\S$-good, we have $\| u_i - \{x\}^\perp \| > r \cdot \eta'(i)$, by~\eqref{eq:tree-partially-good-cond}. Recalling that $u = u_s$, it follows that
$$\| u - \{x\}^\perp \| \ge \| u_i - \{x\}^\perp \| - \| u_i - u_s \| > r \cdot \eta'(i) - \sum_{j = i}^{s-1} \eta'(j) > 0,$$ 
since $\|u_j - u_{j+1}\| = \eta'(j)$ for each $1 \le j \le s - 1$, and since $\eta'$ is decreasing and $s \le r - 1$. This contradicts our assumption that $u \in \{x\}^\perp$, and hence $x \in W^\perp$, as claimed. 

We have shown that if $x \not\in W^\perp$, then $\SS_u \cap \{x\}^\perp = \emptyset$, which in turn implies that $H_x \cap \SS_u \in \{ \emptyset, \SS_u \}$. Now, if $H_x \cap \SS_u = \emptyset$ then we do not need to include $x$ in~\eqref{eq:choice:of:C:proof}, and if $H_x \cap \SS_u = \SS_u$ then we do not need to include $X$ in~\eqref{eq:choice:of:C:proof}.  Therefore
\begin{equation}\label{eq:choice:of:C:refined}
\S \cap \SS_u = \bigcap_{X' \in \,\U'} \bigcup_{x \in X'} \big( H_x \cap \SS_u \big),
\end{equation}
where $\U'$ is the family of sets $X' = \{ x \in X : H_x \cap \SS_u \ne \emptyset \}$ such that $X \in \U$ and $H_x \cap \SS_u \ne \SS_u$ for every $x \in X$. 
In particular, $x \in W^\perp$ for every $x \in X' \in \U'$. 

Applying $\varphi_u$ to both sides of~\eqref{eq:choice:of:C:refined}, we obtain
\[
\T = \varphi_u\big( \S \cap \SS_u \big) = \bigcap_{X' \in \,\U'} \bigcup_{x \in X'} \varphi_u\big( H_x \cap \SS_u \big). 
\]
Now, by~\eqref{obs:projection}, for each $x \in X' \in \U'$ we have  
\begin{equation}\label{eq:new-centres}
\varphi_u\big( H_x \cap \SS_u \big) = \big\{ \varphi_u(v) : v \in \SS_u, \, \<x,v\> \geq 0 \big\} = \big\{ z \in \SS(\{u\}) : \<x,z\> \geq 0 \big\},
\end{equation}
since $\varphi_u(v) = \pi(v,\{u\}^\perp)$ and $x \in \{u\}^\perp$. It follows that $\T$ is an $\SS(\{u\})$-stable set, and that there exists a valid choice for $C(\SS(\{u\});\T)$ with 
\begin{equation}\label{eq:choice:of:C}
C(\SS(\{u\});\T) \subset \bigcup_{X' \in\, \U'} X' \subset W^\perp,
\end{equation}
as claimed. By Lemma~\ref{lem:memphis-subspace}, it follows that there exists an $\SS(\{u\})$-bounding set $B'$ with 
\begin{equation}\label{eq:memphis-subspace-rho:repeat}
\rho^{d-2}\big( \SS(\{u\}); \T, v \big) \geq r^{d-1}\big( \SS(\{u\}) ;\T \big) - 1
\end{equation}
for every $v \in B'$, and such that $B' \cap \<W\> = \emptyset$.

Observe that $B := \varphi_u^{-1}(B')$ is an $\SS_u$-bounding set, since $B'$ is an $\SS(\{u\})$-bounding set and $\varphi_u$ maps hemispheres of $\SS_u$ to hemispheres of $\varphi_u(\SS_u) = \SS(\{u\})$ (as in~\eqref{eq:new-centres}). Moreover,~\eqref{eq:memphis-rho} holds for every $v \in B$, since~\eqref{eq:memphis-subspace-rho:repeat} holds for every $v \in B'$, and $\varphi_u$ is a homothety mapping $\S \cap \SS_u$ to $\T$, so
$$\rho^{d-2}\big( \SS_u; \S, v \big) = \rho^{d-2}\big( \varphi_u(\SS_u); \T, \varphi_u(v) \big) \geq r^{d-1}\big( \varphi_u(\SS_u) ;\T \big) - 1 = r^{d-1}\big( \SS_u; \S \big) - 1$$
for every $v \in B$. Observe also that if $x \in B \cap \<W\>$, then 
$$\varphi_u(x) = \pi(x,\{u\}^\perp) = w + \lambda x \in B' \cap \<W\>$$
for some $w \in \<u\>$ and $\lambda > 0$, since $u \in W$. But $B' \cap \<W\> = \emptyset$, and thus also $B \cap \<W\> = \emptyset$.

Finally, since $B$ is an $\SS_u$-bounding set, the set 
$$\Big\{ \big\{ w \in \SS_u : \< w,v \> > 0 \big\} \,:\, v \in B \Big\}$$
is an open cover of $\SS_u$. Choose a finite sub-cover, and let $B_u$ be the set of centres of the open hemispheres of $\SS_u$ in this sub-cover. Since $B_u \subset B$, it follows that $B_u$ is a finite $\SS_u$-bounding set, $B_u \cap \<W\> = \emptyset$, and~\eqref{eq:memphis-rho} holds for every $v \in B_u$, as required. 
\end{clmproof}

Now, define $\B' = (T',\eta')$ by adding, for each leaf $u$ of $T$, the set $B_u$ as the out-neighbourhood of $u$ in $T'$. We claim that $\B'$ is a partially $\S$-good bounding tree. To see that $\B'$ is a bounding tree, observe that, by Definition~\ref{def:tree}, we only need to check that if $\0 \to u_1 \to \cdots \to u_s = u$ is a directed path in $T$, then the set $\{u_1,\ldots,u_s,v\}$ is linearly independent for every $v \in B_u$, which follows because $B_u \cap \< \{u_1,\dots,u_s\} \> = \emptyset$, and that $B_u$ is a finite $\SS_u$-bounding set, which is the case by construction. To see that $\B'$ is partially $\S$-good, observe first that, by~\eqref{eq:trees:induction:bounds} and~\eqref{eq:memphis-rho}, we have
$$\rho^{d-2}( \SS_u; \S, v ) \ge r - s - 1$$
for every $v \in B_u \subset \SS_u$. Now, by Lemma~\ref{lem:memphis}, we have 
$$\rho^{d-1} \big( \SS^{d-1}; \S, v \big) \ge \rho^{d-2} \big( \SS_u; \S, v \big),$$
and hence~\eqref{eq:condition:d:repeat} holds for $v$. Since $\B$ was assumed to be partially $\S$-good and $\eta'(s)$ was chosen sufficiently small, it follows that $\B'$ is partially $\S$-good, as required.  

Finally, to deduce that there exists an $\S$-good bounding tree, we simply apply the first part of the lemma $r - 2$ times to the partially $\S$-good bounding tree of depth 1 given by Lemma~\ref{lem:trees-induction-basecase}. We obtain a partially $\S$-good bounding tree of depth $r - 1$, which is therefore $\S$-good, as required.
\end{proof}

We can now complete the proof of the main result of this section, which states that 
there exists an $\S$-good bounding tree $\B = (T,\eta)$ with the property that unions of half-spaces along paths in $T$ are $\U$-closed. This property will play an important role in allowing us to control the growth of the $\U$-bootstrap process on the side of a droplet using icebergs (see Section~\ref{sec:icebergs}). In order to construct a tree with this property, we need to choose a magnification function $\eta$ that decreases sufficiently quickly, depending not only on the stable set $\S$, but also on the range $R$ of the update family $\U$ (see~\eqref{eq:R}). The flexibility in Lemma~\ref{lem:trees-induction} of choosing $\eta'(s)$ arbitrarily small will be useful for this purpose. 


In order to prove the next lemma when $\Lambda \ne \Z^d$, it will be useful to write
$$\HH_u^d := \{ x \in \R^d : \< x,u \> < 0 \}$$
for the continuous half-space with normal $u$. We remark that $\HH_u^d$ does not depend on $\Lambda$ and that $\H_u^d = \HH_u^d \cap \Lambda$.

\begin{prop}\label{lem:trees-exist}
There exists an $\S$-good bounding tree $\B = (T,\eta)$, not depending on $\Lambda$, such that, for every directed path $\0\to u_1\to \cdots \to u_{r-1}$ in $T$, and every $X \in \U$, 
$$X \not\subset \HH_{u_1}^d \cup \dots \cup \HH_{u_{r-1}}^d.$$ 
In particular, the set\/ $\H_{u_1}^d \cup \dots \cup \H_{u_{r-1}}^d$ is\/ $\U$-closed.
\end{prop}

\begin{proof}
We shall build the bounding tree $\B$ inductively using Lemma~\ref{lem:trees-induction}. Our induction statement is that, for each $1 \leq s \leq r - 1$, there exists a partially $\S$-good bounding tree $\B_s = (T_s,\eta_s)$ of depth $s$, such that the following holds:
\begin{itemize}
\item[$(+)$] For any path $\0 \to u_1 \to \cdots \to u_s$ in $T_s$, if
\[
x \in \HH_{u_1}^d \cup \dots \cup \HH_{u_s}^d
\]
for some $x \in X \in \U$, then $x \in \HH_{u_s}^d$.
\end{itemize}
(It is perhaps worth emphasizing here that there is no symmetry between the $u_i$, since they are being chosen successively depending on those directions already chosen.)

Let us first deduce the lemma from the induction hypothesis with $s = r - 1$. Note that $\B_{r-1}$ is an $\S$-good bounding tree, since it is partially $\S$-good by construction and has depth $r - 1$. 
In particular, this implies that $\rho^{d-1}(\SS^{d-1}; \S, u_{r-1}) \ge 1$, by~\eqref{eq:condition:d}, and thus that $u_{r-1} \in \S$. Now, let $\0 \to u_1 \to \cdots \to u_{r-1}$ be a directed path in $T_{r-1}$, and suppose that $X \subset \HH_{u_1}^d \cup \dots \cup \HH_{u_{r-1}}^d$ for some $X \in \U$. By~$(+)$, it follows that $X \subset \HH_{u_{r-1}}^d$, which contradicts\footnote{Recall that, by Lemma~\ref{lem:shift:Ssame}, the set $\S$ does not depend on our choice of the lattice $\Lambda$.} the fact that $u_{r-1}$ is a stable direction. 
This contradiction implies that there does not exist an update family $X \in \U$ with $X \subset \HH_{u_1}^d \cup \dots \cup \HH_{u_{r-1}}^d$, as required. 

To deduce that $\H_{u_1}^d \cup \dots \cup \H_{u_{r-1}}^d$ is\, $\U$-closed, observe that if there exists $y \in \Lambda$ with
\begin{equation}\label{eq:transfer-to-Lambda}
y + X \subset \H_{u_1}^d \cup \dots \cup \H_{u_{r-1}}^d \qquad \text{and} \qquad y \not\in \H_{u_1}^d \cup \dots \cup \H_{u_{r-1}}^d,
\end{equation} 
then $X \subset \HH_{u_1}^d \cup \dots \cup \HH_{u_{r-1}}^d$. Indeed, for each $x \in X$ there exists $i \in [r-1]$ such that $\< x + y, u_i \> < 0$ and $\< y, u_i \> \geq 0$, and therefore $\< x, u_i \> < 0$. 

To prove the induction statement, recall first that, by Lemma~\ref{lem:trees-induction-basecase}, there exists a partially $\S$-good bounding tree $\B_1 = (T_1,\eta_1)$ of depth $1$. Since the extra condition of property $(+)$ holds automatically when $s = 1$, this suffices to prove the base case.

Now, let $1 \le s \le r - 2$, and suppose that we have constructed a partially $\S$-good bounding tree $\B_s = (T_s,\eta_s)$ of depth $s$ that satisfies $(+)$. By Lemma~\ref{lem:trees-induction}, there exists a partially $\S$-good bounding tree $\B_{s+1} = (T_{s+1},\eta_{s+1})$ of depth $s+1$ that extends $\B_s$, and with $\eta_{s+1}(s)$ arbitrarily small. Let $L(T_s)$ denote the set of leaves of $T_s$, and observe that
$$\max\Big\{ \< x,u \> \,:\, x \in \HH_u^d \cap B_R(\0), \, u \in L(T_s) \Big\} < 0,$$
where $B_R(\0) := \big\{ x \in \Z^d : \| x \| \le R \big\}$ is the discrete Euclidean ball of radius $R$ centred at~$\0$, since $\< x,u \> < 0$ for each $x \in \HH_u^d$, and since the sets $B_R(\0)$ and $L(T_s)$ are finite.\footnote{We stress that this definition of $B_R(\0)$ as a subset of $\Z^d$ rather than as a subset of $\Lambda$ is intentional, and is because $X$ itself is a subset of $\Z^d$, and our aim is to prove $(+)$. The deduction of the final assertion of the proposition (which transfers from $\HH_u^d$ to $\H_u^d$) has already been proved above, around~\eqref{eq:transfer-to-Lambda}.}
Hence, by the continuity of $\< x,\cdot \>$, we may choose $\delta > 0$ such that
\begin{equation}\label{eq:seaind}
\max \Big\{ \< x,v \> \,:\, x \in \HH_u^d \cap B_R(\0), \, u \in L(T_s), \, v \in \SS^{d-1}, \, \|u - v\| \le \delta \Big\} < 0.
\end{equation}
Choose $\B_{s+1}$ with $\eta_{s+1}(s) \le \delta$, and note that this choice of $\B_{s+1}$ has not depended on $\Lambda$. 

It remains to prove that $\B_{s+1}$ satisfies $(+)$, so let $\0 \to u_1 \to \cdots \to u_{s+1}$ be a path in $T_{s+1}$, and suppose that
$$x \in \HH_{u_1}^d \cup \cdots \cup \HH_{u_{s+1}}^d$$
for some $x \in X \in \U$. 
Then either $x \in \HH^d_{u_{s+1}}$ (in which case we are done), or by the induction hypothesis we have $x \in \HH^d_{u_s}$. Now, note that $x \in B_R(\0)$, by~\eqref{eq:R}, that $u_s \in L(T_s)$, and that $\|u_s - u_{s+1}\| = \eta_{s+1}(s) \le \delta$. Therefore, by~\eqref{eq:seaind}, if $x \in \HH^d_{u_s}$, then $\< x,u_{s+1} \> < 0$, and so $x \in \HH^d_{u_{s+1}}$, as required. This completes the induction step, and hence the proof of the lemma.
\end{proof}

\subsection{The tree $T$}\label{subsec:tree}

{\setstretch{1.12}

Let us fix, for the rest of the paper (and, in particular, for the rest of the proof of Theorem~\ref{thm:lower}), an $\S$-good bounding tree $\B=(T,\eta)$ satisfying the conclusion of Proposition~\ref{lem:trees-exist}, and with $\eta(1)$ sufficiently small. (In fact, with the exception of the proof of Lemma~\ref{lem:icebergs-finite} below, we shall only need the tree, $T$, and can ignore from now on the magnification function $\eta$, which was only needed in the inductive construction of $T$.) 

We shall use the following notation and conventions in connection with $T$:
\begin{itemize}
\item Frequently we stipulate that $\0 \to u_1 \to \cdots \to u_k$ is a path in $T$; when we do so, we always mean a directed path, we allow $k = 0$, and we set $u_0 := \0$.\smallskip
\item We write $d_T(\0,u)$ for the graph distance (or, equivalently, the depth) of $u$ from the root $\0$ in $T$. Thus, $d_T(\0,u) \leq r - 1$ for all $u \in V(T)$.\smallskip
\item If $\0 \to u_1 \to \cdots \to u_k$ is a path in $T$, then we write
\begin{equation}\label{eq:path}
[\0,u_k]_T := \{ \0, u_1, \dots, u_k \}.
\end{equation}
\end{itemize}

The following sets will be used in the next section to define our `icebergs'.

\begin{definition}\label{def:MT}
Set $M_T(\0) := M_T^+(\0) := N_T^\to(\0)$. If $u$ is a non-leaf vertex of $T$ with $u \ne \0$, then define 
\begin{equation}\label{eq:MTu}
M_T(u) := N_T^\to(u) \cup \{-u\}
\end{equation}
and 
\begin{equation}\label{eq:MTplusu}
M_T^+(u) := N_T^\to(u) \cup \{-u_1,\dots,-u_k\},
\end{equation}
where $\0 \to u_1 \to \cdots \to u_k = u$ is the path in $T$ from $\0$ to $u$. 
\end{definition}

The following property of $M_T(u)$ will be useful in Sections~\ref{sec:icebergs} and~\ref{sec:subadd}. 

\begin{lemma}\label{lem:icebergs-finite}
Let $u$ be a non-leaf vertex of $T$. Then $M_T(u)$ is an $\SS^{d-1}$-bounding set. Moreover, if $u$ is not the root of $T$, then $N_T^\to(u)$ is not an $\SS^{d-1}$-bounding set.
\end{lemma}

\begin{proof}
It is straightforward to show that $N_T^\to(u)$ is not an $\SS^{d-1}$-bounding set, since it does not intersect the open hemisphere centred at $-u$. Indeed, if $u$ is not the root, then by~\eqref{eq:tree-nbours} every element of $N_T^\to(u)$ lies on a sphere of radius at most $\eta(1)$ centred at $u$. Since we assumed that $\eta(1)$ is sufficiently small, it follows that $\< u,v \> > 0$ for each $v \in N_T^\to(u)$.

To show that $M_T(u)$ is an $\SS^{d-1}$-bounding set, we let $H := \{ x \in \SS^{d-1} : \< x,v \> > 0 \}$ for some $v \in \SS^{d-1}$ be an arbitrary open hemisphere, and we claim that $M_T(u) \cap H$ is non-empty. Note that we are done immediately if $u = \0$, by Definition~\ref{def:tree}, or if $\< u,v \> < 0$, since $-u \in M_T(u)$. We are also done if $u = v$, because in that case (as observed above) we have $N_T^\to(u) \subset H$, since we assumed that $\eta(1)$ is sufficiently small. 

We may therefore assume that $v = \lambda u + w$, with $\lambda \ge 0$ and $\0 \ne w \in \{u\}^\perp$. Recall from~\eqref{eq:tree-nbours} that $\SS_u = S_\eta(\SS^{d-1},u)$ for some $\eta > 0$. We claim that
\begin{equation}\label{eq:icebergs-finite-u-in}
\big\{ x \in \SS_u : \< x, w \> > 0 \big\} \subset H.
\end{equation}
Note that the left-hand side of~\eqref{eq:icebergs-finite-u-in} is an open hemisphere in $\SS_u$, since $\0 \ne w \in \{u\}^\perp$, and hence has non-empty intersection with $N_T^\to(u)$, since $N_T^\to(u)$ is an $\SS_u$-bounding set, by Definition~\ref{def:tree}. It therefore suffices to prove~\eqref{eq:icebergs-finite-u-in}. To do so, let 
$x = \mu u + y \in \SS_u$, with $y \in \{u\}^\perp$, be such that $\< x, w \> > 0$. Then 
$$\< x,v \> = \mu \lambda + \< x,w \> > 0,$$
since $\lambda,\mu \ge 0$, and hence $x \in H$, as required. This proves~\eqref{eq:icebergs-finite-u-in}, and the lemma follows. 
\end{proof}

}

\section{Icebergs}\label{sec:icebergs}

In the next two sections we shall define a family of $\U$-closed structures (`icebergs') that generalize droplets, and develop a set of tools that we will use to control the growth of the $\U$-bootstrap process acting on the side of a droplet. In particular, in Section~\ref{sec:subadd} we will prove a sub-additivity property for icebergs (see Lemma~\ref{lem:subadd}), and use it to deduce an `Aizenman--Lebowitz-type' lemma and an extremal lemma (Lemmas~\ref{lem:AL} and~\ref{lem:extremal}, respectively). These three results are the basic tools that we will need in order to control the growth of a droplet using icebergs. This general approach can be traced back to the groundbreaking work of Aizenman and Lebowitz~\cite{AL} on the 2-neighbour model. However, we emphasize once again that the results in these sections are not simply technical generalizations of the results in two dimensions (or for the $r$-neighbour model), and that the development of these tools will require a number of very substantial innovations.  

In this section we lay the groundwork for the proofs of these key lemmas by making a number of important definitions, and proving various basic properties. In particular, we define icebergs (see Definition~\ref{def:iceberg}), `iceberg containers' (see Definition~\ref{def:iceberg-container}), and the notion of a `$u$-iceberg spanned iceberg' (see Definition~\ref{def:iceberg spanned}). We shall also define the `$u$-iceberg spanning algorithm' (see Definition~\ref{def:icespanalg}), which will play a role analogous to that of the rectangles process of Aizenman and Lebowitz~\cite{AL}. 

The basic purpose of icebergs is to enable us to control growth on the side of a droplet, which may itself be growing on the side of a (higher dimensional) droplet, and so on. In order to do so, we work in a more `generous' setting in which entire half-spaces are assumed to be already infected; in other words, we consider growth that is `assisted' by certain unions of half-spaces. In this section (and also in Section~\ref{sec:subadd}), these unions of half-spaces will always be of the following form: if $\0 \to u_1 \to \cdots \to u_k$ is a path in $T$,\footnote{Recall that we fixed in Section~\ref{subsec:tree} an $\S$-good bounding tree $\B = (T,\eta)$ satisfying the conclusion of Proposition~\ref{lem:trees-exist}.} 
then we define
$$H_T(u_k) := \bigcup_{i=1}^k \, \H_{u_i}^d = 
\bigcup_{i=1}^k \big\{ x \in \Lambda : \< x,u_i \> < 0 \big\},$$
so, in particular, $H_T(\0) = \emptyset$. We remark that in Section~\ref{one:step:sec} we will need to consider a slight generalization of the objects defined in this section, in which the boundaries of the half-spaces above are all translated by some element of $\R^d$ (which might not be in $\Z^d$). It is for this reason that we work in the shifted lattice $\Lambda$, instead of in $\Z^d$. 

Running the $\U$-bootstrap process with the assistance of the half-spaces in $H_T(u)$ leads us to the following definition. 

\begin{definition}\label{def:u-closure}
Given a vertex $u$ of $T$ and a finite set $K \subset \Lambda \setminus H_T(u)$, let
\[
[K]_u := \big[ K \cup H_T(u) \big]_\U \setminus H_T(u)
\]
be the \emph{$u$-closure} of $K$. We say that $K$ is \emph{$u$-closed} if $[K]_u = K$.
\end{definition}

Recall from Section~\ref{sec:outline} that we control the growth of the infected set using $\T$-droplets, which are $\U$-closed (and hence $\0$-closed) if $\T = N_T^\to(\0)$. For other vertices of $T$ we shall control the growth using icebergs, which are defined as follows.  

\begin{definition}\label{def:iceberg}
Let $u$ be a non-leaf vertex of $T$. A \emph{$u$-iceberg droplet} is a non-empty set of the form $D \setminus H_T(u)$, where $D$ is an $N_T^\to(u)$-droplet. An \emph{iceberg} is a pair $\J = (J,u)$, where $J$ is a $u$-iceberg droplet. 
\end{definition}

If we wish to emphasize $u$, then we shall say that $\J$ is a \emph{$u$-iceberg}, or that it is an \emph{iceberg of type~$u$}. Crucially, $u$-iceberg droplets are $u$-closed, as will be shown shortly in Lemma~\ref{lem:iceberg-drop-closed}. Before that, let us observe that they are finite. 

\begin{lemma}\label{cor:icebergs-finite}
If $u$ is a non-leaf vertex of\/ $T$ and $(J,u)$ is an iceberg, then $J$ is finite.
\end{lemma}

\begin{proof}
Recalling Definition~\ref{def:MT}, 
note that every $u$-iceberg droplet is also an $M_T^+(u)$-droplet, and that $M_T(u) \subset M_T^+(u)$. Since, by Lemma~\ref{lem:icebergs-finite}, $M_T(u)$ is an $\SS^{d-1}$-bounding set, it follows that $u$-iceberg droplets are finite, as claimed.
\end{proof}

Next we prove the key fact that $u$-iceberg droplets are $u$-closed, which follows from the property of the tree $T$ guaranteed by Proposition~\ref{lem:trees-exist}. The proof is slightly complicated by the fact that we are working in $\Lambda$ rather than $\Z^d$. 

\begin{lemma}\label{lem:iceberg-drop-closed}
If\/ $u$ is a non-leaf vertex of\/ $T$ and $(J,u)$ is an iceberg, then $J$ is $u$-closed.
\end{lemma}

\begin{proof}
Suppose, for a contradiction, that we have 
\begin{equation}\label{eq:iceberg:u:closed:contradiction}
x \notin J \cup H_T(u)  \qquad \text{and} \qquad  x + X \subset J \cup H_T(u)
\end{equation}
for some $X \in \U$ and some $x \in \Lambda$, and recall from Section~\ref{subsec:tree} that we chose $T$ so that the conclusion of Proposition~\ref{lem:trees-exist} holds, and in particular so that for every path $\0 \to u_1 \to \cdots \to u_{r-1}$ in $T$, the set $\H_{u_1}^d \cup \dots \cup \H_{u_{r-1}}^d$ is $\U$-closed. Note that, 
by Definitions~\ref{def:droplet} and~\ref{def:iceberg}, we have
\begin{equation}\label{eq:iceberg:halfspaces}
J \cup H_T(u) = \bigg( \bigcap_{v \in N_T^\to(u)} \big( \H_v^d + a_v \big) \bigg) \cup H_T(u) = \bigcap_{v \in N_T^\to(u)} \Big( \big( \H_v^d + a_v \big) \cup H_T(u) \Big),
\end{equation}
for some collection $\{a_v \in \R^d : v \in N_T^\to(u)\}$. It follows from~\eqref{eq:iceberg:u:closed:contradiction} and~\eqref{eq:iceberg:halfspaces} that there exists $v \in N_T^\to(u)$ such that $\big( \H_v^d + a_v \big) \cup H_T(u)$ contains $x + X$ but does not contain $x$.

We claim that if $\0 \to u_1 \to \cdots \to u_{r-1}$ is a path in $T$ containing $u$ and $v$, then
\begin{equation}\label{eq:iceberg-closed-clm}
X \subset \HH_{u_1}^d \cup \dots \cup \HH_{u_{r-1}}^d,
\end{equation}
which contradicts our choice of $T$, by Proposition~\ref{lem:trees-exist}. To prove~\eqref{eq:iceberg-closed-clm}, let $y \in X$, and recall that either $x + y \in \H_v^d + a_v$, or $x + y \in \H_w^d$ for some $\0 \ne w \in [\0,u]_T$. We claim that in the former case we have $y \in \HH_v^d$, and that in the latter case we have $y \in \HH_w^d$. 

Indeed, if $x + y \in \H_v^d + a_v$, then $\< x + y - a_v, v \> < 0$, and also $\< x - a_v, v \> \ge 0$, since $x \notin \H_v^d + a_v$. It follows that $\< y, v \> < 0$, and so $y \in \HH_v^d$. Similarly, if $x + y \in \H_w^d$, then $\< x + y, w \> < 0$, and $\< x, w \> \ge 0$, since $x \notin H_T(u)$. Thus $\< y, w \> < 0$, and so $y \in \HH_w^d$. 

In either case, it follows that $y \in \HH_{u_1}^d \cup \dots \cup \HH_{u_{r-1}}^d$, and so this proves that~\eqref{eq:iceberg-closed-clm} holds. This contradicts our choice of $T$, and hence completes the proof of the lemma.
\end{proof}

We are now ready for an important (and slightly more technical) definition: that of the $u$-iceberg container of a set $K \subset \Lambda \setminus H_T(u)$. The definition is somewhat delicate, and has been chosen so that the corresponding droplets are $u$-closed (see Lemma~\ref{lem:iceberg-cont-closed}), and so that it exhibits various types of `monotonicity' (see Lemmas~\ref{lem:iceberg-cont-mono-type},~\ref{lem:iceberg-cont-reduce} and~\ref{lem:iceberg span-reduce}). We remark that iceberg containers will also play a central role in Sections~\ref{sec:subadd} and~\ref{one:step:sec}. 

\begin{figure}[ht]
  \centering
  \begin{tikzpicture}[>=latex] 
	\draw (-8,0) -- (0,0) -- (6,3);
	\draw (-2,1.2) circle (0.5);
	\draw [dotted, thick] (-2.87,1.7) -- (-1.13,1.7) -- (-2,0.2) -- cycle;
	\draw [densely dashed] (-5.6,0) -- (-2,1.8) -- (1.6,0) -- (0,0);
	\draw [densely dashed, thick] (-4.4,0) -- (-2.25,1.7) -- (4.8,2.4);
	\node at (-2,1.2) {$K$};
	\node at (4,0) {$\H_{u_1}^d \cup \H_{u_2}^d$};
	\draw [->] (-7,0) -- (-7,0.5) node [above] {$u_1$};
	\draw [->] (5.4,2.7) -- (5.2,3.1) node [above left] {$u_2$};
	\draw (-2.6,1.7) -- (-2.6,2.2) node [above] {$J_0$};
	\draw (-4,0.8) -- (-4,1.3) node [above] {$J_1$};
	\draw (-0.3,1.9) -- (-0.3,2.4) node [above] {$J_2$};
  \end{tikzpicture}
  \caption{The $u$-iceberg container of a set $K$. Three iceberg-droplets are shown: the minimal $u_i$-iceberg droplet $J_i$ containing $K$ for each $i \in \{0,1,2\}$. Because $J_0$ is strongly connected to $\HH_{u_1}^d$ and $J_1$ is strongly connected to $\HH_{u_2}^d$, the $u_2$-iceberg container of $K$ is $(J_2,u_2)$.}
  \label{fig:iceberg container}
\end{figure}
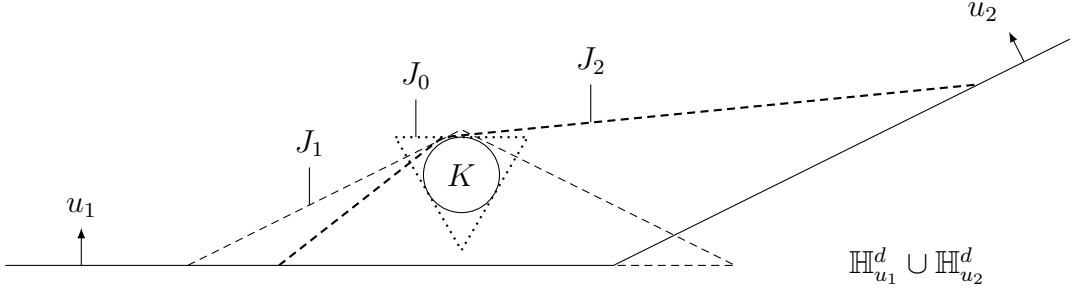

\begin{definition}\label{def:iceberg-container}
Let $\0 \to u_1 \to \cdots \to u_k$ be a path in $T$, with $u_k$ not a leaf, and let $K \subset \Lambda \setminus H_T(u_k)$ be a finite, non-empty set. For each $0 \leq i \leq k$, let $J_i$ be the minimal $u_i$-iceberg droplet containing $K$. We say that a sequence 
$$0 = i_0 < i_1 < \dots < i_\ell \leq k$$ 
is \emph{$(K, u_k)$-permissible} if, for each $0 \leq j < \ell$, the set $J_{i_j}$ is strongly connected to $\HH_{u_{i_{j+1}}}^d$.\footnote{Recall from Definition~\ref{def:strconn} that two sets are strongly connected to each other if they are at distance at most $R$. Recall also from Section~\ref{subsec:tree} that $u_0 = \0$, and note that the sequence $i_0 = 0$ of length 1 is always permissible, and that the type of a $u_k$-iceberg container is an element of the set $[\0,u_k]_T = \{ \0, u_1, \dots, u_k \}$.} Now define $\J_{u_k}(K)$, the \emph{$u_k$-iceberg container of $K$}, to be the iceberg $(J_i,u_i)$, where $i$ is the largest integer that appears in any $(K, u_k)$-permissible sequence.
\end{definition}

One should think of the permissible sequences in Definition~\ref{def:iceberg-container} as allowing us to take the type of the iceberg container $\J_{u_k}(K)$ to lie as far from the root $\0$ as possible, while still allowing us to maintain control over the size of $J_i$. We are able to maintain this control because $J_{i_j}$ being strongly connected to $\HH_{u_{i_{j+1}}}^d$ ensures that $J_{i_{j+1}}$ cannot be too much larger than $J_{i_j}$ (see Lemma~\ref{lem:wd-permissible}). This fact will be crucial in the proof of our sub-additivity lemma, Lemma~\ref{lem:subadd}.


An important property of $u$-iceberg containers is that their iceberg droplets are $u$-closed. This is a simple consequence of Lemma~\ref{lem:iceberg-drop-closed}. 

\begin{lemma}\label{lem:iceberg-cont-closed}
Let $u$ be a non-leaf vertex of $T$, and let $\J_u(K) = (J,v)$ be the $u$-iceberg container of a finite set $K \subset \Lambda \setminus H_T(u)$. Then $J \subset \Lambda \setminus H_T(u)$, and $J$ is $u$-closed.
\end{lemma}

\begin{proof}
Let $\0 \to u_1 \to \cdots \to u_k = u$ be the path in the tree $T$ from the root to $u$, and suppose that $v = u_i$. By Definition~\ref{def:iceberg-container}, it follows that $J = J_i$, the minimal $u_i$-iceberg droplet containing $K$, and in particular $J \subset \Lambda \setminus H_T(u_i)$, by Definition~\ref{def:iceberg}. Moreover, again by Definition~\ref{def:iceberg-container}, there exists a $(K,u)$-permissible path ending in $i$, but there does not exist one ending in $j$ for any $i < j \le k$, and therefore $J$ is not strongly connected to the set $\HH_{u_{i+1}} \cup \dots \cup \HH_{u_k}$. It follows that $J \subset \Lambda \setminus H_T(u)$.

To prove that $J$ is $u$-closed, suppose, for a contradiction, that we have 
$$x \notin J \cup H_T(u)  \qquad \text{and} \qquad  x + X \subset J \cup H_T(u)$$
for some $X \in \U$ and some $x \in \Lambda$. Recall that (by our choice of $T$) the set $H_T(u)$ is $\U$-closed, and therefore does not contain the set $x + X$. It follows that $x + X$ intersects $J$, and hence, since $J$ is not strongly connected to $\HH_{u_{i+1}} \cup \dots \cup \HH_{u_k}$,  
$$\big( x + X \big) \cap \big( \H_{u_{i+1}} \cup \dots \cup \H_{u_k} \big) = \emptyset.$$
But it now follows that $x + X \subset J \cup H_T(u_i)$, which means that $J$ is not $u_i$-closed. This contradicts Lemma~\ref{lem:iceberg-drop-closed}, and therefore proves the lemma.
\end{proof}

The following simple but very useful monotonicity property provides further motivation for Definition~\ref{def:iceberg-container}. 

\begin{lemma}\label{lem:iceberg-cont-mono-type}
Let $u$ be a non-leaf vertex of $T$, and let $K \subset K' \subset \Lambda \setminus H_T(u)$ be finite sets. If $\J_u(K)$ is a $v$-iceberg and $\J_u(K')$ is a $w$-iceberg, then $d_T(\0,v) \le d_T(\0,w)$.
\end{lemma}

\begin{proof}
The lemma follows immediately from the observation that any $(K,u)$-permissible sequence is also a $(K',u)$-permissible sequence, which in turn follows from the fact that for any $u' \in [\0,u]_T$, the minimal $u'$-iceberg droplet containing $K$ is a subset of the minimal $u'$-iceberg droplet containing $K'$. 
\end{proof}

We shall also need the following monotonicity-type properties. 


\begin{lemma}\label{lem:iceberg-cont-reduce}
Let $\0 \to u_1 \to \cdots \to u_k$ be a path in $T$, with $u_k$ not a leaf, and let $K \subset \Lambda \setminus H_T(u_k)$ be a finite set.
\begin{itemize}
\item[$(a)$] If $\J_{u_k}(K)$ is a $u_i$-iceberg, then $\J_{u_j}(K) = \J_{u_k}(K)$ for every $i \le j \le k$.\smallskip
\item[$(b)$] If $\J_{u_k}\big( [K]_{u_k} \big)$ is a $u_i$-iceberg, then 
$$[K]_{u_j} = [K]_{u_k} \qquad \text{and} \qquad \J_{u_j}\big( [K]_{u_j} \big) = \J_{u_k}\big( [K]_{u_k} \big)$$ 
for every $i \le j \le k$.
\end{itemize}
\end{lemma}

\begin{proof}
In order to prove $(a)$, we shall show that the $(K,u_j)$-permissible sequences are the same for every $i \le j \le k$. Indeed, if the sequence $0 = i_0 < \dots < i_\ell \leq k$ is $(K,u_k)$-permissible, then $i_\ell \le i \le j$, since $\J_{u_k}(K)$ is a $u_i$-iceberg. Since the iceberg droplets $J_0,\ldots,J_j$ are the same in both cases, it follows that a sequence is $(K,u_k)$-permissible  if and only if it is $(K,u_j)$-permissible. It follows immediately that $\J_{u_j}(K) = \J_{u_k}(K)$.

To deduce~$(b)$ from~$(a)$, it suffices to prove that $[K]_{u_j} = [K]_{u_k}$ for every $i \le j \le k$. To do so, we need to show that the set $[K]_{u_j} \cup H_T(u_k)$ is $\U$-closed, and that the set $[K]_{u_j} \cap H_T(u_k)$ is empty. Both of these properties will follow from the following claim: 
\begin{equation}\label{claim:not:strongly:connected}
[K]_{u_j} \text{ is not strongly connected to }\HH_{u_{j+1}}^d \cup \dots \cup \HH_{u_k}^d.
\end{equation} 
To prove~\eqref{claim:not:strongly:connected}, let $\J_{u_k}\big( [K]_{u_k} \big) = (J,u_i)$, and observe first that $\J_{u_j}\big( [K]_{u_k} \big) = (J,u_i)$, by part~$(a)$, and therefore $J$ is $u_j$-closed, by Lemma~\ref{lem:iceberg-cont-closed}. Since $K \subset J$, it follows that $[K]_{u_j} \subset [J]_{u_j} = J$. Now, to deduce~\eqref{claim:not:strongly:connected} simply note that $J$ is not strongly connected to $\HH_{u_{j+1}}^d \cup \dots \cup \HH_{u_k}^d$, by Definition~\ref{def:iceberg-container}, since otherwise there would be a $\big( [K]_{u_k}, u_k \big)$-permissible sequence containing an element strictly larger than $i$. 

It follows immediately from~\eqref{claim:not:strongly:connected} and Definition~\ref{def:u-closure} that $[K]_{u_j} \cap H_T(u_k)$ is empty. Finally, to show that $[K]_{u_j} \cup H_T(u_k)$ is $\U$-closed, suppose for a contradiction that 
$$x \not\in [K]_{u_j} \cup H_T(u_k) \qquad \text{and} \qquad x + X \subset [K]_{u_j} \cup H_T(u_k)$$
for some $X \in \U$ and $x \in \Lambda$. Recall that $[K]_{u_j} \cup H_T(u_j)$ is $\U$-closed (by definition) and that $H_T(u_k)$ is $\U$-closed (by our choice of $T$), and observe that therefore $x + X$ must intersect both $[K]_{u_j}$ and $\H_{u_{j+1}}^d \cup \dots \cup \H_{u_k}^d$. But this contradicts~\eqref{claim:not:strongly:connected}, and hence proves that $[K]_{u_j} \cup H_T(u_k)$ is $\U$-closed, as claimed. 

As noted above, this implies that $[K]_{u_j} = [K]_{u_k}$, and applying part~$(a)$ it follows that we also have $\J_{u_j}\big( [K]_{u_j} \big) = \J_{u_k}\big( [K]_{u_k} \big)$, as required.
\end{proof}

Having introduced iceberg containers, we are now ready for the following definition, which generalizes to the iceberg setting the property of being internally $\T$-spanned (see Definition~\ref{def:intspan}), and which will play a central role in Sections~\ref{one:step:sec} and~\ref{sec:ind}. 

\begin{definition}\label{def:iceberg spanned}
Let $u$ be a non-leaf vertex of $T$. An iceberg $\J = (J,v)$ is said to be \emph{$u$-iceberg spanned} if there exists $K \subset J \cap A$ such that $[K]_u$ is strongly connected and $\J = \J_u\big( [K]_u \big)$. We denote the event that $\J$ is $u$-iceberg spanned by $\ice_u(\J)$.
\end{definition}

Note that $\ice_u(\J)$ is an increasing event, and that if $\J = (J,v)$ is $u$-iceberg spanned then $v \in [\0,u]_T$, since $v$ is the type of a $u$-iceberg container. Note also that we may assume that $K \subset \Lambda \setminus H_T(u)$, since this does not affect the $u$-closure $[K]_u$. 

Our next task is to define a \emph{spanning algorithm}, which will be one our tools for understanding the event that an iceberg is $u$-iceberg spanned. The following algorithm generalizes the spanning algorithm of~\cite{BDMS} to the iceberg setting.\footnote{More precisely, setting $u = \0$ in the algorithm below gives~\cite[Definition~6.14]{BDMS}.}


\begin{definition}[The $u$-iceberg spanning algorithm]\label{def:icespanalg}
Let $u$ be a non-leaf vertex of $T$, let $K$ be a finite subset of $\Lambda \setminus H_T(u)$, and set $\K^0 := \big\{ \{x\} : x \in K \big\}$. Set $t := 0$, and repeat the following steps until STOP: 
\begin{itemize}
\item If there are two sets $K_1,K_2 \in \K^t$ such that the set
$$[K_1]_u \cup [K_2]_u$$
is strongly connected, then set
\[
\K^{t+1} := \big( \K^t \setminus \big\{ K_1, K_2 \big\} \big) \cup \big\{ K_1 \cup K_2 \big\},
\]
and set $t := t + 1$. 
\item Otherwise set $T := t$ and STOP.
\end{itemize}
The output of the algorithm is the \emph{$u$-iceberg span} of $K$, 
$$\< K \>_u := \big\{ \J_u\big( [K']_u \big) : K' \in \K^T \big\}.$$
\end{definition}

The alert reader may perhaps already have guessed that the $u$-iceberg span of $K$ is simply the collection of $u$-iceberg containers of the strongly connected components of~$[K]_u$. The following lemma records this simple but important fact.


\begin{lemma}\label{lem:output-of-span}
Let $u$ be a non-leaf vertex of $T$, and let $K \subset \Lambda\setminus H_T(u)$ be a finite set. Then 
\begin{equation}\label{eq:spanofK}
\< K \>_u = \big\{ \J_u(L_1), \ldots, \J_u(L_\ell) \big\},
\end{equation}
where $L_1,\ldots,L_\ell$ are the strongly connected components of $[K]_u$.
\end{lemma}

Before proving Lemma~\ref{lem:output-of-span}, let us first observe that if $K \subset \Lambda \setminus H_T(u)$ is strongly connected, then $[K]_u$ is also strongly connected. Indeed, this follows from the fact that $H_T(u)$ is $\U$-closed, by Proposition~\ref{lem:trees-exist}, and therefore every vertex of $\Lambda \setminus H_T(u)$ that is infected by the $\U$-process, with initial set $K \cup H_T(u)$, must be strongly connected to some already-infected vertex of $\Lambda \setminus H_T(u)$. We shall use this fact to show that every set $[K']_u$ that occurs in the $u$-iceberg spanning algorithm is strongly connected.
 
\begin{proof}[Proof of Lemma~\ref{lem:output-of-span}]
To prove the lemma, we need to show that the sets $\big\{ [K']_u : K' \in \K^T \big\}$ at the end of the algorithm are precisely the strongly connected components of $[K]_u$. To see this, note first that $[K']_u$ is strongly connected for every $0 \le t \le T$ and $K' \in \K^t$, since this is true at time $t = 0$, and if the set $[K_1]_u \cup [K_2]_u$ is strongly connected for some $K_1,K_2 \in \K^t$, then 
the set $[K_1 \cup K_2]_u$ is also strongly connected, by the observation above. 
Moreover, for any two (distinct) sets $K',K'' \in \K^T$, the sets $[K']_u$ and $[K'']_u$ are not strongly connected to each other, since the algorithm stopped at step $T$. 

It therefore only remains to show that $[K]_u = \bigcup_{i=1}^k [K_i^T]_u$. Clearly $[K]_u \supset \bigcup_{i=1}^k [K_i^T]_u$, so suppose for a contradiction that $x \in [K]_u$ for some $x \not\in \bigcup_{i=1}^k [K_i^T]_u$. It follows, by Definition~\ref{def:u-closure}, that $x \not\in H_T(u)$, and that $x + X \subset \bigcup_{i=1}^k [K_i^T]_u \cup H_T(u)$ for some $X \in \U$. Since $[K_i^T]_u \cup H_T(u)$ is $\U$-closed for each $i$ (by definition), the set $x + X$ must intersect at least two of the sets $[K_i^T]_u$. But this is impossible, by the definitions of $R$ and of being strongly connected, 
so we have the required contradiction.
\end{proof}



Let us also record, for future convenience, the following connection between Definitions~\ref{def:iceberg spanned} and~\ref{def:icespanalg}, which is an immediate consequence of Lemma~\ref{lem:output-of-span}. 

\begin{lemma}\label{lem:icebergspanned}
Let $u$ be a non-leaf vertex of $T$. Then an iceberg $\J = (J,v)$ is $u$-iceberg spanned if and only if $\< K \>_u = \{ \J \}$ for some $K \subset J \cap A$.
\end{lemma}

\begin{proof}
By Lemma~\ref{lem:output-of-span}, it follows that $\< K \>_u = \{ \J \}$ if and only if $[K]_u$ is strongly connected and $\J = \J_u\big( [K]_u \big)$. The fact that this holds for some $K \subset J \cap A$ is precisely the definition of the event that $\J$ is $u$-iceberg spanned.
\end{proof}

The following lemma is another easy (but important) consequence of the $u$-iceberg spanning algorithm. The proof is essentially identical to that of~\cite[Lemma~6.20]{BDMS}.

\begin{lemma}\label{lem:penultimate}
Let $u$ be a non-leaf vertex of\/ $T$ and let $K \subset \Lambda \setminus H_T(u)$ be a finite set with $|K| \ge 2$. If $[K]_u$ is strongly connected, then there exists a partition $K = K_1 \cup K_2$ into non-empty sets such that $[K_1]_u$, $[K_2]_u$ and $[K_1]_u \cup [K_2]_u$ are all strongly connected.
\end{lemma}

\begin{proof}
Run the $u$-iceberg spanning algorithm on $K$ and consider the penultimate step. Observe that, for each $0 \le t \le T$, the sets in $\K^t$ form a partition of $K$, and that the set $[K']_u$ is strongly connected for every $K' \in \K^t$ (see the proof of Lemma~\ref{lem:output-of-span}).

Now, since $[K]_u$ is strongly connected, by Lemma~\ref{lem:output-of-span} we have $\< K \>_u = \big\{ \J_u\big( [K]_u \big) \big\}$, and therefore 
\[
\K^{T-1} = \big\{ K_1, K_2 \big\}
\]
for some disjoint, non-empty sets $K_1,K_2 \subset K$ with $K_1 \cup K_2 = K$. By the observations above, both $[K_1]_u$ and $[K_2]_u$ are strongly connected. Finally, since $K_1$ and $K_2$ combine at the final step, the set $[K_1]_u \cup [K_2]_u$ must also be strongly connected, as required.
\end{proof}

To finish the section, let us note another monotonicity property, which will be used in Section~\ref{sec:ind}, and which is a straightforward corollary of Lemma~\ref{lem:iceberg-cont-reduce}. 

\begin{lemma}\label{lem:iceberg span-reduce}
Let $\0 \to u_1 \to \cdots \to u_k$ be a path in $T$, with $u_k$ not a leaf. If the iceberg $\J = (J,u_i)$ is $u_k$-iceberg spanned, then $\J$ is also $u_j$-iceberg spanned for every $i \leq j \leq k$.
\end{lemma}

\begin{proof}
By Definition~\ref{def:iceberg spanned}, if $\J$ is $u_k$-iceberg spanned, then there exists a set $K \subset J \cap A$ such that $[K]_{u_k}$ is strongly connected and $\J = \J_{u_k}\big( [K]_{u_k} \big) = (J,u_i)$. But then, by Lemma~\ref{lem:iceberg-cont-reduce}, it follows that $[K]_{u_j} = [K]_{u_k}$ is strongly connected and $\J_{u_j}\big( [K]_{u_j} \big) = \J$ for every $i \leq j \leq k$, and hence $\J$ is $u_j$-iceberg spanned, as claimed.
\end{proof}

\section{The size of an iceberg}\label{sec:subadd}

In this section we shall prove our main sub-additivity lemma for icebergs, Lemma~\ref{lem:subadd}. This lemma is an analogue for icebergs of the following inequality from~\cite{BDMS,BSU}: if $\T \subset \SS^1$ is a finite set, and if $D_1, D_2 \subset \Z^2$ are two $\T$-droplets with $D_1 \cap D_2 \neq \emptyset$, then
\begin{equation}\label{eq:2d-subadd}
\diam(D) \leq \diam(D_1) + \diam(D_2),
\end{equation}
where $D$ is the minimal $\T$-droplet containing $D_1 \cup D_2$. 
This inequality has several crucial consequences for the study of two-dimensional (critical) update families, including the so-called `Aizenman--Lebowitz lemma', which implies the existence of an internally spanned `critical' droplet, and the `extremal lemma', which gives a lower bound on $|D \cap A|$ for a droplet $D$ that is internally spanned by the set $A$. 

Variants of these lemmas that hold in the more general setting of icebergs (see Lemmas~\ref{lem:AL} and~\ref{lem:extremal}) will be required in Section~\ref{sec:ind}. Unfortunately, proving an analogue of~\eqref{eq:2d-subadd} for icebergs turns out to be significantly more difficult when $d \ge 3$; we therefore begin by briefly discussing some of the challenges that will be faced in this section, and how we shall go about overcoming them. 

The first obvious problem we encounter is that icebergs may have different types, and icebergs of different types have different shapes, and in such cases it is not obvious what the generalization of~\eqref{eq:2d-subadd} should be. However, there is a more fundamental problem, which is that, even restricting to the setting of droplets, the direct generalization of~\eqref{eq:2d-subadd} to dimensions 3 and higher is false: that is, there exists a finite set $\T \subset \SS^2$ and a pair $D_1, D_2 \subset \Z^3$ of intersecting $\T$-droplets such that
\begin{equation}\label{eq:diam-subadd-false}
\diam(D) > \diam(D_1) + \diam(D_2),
\end{equation}
where $D$ is the minimal $\T$-droplet containing $D_1 \cup D_2$.\footnote{To give a concrete example in three dimensions, consider tetrahedra $D_1$ and $D_2$ with corners $\{ (1,\pm 2c,0),(-1,0,\pm 2c) \}$, and $\{ (1,0,\pm 2c), (-1,\pm 2c,0) \}$. Each has diameter close to $2$ when $c$ is small, but $(3,0,0)$ and $(-3,0,0)$ are both contained in the smallest $\T$-droplet containing $D_1 \cup D_2$.}

In order to address this second problem, we shall introduce in Definition~\ref{def:weighted-diam} a new notion of the `size' of a $\T$-droplet $D$, which we denote by $\wdiam_\T(D)$ and call the \emph{$\T$-weighted diameter} of $D$. This function will be chosen so that $\wdiam_\T(D) = \Theta_\T\big( \diam(D) \big)$ for every droplet $D$ (see Lemma~\ref{lem:wd-prop-diam}), and so that a sub-additivity inequality similar to~\eqref{eq:2d-subadd} holds (see Lemma~\ref{lem:sub-add-droplets}). We shall then use the function $\wdiam_\T$ to define a measure, $\diam^*(\J)$, of the size of an iceberg $\J$ (see Definition~\ref{def:diamstar}). 

The main result of the section, Lemma~\ref{lem:subadd}, is a sub-additivity inequality for the measure $\diam^*(\J)$, but the statement is (necessarily) rather subtle, involving the $u$-iceberg containers introduced in the previous section (see Definition~\ref{def:iceberg-container}), as well as another notion of container (see Definition~\ref{def:FuGu}) which is carefully chosen to work well with $u$-closures and $(K,u)$-permissible sequences (see Lemmas~\ref{lem:wd-permissible} and~\ref{lem:wd-change-closures}). Once we have proved Lemma~\ref{lem:subadd}, it will be relatively straightforward\footnote{Unlike in the case $d = 2$, however, to do so we will require one additional technical lemma, whose proof is surprisingly challenging, see Lemma~\ref{lem:extremal:basecase}.} to deduce Lemmas~\ref{lem:AL} and~\ref{lem:extremal}, our `Aizenman--Lebowitz' and `extremal' lemmas for icebergs, respectively. 

In order to define $\wdiam_\T$, we need the following simple geometric lemma. Recall that an $\SS$-bounding set is a set $\T \subset \SS$ that intersects every open hemisphere of $\SS$. 

\begin{lemma}\label{lem:sub-add-convex}
Let\/ $\T = \{u_1,\dots,u_k\} \subset \SS^{d-1}$ be an\/ $\SS^{d-1}$-bounding set. Then there exist real numbers $\lambda_1,\dots,\lambda_k > 0$ such that\/ $\sum_{i=1}^k \lambda_i u_i = \0$ and\/ $\sum_{i=1}^k \lambda_i = 1$.
\end{lemma}

\begin{proof}
Suppose, for a contradiction, that $\0$ is not in the interior of $\conv(\T)$, the convex hull of $\T$. Since the interior of $\conv(\T)$ is also convex, it follows by the hyperplane separation theorem that there exists $v \in \SS^{d-1}$ such that $\< u,v \> \leq 0$ for all $u \in \T$, contradicting our assumption that $\T$ meets every open hemisphere of $\SS^{d-1}$. It follows that $\0$ is in the interior of $\conv(\T)$, and the conclusion of the lemma follows easily.
\end{proof}

For each finite $\SS^{d-1}$-bounding set $\T \subset \SS^{d-1}$, 
we may use Lemma~\ref{lem:sub-add-convex} to choose a set $\Xi_\T = \{ \lambda_v : v \in \T \}$, with $\lambda_v > 0$ for all $v \in \T$, such that
\begin{equation}\label{def:Lambda}
\sum_{v \in \T} \lambda_v v = \0 \qquad \text{and} \qquad \sum_{v \in \T} \lambda_v = 1.
\end{equation}
We consider the sets $\Xi_\T$ to be fixed for the rest of the paper. Moreover, for each non-leaf vertex $u$ of $T$, define
$$\Xi_u := \Xi_{M_T(u)},$$
where the set $M_T(u)$, defined in~\eqref{eq:MTu}, is an $\SS^{d-1}$-bounding set by Lemma~\ref{lem:icebergs-finite}.

\begin{definition}\label{def:weighted-diam}
Let $\T \subset \SS^{d-1}$ be a finite $\SS^{d-1}$-bounding set. For each bounded, non-empty set $\X \subset \R^d$, define the \emph{$\T$-weighted diameter} of $\X$ to be
\begin{equation}\label{eq:weighted-diam}
\wdiam_\T(\X) := \sum_{v \in \T} \lambda_v \cdot d_v(\X),
\end{equation}
where $d_v(\X) := \sup_{y \in \X} \< y, v \>$, and $\Xi_\T = \{ \lambda_v : v \in \T \}$ is the set chosen in~\eqref{def:Lambda}.

Moreover, for each non-leaf vertex $u$ of $T$, define $\wdiam_u := \wdiam_{M_T(u)}$. 
\end{definition}

In order to motivate the definition above, let us note first that if $\X \subset \Lambda$ and $D$ is the minimal $\T$-droplet containing $\X$, then $\wdiam_\T(D) = \wdiam_\T(\X)$. Indeed, $D$ consists of exactly the points $y \in \Lambda$ such that $\< y, v \> \le d_v(\X)$ for every $v \in \T$. Observe also that $\wdiam_\T$ is increasing, in the sense that if $\X \subset \Y$ then $\wdiam_\T(\X) \le \wdiam_\T(\Y)$. 

As further motivation, 
the following lemma shows that the $\T$-weighted diameter is (as one would hope) invariant under translations. This property follows easily from~\eqref{def:Lambda}. 

\pagebreak

\begin{lemma}\label{cor:weighted-diam-any-x}
Let $\T \subset \SS^{d-1}$ be a finite $\SS^{d-1}$-bounding set, let $\X \subset \R^d$ be bounded and non-empty, and let $x \in \R^d$. Then $\wdiam_\T(\X) = \wdiam_\T(\X + x)$.
\end{lemma}

\begin{proof}
Define $f(x) := \wdiam_\T(\X + x)$, and observe that
$$f(x) = \sum_{v \in \T} \lambda_v \big( d_v(\X) + \< x,v \> \big),$$
where $d_v(\X) = \sup_{y \in \X} \< y, v \>$ and $\Xi_\T = \{ \lambda_v : v \in \T \}$, as in Definition~\ref{def:weighted-diam}. Now, recalling that $\sum_{v \in \T} \lambda_v v = \0$, by~\eqref{def:Lambda}, it follows that $f(x) = f(\0)$, as claimed.
\end{proof}

The next useful property of the $\T$-weighted diameter is that $\wdiam_u(\X) = \Theta\big( \diam(\X) \big)$ for every non-leaf vertex $u$ of $T$, where the implicit constants depend only on $T$.  

\begin{lemma}\label{lem:wd-prop-diam}
There exists\/ $C > 0$ depending only on\/ $T$ such that the following holds. If\/~$u$~is a non-leaf vertex of\/ $T$, then
$$\frac{1}{C} \cdot \diam(\X) \le \wdiam_u(\X) \le \diam(\X)$$
for every bounded and non-empty set $\X \subset \R^d$. 
\end{lemma}

\begin{proof}
By Lemma~\ref{cor:weighted-diam-any-x}, we may assume without loss of generality that $\0 \in \X$. It follows that $d_v(\X) = \sup_{y \in \X} \< y, v \> \le \diam(\X)$, and hence 
\[
\wdiam_u(\X) = \sum_{v \in M_T(u)} \lambda_v \cdot d_v(\X) \leq \bigg( \sum_{v \in M_T(u)} \lambda_v \bigg) \cdot \diam(\X) = \diam(\X),
\]
since $\lambda_v > 0$ for all $v \in M_T(u)$ and $\sum_{v \in M_T(u)} \lambda_v = 1$, by~\eqref{def:Lambda}.

To prove the lower bound on $\wdiam_u(\X)$, observe first that 
\begin{equation}\label{eq:wd:lower:bound}
\wdiam_u(\X) = \sum_{v \in M_T(u)} \lambda_v \cdot d_v(\X) \ge \Delta_u(\X) \cdot \min \Xi_u.
\end{equation}
where $\Delta_u(\X) = \max\{ d_v(\X) : v \in M_T(u) \}$, and that
\begin{equation}\label{eq:X:subset:dcube}
\X \subset \bigcap_{v \in M_T(u)} \big\{ x \in \R^d : \< x,v \> \le \Delta_u(\X) \big\}.
\end{equation}
Note that $\Delta_u(\X) \ge 0$, since $\0 \in \X$, and set
\[
F_u := \bigcap_{v \in M_T(u)} \big\{ x \in \R^d \,:\, \< x,v \> \leq 1 \big\}.
\]
It follows from~\eqref{eq:wd:lower:bound} and~\eqref{eq:X:subset:dcube} that
\[
\diam(\X) \le \Delta_u(\X) \cdot \diam(F_u) \leq \frac{\diam(F_u)}{\min \Xi_u} \cdot \wdiam_u(\X) \leq C \cdot \wdiam_u(\X),
\]
where
\[
C := \max_{u \in T} \frac{\diam(F_u)}{\min \Xi_u}
\]
is finite, by Lemma~\ref{lem:icebergs-finite}, and depends only on $T$, as required.
\end{proof}

Let us note here, for future reference, the following immediate corollary of Lemma~\ref{lem:wd-prop-diam}. 

\begin{lemma}\label{lem:wd-prop1}
There exists\/ $C > 0$ depending only on\/ $T$ such that the following holds. If\/~$u$~and\/ $v$ are non-leaf vertices of\/ $T$, then
$$\wdiam_u(\X) \le C \cdot \wdiam_v(\X)$$
for every bounded and non-empty set $\X \subset \R^d$. 
\end{lemma}

\begin{proof}
By two applications of Lemma~\ref{lem:wd-prop-diam}, we have
\[
\wdiam_u(\X) \leq \diam(\X) \leq C \cdot \wdiam_v(\X),
\]
which proves the lemma.
\end{proof}

In the proofs of the sub-additivity lemmas it will be convenient to work with the following continuous analogue of droplets, which we shall refer to as `regions'.

\begin{definition}\label{def:FuGu}
Let $\T \subset \SS^{d-1}$ be a finite and non-empty set. A \emph{$\T$-region} is a non-empty set of the form
$$\bigcap_{u \in \T} \big\{ x \in \R^d : \< x,u \> \le a_u \big\}$$
for some collection $\{ a_u \in \R \,:\, u \in \T \}$.

Now, given a non-leaf vertex $u$ of $T$ and a bounded, non-empty set $\X \subset \R^d$, define 
\begin{enumerate}
\item $F_u(\X)$ to be the minimal $M_T(u)$-region containing $\X$;\smallskip
\item $G_\0(\X) := F_\0(\X)$, and if $u$ is not the root then
$$G_u(\X) := \big\{ x \in F^*_u(\X) : \< x,u\> \ge 0 \big\},$$
where $F^*_u(\X)$ is the minimal $N_T^\to(u)$-region containing $\X$.
\end{enumerate}
\end{definition}

Note that if $\{ x \in \X : \< x,u \> \ge 0 \}$ is non-empty then both $F_u(\X)$ and $G_u(\X)$ are bounded $M_T(u)$-regions. Observe also that if $\X \cap H_T(u) = \emptyset$ and $J$ is the minimal $u$-iceberg droplet containing $\X$, then by Definition~\ref{def:iceberg} we have $J \subset G_u(\X)$.

The following is our first sub-additivity lemma for $d$-dimensional droplets. Recalling from~\eqref{eq:diam-subadd-false} that the corresponding statement for $u = \0$ with $\wdiam_\0$ replaced by $\diam$ is false in general, the lemma provides further motivation for Definition~\ref{def:weighted-diam}.

\begin{lemma}\label{lem:sub-add-droplets}
Let $u$ be a non-leaf vertex of $T$, and let $\X_1,\X_2 \subset \R^d$ be bounded sets such that $\X_1 \cap \X_2 \neq \emptyset$. Then
\begin{equation}\label{eq:sub-add-droplets}
\wdiam_u\big( F_u( \X_1 \cup \X_2 ) \big) \leq \wdiam_u(\X_1) + \wdiam_u(\X_2).
\end{equation}
\end{lemma}

\begin{proof}
By Lemma~\ref{cor:weighted-diam-any-x}, we may assume that $\0 \in \X_1 \cap \X_2$, and hence $d_v(\X_1) \ge 0$ and $d_v(\X_2) \ge 0$ for each $v \in M_T(u)$. Let $\X := F_u( \X_1 \cup \X_2 )$, and observe that 
$$d_v(\X) = \sup_{x \in \X} \, \< x, v \> = \max\big\{ d_v(\X_1), d_v(\X_2) \big\} \le d_v(\X_1) + d_v(\X_2),$$
for each $v \in M_T(u)$. By Definition~\ref{def:weighted-diam}, it follows that
$$\wdiam_u(\X) 
\le \sum_{v \in M_T(u)} \lambda_v \cdot \big( d_v(\X_1) + d_v(\X_2) \big) = \wdiam_u(\X_1) + \wdiam_u(\X_2),$$
as required.
\end{proof}

\begin{figure}[ht]
  \centering
  \begin{tikzpicture}[>=latex] 
	\draw (-8.5,0) -- (0,0) -- (5.4,2.7);
	\draw (-1.7,0.4) -- (-0.9,2) -- (-2.5,2) -- cycle;
	\draw (-0.5,0) -- (-3.7,1.6) -- (-6.9,0);
	\draw [densely dashed] (-3.7,1.6) -- (-1.9,2.5) -- (3.1,0) -- (0,0);
	\draw [densely dashed] (-8.5,-2.4) -- (-6.9,0) -- (-5.9,1.5);
	\draw [densely dashed] (-5.9,1.5) -- (-2.5,2) -- (3.1,2.8);
	\draw [densely dashed] (0,0) -- (-4,-2);
	\node at (-3.7,0.6) {$J_1$};
	\node at (-1.7,1.5) {$J_2$};
	\node at (1,-1) {$\H_{u_1}^d \cup \H_{u_2}^d$};
	\draw [->] (-7.8,0) -- (-7.8,0.5) node [above] {$u_1$};
	\draw [->] (4.6,2.3) -- (4.4,2.7) node [above left] {$u_2$};
	\draw (0.5,1.3) -- (0.5,3) node [above] {$G_{u_1}(J_1 \cup J_2)$};
	\draw (-4.6,1.7) -- (-4.6,3) node [above] {$G_{u_2}(J_1 \cup J_2)$};
  \end{tikzpicture}
  \caption{A $u_1$-iceberg droplet $J_1$ and a $\0$-iceberg droplet $J_2$. Here, $\0 \to u_1 \to u_2$ is a path in $T$. Also shown are the regions $G_{u_1}(J_1 \cup J_2)$ and $G_{u_2}(J_1 \cup J_2)$ (the latter being clipped at its extremities for space).}
  \label{fig:subadd}
\end{figure}
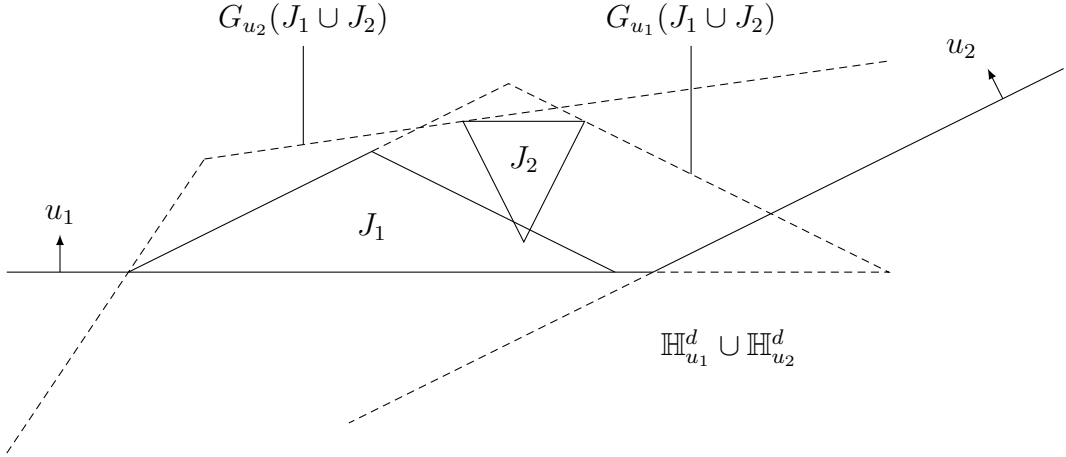

We shall use Lemma~\ref{lem:sub-add-droplets} in the (significantly more difficult) proof of our main sub-additivity lemma for icebergs, Lemma~\ref{lem:subadd}. We shall also use the following bound on the $M_T(u)$-weighted diameter of $F_v(\X)$, which is an easy consequence of Lemma~\ref{lem:wd-prop1}. 

\begin{lemma}\label{lem:wd-prop2}
There exists\/ $C > 0$ depending only on\/ $T$ such that the following holds. If\/~$u$~and\/ $v$ are non-leaf vertices of\/ $T$, then
$$\wdiam_u\big( F_v(\X) \big) \leq C \cdot \wdiam_u(\X)$$
for every bounded and non-empty set $\X \subset \R^d$. 
\end{lemma}

\begin{proof}
Note that $F_v(\X)$ is bounded, since $M_T(v)$ is an $\SS^{d-1}$-bounding set, by Lemma~\ref{lem:icebergs-finite}. Therefore, by Lemma~\ref{lem:wd-prop1}, we have
$$\wdiam_u\big( F_v(\X) \big) \le C \cdot \wdiam_v\big( F_v(\X) \big).$$
Now, by Definitions~\ref{def:weighted-diam} and~\ref{def:FuGu}, we have $\wdiam_v\big( F_v(\X) \big) = \wdiam_v(\X)$, since $F_v(\X)$ consists of exactly the points $y \in \R^d$ such that $\< y, v \> \le d_v(\X)$ for every $v \in M_T(v)$. Hence, applying Lemma~\ref{lem:wd-prop1} again, we obtain
$$\wdiam_v\big( F_v(\X) \big) = \wdiam_v(\X) \leq C \cdot \wdiam_u(\X),$$
as required.
\end{proof}

The following property of the function $G_u$, which follows from Lemma~\ref{lem:iceberg-drop-closed}, will also be useful in the proof of Lemma~\ref{lem:subadd}. 


\begin{lemma}\label{lem:Gu-closure}
Let $u$ be a non-leaf vertex of $T$, and let $K \subset \Lambda \setminus H_T(u)$ be a finite set. Then $G_u\big( [K]_u \big) = G_u(K)$, and in particular $[K]_u \subset G_u(K)$. 
\end{lemma}

\begin{proof}
By Definition~\ref{def:FuGu}, we are required to show that
\begin{equation}\label{eq:Gu-closure-suff}
\max\big\{ \< x,v \> \,:\, x \in [K]_u \big\} = \max\big\{ \< x,v \> \,:\, x \in K \big\}
\end{equation}
for each $v \in N_T^\to(u)$. 
To do so, let $J$ be the minimal $u$-iceberg droplet containing $K$, and (recalling Definition~\ref{def:iceberg}) observe that
$$\max\big\{ \< x,v \> : x \in J \big\} = \max\big\{ \< x,v \> : x \in K \big\}$$
for every $v \in N_T^\to(u)$. 
Now, since $J$ is $u$-closed, by Lemma~\ref{lem:iceberg-drop-closed}, and $K \subset J$, we have
$$\max\big\{ \< x,v \> : x \in [K]_u \big\} \le \max\big\{ \< x,v \> : x \in [J]_u \big\} = \max\big\{ \< x,v \> : x \in J \big\},$$
and~\eqref{eq:Gu-closure-suff} follows, since $K \subset [K]_u$. This proves that $G_u\big( [K]_u \big) = G_u(K)$, and it follows that $[K]_u \subset G_u(K)$, since $\< x,u\> \ge 0$ for every $x \in [K]_u$. 
\end{proof}

We need two more technical lemmas for the proof of Lemma~\ref{lem:subadd}. The following preliminary inequality will be used in the proof of each of them. We shall write $O(1)$ to denote a function that is bounded by a constant depending on both $T$ and $R$. 


\begin{lemma}\label{lem:wd-close-to-H}
Let\/ $u$ be a non-root, non-leaf vertex of\/ $T$, let\/ $\X \subset \R^d$ be a bounded and non-empty set, and suppose that\/ $\X$ is strongly connected to\/ $\HH_u^d$. Then
\[
\wdiam_u\big( G_u(\X) \big) \le \wdiam_u(\X) + O(1).
\]
\end{lemma}

\begin{proof}
Recall from~\eqref{eq:MTu} that 
$-u \in M_T(u)$. By Definition~\ref{def:FuGu}, and since $\X$ is strongly connected to $\HH_u^d$, it follows that
\[
d_{-u}\big( F_u(\X) \big) = d_{-u}(\X) \ge -R,
\]
where $d_v(\X) = \sup_{y \in \X} \< y, v \>$, as in Definition~\ref{def:weighted-diam}. Hence, since $d_v\big( G_u(\X) \big) \le d_v\big( F_u(\X) \big)$ for all $v \in N_T^\to(u)$, also by Definition~\ref{def:FuGu}, it follows from Definition~\ref{def:weighted-diam} that
\[
\wdiam_u\big( G_u(\X) \big) \leq \wdiam_u\big( F_u(\X) \big) + O(1).
\]
But $\wdiam_u\big( F_u(\X) \big) = \wdiam_u(\X)$, so this completes the proof.
\end{proof}

The first of our two technical lemmas allows us to control the change in $\wdiam_v(G_v(K))$ as $v$ varies along a $(K,u)$-permissible sequence. 

\begin{lemma}\label{lem:wd-permissible}
There exists $C > 0$ depending only on $T$ such that the following holds. Let $\0 \to u_1 \to \cdots \to u_k$ be a path in $T$, with $u_k$ not a leaf, let $K \subset \Lambda \setminus H_T(u_k)$ be a finite, non-empty set, and let $0 = i_0 < \dots < i_\ell \leq k$ be a $(K,u_k)$-permissible sequence. Then
\begin{equation}\label{eq:wd-permissible1}
\wdiam_{u_{i_j}}(G_{i_j}) \leq C \cdot \wdiam_{u_{i_{j+1}}}(G_{i_{j+1}}) + O(1)
\end{equation}
and
\begin{equation}\label{eq:wd-permissible2}
\wdiam_{u_{i_{j+1}}}(G_{i_{j+1}}) \leq C \cdot \wdiam_{u_{i_j}}(G_{i_j}) + O(1)
\end{equation}
for each $0 \leq j < \ell$, where $G_i := G_{u_i}(K)$ for each $0 \leq i \le k$. 
\end{lemma}

\begin{proof}
Starting with~\eqref{eq:wd-permissible2}, the key fact we need is that $J_{i_j}$, the minimal $u_{i_j}$-iceberg droplet containing $K$, is strongly connected to $\HH_{u_{i_{j+1}}}^d$, which follows from Definition~\ref{def:iceberg-container} since $0 = i_0 < \dots < i_\ell \leq k$ is a $(K,u_k)$-permissible sequence. Since $J_{i_j} \subset G_{i_j}$, by Definitions~\ref{def:iceberg} and~\ref{def:FuGu}, this means that $G_{i_j}$ is also strongly connected to $\HH_{u_{i_{j+1}}}^d$. 

By Lemma~\ref{lem:wd-close-to-H} applied with $\X = G_{i_j}$, it follows that
\begin{equation}\label{eq:wd-perm1}
\wdiam_{u_{i_{j+1}}}\big( G_{u_{i_{j+1}}}(G_{i_j}) \big) \leq \wdiam_{u_{i_{j+1}}}(G_{i_j}) + O(1).
\end{equation}
Now, note that $K \subset J_{i_j} \subset G_{i_j}$, and so
\begin{equation}\label{eq:wd-perm2}
\wdiam_{u_{i_{j+1}}}(G_{i_{j+1}}) = \wdiam_{u_{i_{j+1}}}\big( G_{u_{i_{j+1}}}(K) \big) \leq \wdiam_{u_{i_{j+1}}}\big( G_{u_{i_{j+1}}}(G_{i_j}) \big).
\end{equation}
Moreover, by Lemma~\ref{lem:wd-prop1} we have
\begin{equation}\label{eq:wd-perm3}
\wdiam_{u_{i_{j+1}}}(G_{i_j}) \leq C \cdot \wdiam_{u_{i_j}}\big( G_{i_j} \big),
\end{equation}
for some $C > 0$ depending only on $T$. Combining~\eqref{eq:wd-perm1},~\eqref{eq:wd-perm2} and~\eqref{eq:wd-perm3} gives~\eqref{eq:wd-permissible2}.

To prove~\eqref{eq:wd-permissible1}, we shall compare $\wdiam_{u_{i_{j+1}}}(G_{i_{j+1}})$ to $\wdiam_{u_0}(G_0)$, and then use the fact that we have already proved~\eqref{eq:wd-permissible2} for all $0 \leq j < \ell$. To be precise, observe first that
\begin{equation}\label{eq:wd-permissible1:claim}
\wdiam_{u_{i_{j+1}}}(G_{i_{j+1}}) \geq \wdiam_{u_{i_{j+1}}}(K) \geq \frac{1}{C'} \cdot \wdiam_{\0}(K) = \frac{1}{C'} \cdot \wdiam_{\0}(G_0),
\end{equation}
for some $C' > 0$ depending only on $T$. Indeed, the first inequality holds because (as noted above) $K \subset G_{i_{j+1}}$, and the second inequality holds by Lemma~\ref{lem:wd-prop1}. For the final step, recall from Definition~\ref{def:FuGu} that $G_0 = G_\0(K)$ is the minimal $N_T^\to(\0)$-region containing $K$, and that $M_T(\0) = N_T^\to(\0)$, so the equality holds by Definition~\ref{def:weighted-diam}. Now, by~\eqref{eq:wd-permissible2} (applied $j$ times), we have
\begin{equation}\label{eq:wd-permissible1:using2}
\wdiam_{u_{i_j}}(G_{i_j}) \le C^j \cdot \wdiam_{\0}(G_0) + O(1),
\end{equation}
and combining~\eqref{eq:wd-permissible1:claim} and~\eqref{eq:wd-permissible1:using2} gives~\eqref{eq:wd-permissible1}, as required.
\end{proof}

The final lemma that we need for the proof of Lemma~\ref{lem:subadd} allows us to bound the `size' of the $u_k$-closure of a set in terms of the `size' of its $u_i$-closure.

\begin{lemma}\label{lem:wd-change-closures}
There exists $C > 0$ depending only on $T$ such that the following holds. Let $\0 \to u_1 \to \cdots \to u_k$ be a path in $T$, with $u_k$ not a leaf, and let $K \subset \Lambda \setminus H_T(u_k)$ be a finite, non-empty set. Then
\[
\wdiam_{u_i}\big( G_{u_i}\big( [K]_{u_k} \big) \big) \leq C \cdot \wdiam_{u_i}\big( G_{u_i}\big( [K]_{u_i} \big) \big) + O(1)
\]
for all $0 \leq i \leq k$.
\end{lemma}

\begin{proof}
We shall in fact prove that
\begin{equation}\label{eq:wd-change-closures-ind}
\wdiam_{u_i}\big( G_{u_i}\big( [K]_{u_{j+1}} \big) \big) \leq C \cdot \wdiam_{u_i}\big( G_{u_i}\big( [K]_{u_j} \big) \big) + O(1)
\end{equation}
for each $i \leq j < k$, which is clearly sufficient to prove the lemma.

If $[K]_{u_{j+1}} = [K]_{u_j}$ then there is nothing to prove, so assume that $[K]_{u_{j+1}} \ne [K]_{u_j}$. We claim that in this case $[K]_{u_j}$ is strongly connected to $\HH_{u_{j+1}}^d$. To see this, note first that, by Definition~\ref{def:u-closure}, if $[K]_{u_j}$ does not intersect $\HH_{u_{j+1}}^d$, then $[K]_{u_j} \subset [K]_{u_{j+1}}$. It follows that if $[K]_{u_{j+1}} \ne [K]_{u_j}$, then $x + X \subset [K]_{u_j} \cup H_T(u_{j+1})$ for some $x \not\in [K]_{u_j} \cup H_T(u_{j+1})$ and $X \in \U$. Now, the sets $H_T(u_{j+1})$ and $[K]_{u_j} \cup H_T(u_j)$ are $\U$-closed, by Proposition~\ref{lem:trees-exist} and Definition~\ref{def:u-closure}, so $x + X$ must intersect both $[K]_{u_j}$ and $\H_{u_{j+1}}^d$. By Definition~\ref{def:strconn}, this implies that $[K]_{u_j}$ is strongly connected to $\HH_{u_{j+1}}^d$, as claimed.

Now, observe that $K \subset G_{u_i}\big( [K]_{u_j} \big)$, and therefore $[K]_{u_{j+1}} \subset G_{u_{j+1}}\big( G_{u_i}\big( [K]_{u_j} \big) \big)$, by Lemma~\ref{lem:Gu-closure}, so
\begin{equation}\label{eq:wd-change-closures-1}
G_{u_i}\big( [K]_{u_{j+1}} \big) \subset G_{u_i}\big( G_{u_{j+1}}\big( G_{u_i}( [K]_{u_j} ) \big) \big) =: \X.
\end{equation}
In order to bound $\wdiam_{u_i}(\X)$, we claim that, for each $v \in N_T(u_i)$, we have
\begin{equation}\label{eq:wd-change-closures-dbound}
d_v(\X) \le C \cdot \wdiam_{u_i}\big( G_{u_i}( [K]_{u_j} ) \big) + O(1).
\end{equation}
To prove~\eqref{eq:wd-change-closures-dbound}, let us write $\Y := G_{u_{j+1}}\big( G_{u_i}( [K]_{u_j} ) \big)$, so that $\X = G_{u_i}(\Y)$. Now, observe that $d_v(\X) \le d_v(\Y)$, since $\X$ is contained in the minimal $N_T^\to(u_i)$-region containing $\Y$, by Definition~\ref{def:FuGu}. It follows that
\begin{equation}\label{eq:wd-change-closures-2}
d_v(\X) \le d_v(\Y) \le C \cdot \wdiam_{u_i}(\Y) \le C^2 \cdot \wdiam_{u_{j+1}}(\Y),
\end{equation}
by Definition~\ref{def:weighted-diam} and Lemma~\ref{lem:wd-prop1}, for some $C$ depending only on $T$. 

Now, since 
$[K]_{u_j}$ is strongly connected to $\HH_{u_{j+1}}^d$ and $i \le j$, so $\< x,u_i\> \ge 0$ for every $x \in [K]_{u_j}$, it follows that $G_{u_i}( [K]_{u_j} )$ is also strongly connected to $\HH_{u_{j+1}}^d$, and thus
$$\wdiam_{u_{j+1}}\big( G_{u_{j+1}}\big( G_{u_i}( [K]_{u_j} ) \big) \big) \le \wdiam_{u_{j+1}}\big( G_{u_i}( [K]_{u_j} ) \big) + O(1),$$
by Lemma~\ref{lem:wd-close-to-H}. Hence, by Lemma~\ref{lem:wd-prop1}, we obtain
$$\wdiam_{u_{j+1}}(\Y) = \wdiam_{u_{j+1}}\big( G_{u_{j+1}}\big( G_{u_i}( [K]_{u_j} ) \big) \big) \le C \cdot \wdiam_{u_i}\big( G_{u_i}( [K]_{u_j} ) \big) + O(1),$$
which, together with~\eqref{eq:wd-change-closures-2}, proves~\eqref{eq:wd-change-closures-dbound}.

If $u_i = \0$ then~\eqref{eq:wd-change-closures-ind} now follows from~\eqref{eq:wd-change-closures-1} and~\eqref{eq:wd-change-closures-dbound}. On the other hand, if $u_i \ne \0$ then observe that
$$d_{-u_i}(\X) = d_{-u_i}\big( G_{u_i}( \Y ) \big) \le 0,$$ 
since $\< x,u_i\> \ge 0$ for every $x \in G_{u_i}(\Y)$, by Definition~\ref{def:FuGu}. 
Together with~\eqref{eq:wd-change-closures-1} and~\eqref{eq:wd-change-closures-dbound}, this implies~\eqref{eq:wd-change-closures-ind}, and hence completes the proof of the lemma. 
\end{proof}

We are almost ready to state our sub-additivity lemma for icebergs. In order to do so, we need one more important definition. Let $\delta > 0$ be a sufficiently small constant, depending on $T$, and fixed for the rest of the paper. 

\begin{definition}\label{def:diamstar}
Let $\J = (J,u)$ be an iceberg and let $k := d_T(\0,u)$. Then we define $\diam^*(\J)$, the \emph{iceberg diameter} of $\J$, as follows:
\begin{equation}\label{eq:diamstar}
\diam^*(\J) := \delta^k \wdiam_u\big( G_u(J) \big).
\end{equation}
\end{definition}

The next lemma is our sub-additivity lemma for icebergs. The lemma will be used in the proofs of both the Aizenman--Lebowitz lemma for icebergs (Lemma~\ref{lem:AL}) and the extremal lemma for icebergs (Lemma~\ref{lem:extremal}), and also in the proof of Lemma~\ref{lem:goodsat}.


\begin{lemma}\label{lem:subadd}
Let $u$ be a non-leaf vertex of\/ $T$, and let $K_1$ and $K_2$ be finite non-empty subsets of\/ $\Lambda \setminus H_T(u)$ such that\/ $[K_1]_u$, $[K_2]_u$ and\/ $[K_1]_u \cup [K_2]_u$ are all strongly connected. Then
\begin{equation}\label{eq:subadd}
\diam^*\Big( \J_u\big( [K_1 \cup K_2]_u \big) \Big) \leq \diam^*\Big( \J_u\big( [K_1]_u \big) \Big) + \diam^*\Big( \J_u\big( [K_2]_u \big) \Big) + O(1).
\end{equation}
\end{lemma}

\begin{proof}
Let $\0 \to u_1 \to \cdots \to u_k = u$ be the path from $\0$ to $u$ in $T$. Let
\begin{equation}\label{def:J1J2:containers}
\J_1 := \J_u\big( [K_1]_u \big) = (J_1,u_{i_1}) \quad \text{and} \quad \J_2 := \J_u\big( [K_2]_u \big) = (J_2, u_{i_2}) 
\end{equation}
for some $0 \leq i_1,i_2 \leq k$, and assume without loss of generality that $i_1 \geq i_2$. Now let
\[
\J = \J_u\big([K_1 \cup K_2]_u\big) = (J,u_\ell),
\]
and observe, by Lemma~\ref{lem:iceberg-cont-mono-type}, that $\ell \geq \max\{i_1,i_2\}$, since $[K_1]_u \cup [K_2]_u \subset [K_1 \cup K_2]_u$. We~claim that we may in fact assume that $\ell = k$. Indeed, this follows from Lemma~\ref{lem:iceberg-cont-reduce}, which implies that
$$\J_{u_\ell}\big( [K_1 \cup K_2]_{u_\ell} \big) = \J_{u_k}\big( [K_1 \cup K_2]_{u_k} \big),$$
since $\J$ is a $u_\ell$-iceberg, and similarly that
$$\J_{u_\ell}\big( [K_j]_{u_\ell} \big) = \J_{u_k}\big( [K_j]_{u_k} \big) \qquad \text{and} \qquad [K_j]_{u_\ell} = [K_j]_{u_k}$$
for each $j \in \{1,2\}$, so $[K_1]_{u_\ell}$, $[K_2]_{u_\ell}$ and $[K_1]_{u_\ell} \cup [K_2]_{u_\ell}$ are all strongly connected.

Recalling Definition~\ref{def:diamstar}, it follows that~\eqref{eq:subadd} is equivalent to
\begin{equation}\label{eq:subadd-step1}
\delta^k \wdiam_{u_k}\big( G_{u_k}( J ) \big) \le \delta^{i_1} \wdiam_{u_{i_1}}\big( G_{u_{i_1}}( J_1 ) \big) + \delta^{i_2} \wdiam_{u_{i_2}}\big( G_{u_{i_2}}( J_2 ) \big) + O(1).
\end{equation}
We divide the proof of~\eqref{eq:subadd-step1} into two claims, using the set $[K_1 \cup K_2]_{u_{i_1}}$ as a stepping stone 
between the sets $[K_1]_{u_{i_1}}$ and $[K_2]_{u_{i_2}}$, and the set $[K_1 \cup K_2]_{u_k}$.

\begin{claim}\label{clm:subadd1}
$$\delta^{i_1} \wdiam_{u_{i_1}}\big( G_{u_{i_1}}\big( [K_1 \cup K_2]_{u_{i_1}} \big) \big) 
\leq \delta^{i_1} \wdiam_{u_{i_1}}\big( G_{u_{i_1}}(J_1) \big) + \delta^{i_2} \wdiam_{u_{i_2}}\big( G_{u_{i_2}}(J_2) \big) + O(1).$$
\end{claim}

\smallskip

\begin{claim}\label{clm:subadd2}
$$\delta^k \wdiam_{u_k}\big( G_{u_k}(J) \big) \leq \delta^{i_1} \wdiam_{u_{i_1}}\big( G_{u_{i_1}}\big( [K_1 \cup K_2]_{u_{i_1}} \big) \big) + O(1).$$
\end{claim}

Note that~\eqref{eq:subadd-step1} follows immediately from Claims~\ref{clm:subadd1} and~\ref{clm:subadd2}, and therefore, by the observations above, proving these two claims is sufficient to prove Lemma~\ref{lem:subadd}.  

\begin{clmproof}{clm:subadd1}
Recalling Definition~\ref{def:FuGu}, set 
\[
G_1 := G_{u_{i_1}}(J_1), \quad G_2 := G_{u_{i_2}}(J_2), \quad \text{and} \quad G_2' := F_{u_{i_1}}(G_2) = F_{u_{i_1}}\big( G_{u_{i_2}}(J_2) \big).
\]
Our first sub-claim is that
\begin{equation}\label{eq:subadd1-subclm1a}
\delta^{i_2} \wdiam_{u_{i_2}}( G_2 ) \ge \delta^{i_1} \wdiam_{u_{i_1}}( G_2' ).
\end{equation}
If $i_1 = i_2$ then $G_2' = G_2$ because in this case $G_2$ is already an $M_T(u_{i_1})$-region, and~\eqref{eq:subadd1-subclm1a} is immediate. On the other hand, if $i_1 > i_2$ then we have
\[
\delta^{i_2} \wdiam_{u_{i_2}}( G_2 ) \ge \delta^{i_2 + 1/2} \wdiam_{u_{i_2}}( G_2' ) \ge \delta^{i_2 + 1} \wdiam_{u_{i_1}}( G_2' ) \ge \delta^{i_1} \wdiam_{u_{i_1}}( G_2' ),
\]
where the first inequality follows from Lemma~\ref{lem:wd-prop2} and the second from Lemma~\ref{lem:wd-prop1}, and since $\delta$ is sufficiently small (depending on $T$). Thus in either case~\eqref{eq:subadd1-subclm1a} holds, and hence 
\begin{equation}\label{eq:subadd1-subclm1}
\delta^{i_1} \wdiam_{u_{i_1}}( G_1 ) + \delta^{i_2} \wdiam_{u_{i_2}}( G_2 ) \geq \delta^{i_1}\big( \wdiam_{u_{i_1}}( G_1 ) + \wdiam_{u_{i_1}}( G_2' ) \big).
\end{equation}

Our second sub-claim is that
\begin{equation}\label{eq:subadd1-subclm2}
\wdiam_{u_{i_1}}( G_1 ) + \wdiam_{u_{i_1}}( G_2' ) \geq \wdiam_{u_{i_1}}\big( F_{u_{i_1}}( G_1 \cup G_2') \big) - O(1).
\end{equation}
This will follow from two applications of Lemma~\ref{lem:sub-add-droplets}, and our assumption that the sets $[K_1]_u$ and $[K_2]_u$ are strongly connected to each other. Observe first that $[K_1]_u \subset J_1 \subset G_1$, since $J_1$ is, by~\eqref{def:J1J2:containers} and Definition~\ref{def:iceberg-container}, the minimal $u_{i_1}$-iceberg droplet containing $[K_1]_u$, and therefore does not intersect $\H_{u_{i_1}}^d$. Similarly, we have $[K_2]_u \subset J_2 \subset G_2 \subset G_2'$, and so $G_1$ and $G_2'$ are strongly connected to each other, since $[K_1]_u$ and $[K_2]_u$ are strongly connected to each other by assumption. Let $B$ be a set of diameter at most $R$ such that $G_1 \cap B \neq \emptyset$ and $G_2' \cap B \neq \emptyset$. Applying Lemma~\ref{lem:sub-add-droplets} to the sets $G_1$ and $B$ gives
$$\wdiam_{u_{i_1}}\big( F_{u_{i_1}}( G_1 \cup B ) \big) \le \wdiam_{u_{i_1}}(G_1) + \wdiam_{u_{i_1}}(B),$$
and applying Lemma~\ref{lem:sub-add-droplets} to the sets $F_{u_{i_1}}(G_1 \cup B)$ and $G_2'$ gives 
$$\wdiam_{u_{i_1}}\big( F_{u_{i_1}}( G_1 \cup G'_2 ) \big) \le \wdiam_{u_{i_1}}\big( F_{u_{i_1}}( G_1 \cup B ) \big) + \wdiam_{u_{i_1}}(G'_2),$$
since $G_1 \cup G'_2 \subset F_{u_{i_1}}( G_1 \cup B ) \cup G_2'$. Combining these two inequalities, and noting that $\wdiam_{u_{i_1}}(B) \le R$, by Lemma~\ref{lem:wd-prop-diam}, we obtain the inequality~\eqref{eq:subadd1-subclm2}.

Our third sub-claim, which will complete the proof of Claim~\ref{clm:subadd1}, is that
\begin{equation}\label{eq:subadd1-subclm3-pre}
G_{u_{i_1}}\big( [K_1 \cup K_2]_{u_{i_1}} \big) \subset F_{u_{i_1}}\big( G_1 \cup G_2' \big),
\end{equation}
and hence that
\begin{equation}\label{eq:subadd1-subclm3}
\wdiam_{u_{i_1}} \Big( G_{u_{i_1}}\big( [K_1 \cup K_2]_{u_{i_1}} \big) \Big) \leq \wdiam_{u_{i_1}} \Big( F_{u_{i_1}}\big( G_1 \cup G_2' \big) \Big).
\end{equation}
To prove~\eqref{eq:subadd1-subclm3-pre}, first we show that
\begin{equation}\label{eq:subadd1-subclm3-1}
G_{u_{i_1}}\big( G_1 \cup G_2' \big) \subset F_{u_{i_1}}\big( G_1 \cup G_2' \big).
\end{equation}
If $u_{i_1} = \0$, then $F_{u_{i_1}} = G_{u_{i_1}}$ and~\eqref{eq:subadd1-subclm3-1} holds trivially, so suppose that $u_{i_1}$ is not the root. Now, recalling Definition~\ref{def:FuGu}, observe that $G_1 = G_{u_{i_1}}(J_1)$ is bounded and $F^*_{u_{i_1}}(J_1)$ is unbounded, by Lemma~\ref{lem:icebergs-finite}, and that therefore
$$G_1 \cap \big\{ x \in \R^d \,:\, \< x,u_{i_1} \> = 0 \big\} \neq \emptyset.$$
By Definition~\ref{def:FuGu}, this is enough to prove~\eqref{eq:subadd1-subclm3-1}. 

To deduce~\eqref{eq:subadd1-subclm3-pre} from~\eqref{eq:subadd1-subclm3-1}, it suffices to observe that
\[
G_{u_{i_1}}\big( [K_1 \cup K_2]_{u_{i_1}} \big) = G_{u_{i_1}}\big( K_1 \cup K_2 \big) \subset G_{u_{i_1}}\big( G_1 \cup G_2' \big),
\]
where the first step follows from Lemma~\ref{lem:Gu-closure}, and the second step follows since (as we observed in the proof of the second sub-claim) we have $K_1 \subset G_1$ and $K_2 \subset G_2'$.

Claim~\ref{clm:subadd1} now follows by combining~\eqref{eq:subadd1-subclm1},~\eqref{eq:subadd1-subclm2}, and~\eqref{eq:subadd1-subclm3}.
\end{clmproof}

Claim~\ref{clm:subadd2} is a relatively straightforward application of Lemmas~\ref{lem:wd-permissible} and~\ref{lem:wd-change-closures}.

\begin{clmproof}{clm:subadd2}
Recall that $J$ is, by Definition~\ref{def:iceberg-container}, the minimal $u_k$-iceberg droplet containing $[K_1 \cup K_2]_{u_k}$, and by Definition~\ref{def:iceberg} is therefore contained in the minimal $N_T^\to(u_k)$-droplet containing $[K_1 \cup K_2]_{u_k}$. By Definition~\ref{def:FuGu}, it follows that $G_{u_k}(J) \subset G_{u_k}\big( [K_1 \cup K_2]_{u_k} \big)$, and hence
\begin{equation}\label{eq:subadd2-starter}
\wdiam_{u_k}\big( G_{u_k}(J) \big) \leq \wdiam_{u_k}\Big( G_{u_k}\big( [K_1 \cup K_2]_{u_k} \big) \Big).
\end{equation}
This proves the claim if $i_1 = k$, so let us assume that $i_1 < k$. Observe that
\begin{equation}\label{eq:subadd2-starter2}
\wdiam_{u_{i_1}}\Big( G_{u_{i_1}}\big( [K_1 \cup K_2]_{u_k} \big) \Big) \leq C \cdot \wdiam_{u_{i_1}}\Big( G_{u_{i_1}}\big( [K_1 \cup K_2]_{u_{i_1}} \big) \Big) + O(1),
\end{equation}
for some $C$ depending only on $T$, by Lemma~\ref{lem:wd-change-closures}. Now define
\[
G^{(j)} := G_{u_j}\big( [K_1 \cup K_2]_{u_k} \big)
\]
for each $0 \leq j \leq k$. Since we chose $\delta$ sufficiently small depending on $T$, and since $i_1 < k$, by~\eqref{eq:subadd2-starter} and~\eqref{eq:subadd2-starter2} it will suffice to show that
\begin{equation}\label{eq:subadd2:claim2:suff}
\delta^{1/2} \wdiam_{u_k}\big( G^{(k)} \big) \le \wdiam_{u_{i_1}}\big( G^{(i_1)} \big) + O(1).
\end{equation}
To prove~\eqref{eq:subadd2:claim2:suff} we shall apply Lemma~\ref{lem:wd-permissible} along two $([K_1 \cup K_2]_{u_k}, u_k)$-permissible sequences leading to $i_1$ and $k$ respectively. 

To be precise, by Definition~\ref{def:iceberg-container} and since $\J_{u_k}\big( [K_1 \cup K_2]_{u_k} \big) = (J,u_k)$, there exists a $([K_1 \cup K_2]_{u_k}, u_k)$-permissible sequence containing $k$. Moreover, since $\J_{u_k}\big( [K_1]_{u_k} \big) = (J_1,u_{i_1})$, there exists a $([K_1]_{u_k}, u_k)$-permissible sequence containing $i_1$. Now, recall that if $K \subset K'$ then any $(K,u_k)$-permissible sequence is also a $(K',u_k)$-permissible sequence (see the proof of Lemma~\ref{lem:iceberg-cont-mono-type}), so there also exists a $([K_1 \cup K_2]_{u_k}, u_k)$-permissible sequence containing $i_1$. By Lemma~\ref{lem:wd-permissible}, it follows that
\[
\wdiam_{u_k}\big( G^{(k)} \big) \leq C^k \wdiam_{\0}\big( G^{(0)} \big) + O(1) \leq C^{k + i_1} \wdiam_{u_{i_1}}\big( G^{(i_1)} \big) + O(1).
\]
Since we chose $\delta$ sufficiently small depending on $T$, this proves~\eqref{eq:subadd2:claim2:suff}, and hence, by~\eqref{eq:subadd2-starter} and~\eqref{eq:subadd2-starter2}, this suffices to prove Claim~\ref{clm:subadd2}.
\end{clmproof}

As observed earlier, Claims~\ref{clm:subadd1} and~\ref{clm:subadd2} together imply~\eqref{eq:subadd-step1}, and thus together they complete the proof of Lemma~\ref{lem:subadd}. 
\end{proof}

We shall next deduce our Aizenman--Lebowitz-type and extremal lemmas for icebergs from the $u$-iceberg spanning algorithm and Lemma~\ref{lem:subadd}. In order to do so, we need the following lemma, which provides a uniform upper bound on the size of the $u$-iceberg container of the $u$-closure of a \emph{single} vertex. The proof of the lemma is (perhaps surprisingly) non-trivial, and relies on Lemma~\ref{lem:wd-permissible}. 

\begin{lemma}\label{lem:extremal:basecase}
There exists a constant\/ $C > 0$ 
such that the following holds. If\/ $u$ is a non-leaf vertex of\/ $T$ and\/ $x \in \Lambda \setminus H_T(u)$, then\/ $\diam^*\big( \J_u\big( [\{x\}]_u \big) \big) \le C$.
\end{lemma}

\begin{proof}
Set $K := [\{x\}]_u$, and note that $K \subset \Lambda \setminus H_T(u)$, by Definition~\ref{def:u-closure}, and that $K$ is non-empty and finite, by Lemmas~\ref{cor:icebergs-finite} and~\ref{lem:iceberg-drop-closed}, since it is contained in the minimal $u$-iceberg droplet containing $x$. Let $\J_u(K) = (J,v)$ so, by Definition~\ref{def:diamstar}, we are required to bound $\wdiam_v\big( G_v(J) \big)$. 

Observe that $G_v(J) = G_v(K)$, since $J$ is the minimal $v$-iceberg droplet containing $K$. Now, since by Definition~\ref{def:iceberg-container} there exists a $(K,u)$-permissible sequence ending in $v$, we have   
$$\wdiam_v\big( G_v(J) \big) = \wdiam_v\big( G_v(K) \big) \le C^k \wdiam_\0\big( G_\0(K) \big) + O(1),$$
by Lemma~\ref{lem:wd-permissible}. Since $G_\0(K) = F_\0(K)$ and $\wdiam_\0\big( F_\0(K) \big) = \wdiam_\0(K)$, the lemma will hold if we can prove a uniform bound on the $\0$-weighted diameter of $K = [\{x\}]_u$. 

To bound $\wdiam_\0(K)$, let $\J_u(\{x\}) = (J',w)$ be the $u$-iceberg container of $\{x\}$, and observe that $J'$ is $u$-closed, by Lemma~\ref{lem:iceberg-cont-closed}, and therefore that $K \subset J'$. Note also that $J' \subset G_w(J') = G_w(\{x\})$, since $J'$ is the minimal $w$-iceberg droplet containing $x$. Now, since there exists an $(\{x\},u)$-permissible sequence ending in $w$, we have   
$$\wdiam_w(K) \le \wdiam_w\big( G_w(\{x\}) \big) \le C^k \wdiam_\0\big( G_\0(\{x\}) \big) + O(1),$$
by Lemma~\ref{lem:wd-permissible}. Noting that $G_\0(\{x\}) = F_\0(\{x\})$ and $\wdiam_\0\big( F_\0(\{x\}) \big) = 0$, since $\{x\}$ is an $N_T^\to(\0)$-region, it follows, using Lemma~\ref{lem:wd-prop1}, that 
$$\wdiam_\0(K) \le C \cdot \wdiam_w(K) = O(1),$$
as required.
\end{proof}

We can now prove our Aizenman--Lebowitz-type lemma for icebergs. We remark that, while the lemma is a straightforward consequence of the $u$-iceberg spanning algorithm and Lemmas~\ref{lem:output-of-span},~\ref{lem:subadd} and~\ref{lem:extremal:basecase}, the statement is slightly more complicated than that of the corresponding lemma in~\cite{BDMS} (which was about $\T$-internally spanned $\T$-droplets). This is because, while an iceberg being $u$-iceberg spanned implies the existence of similar such icebergs at smaller scales, it does \emph{not} imply that the corresponding iceberg droplets are nested if the icebergs are of different types. However, the lemma gives us control over the location of the elements of $A$ that are used to span our icebergs, and this will suffice for our purposes. Let $\lambda > 0$ be a sufficiently large constant (in particular, larger than the bound in Lemma~\ref{lem:extremal:basecase}, and the constant implicit in the $O(1)$-term in Lemma~\ref{lem:subadd}). 

\begin{lemma}[Aizenman--Lebowitz lemma for icebergs]\label{lem:AL}
Let $u$ be a non-leaf vertex of\/~$T$, and let $K \subset \Lambda \setminus H_T(u)$ be a finite set. Suppose that $\< K \>_u = \{ \J \}$ for some iceberg $\J$, and let\/ $\lambda \leq m \leq \diam^*(\J)$. Then there exists a set $K' \subset K$ and an iceberg $\J'$, with 
$$m \leq \diam^*(\J') \leq 3m,$$
such that $\< K' \>_u = \{ \J' \}$.
\end{lemma}

\begin{proof}
Apply the $u$-iceberg spanning algorithm (Definition~\ref{def:icespanalg}) to $K$. We claim that
$$f(t) := \max\Big\{ \diam^*\big( \J_u\big( [L]_u \big) \big) : L \in \K^t \Big\}$$
satisfies $f(t) \le 2f(t-1) + O(1)$ for every $t \in [T]$, and hence that $f(t) \le 3f(t-1)$ if $f(t-1) \ge \lambda$, by our choice of $\lambda$. To see this, suppose that the sets $K_1$ and $K_2$ are united in step $t$ of the algorithm, and observe that $K_1 \cup K_2 \subset K \subset \Lambda \setminus H_T(u)$, and that the sets $[K_1]_u$, $[K_2]_u$ and $[K_1]_u \cup [K_2]_u$ are all strongly connected (see the proof of Lemma~\ref{lem:output-of-span}). It follows, by Lemma~\ref{lem:subadd}, that
\begin{equation}\label{eq:subadd:AL:app}
\diam^*\big( \J_u\big( [K_1 \cup K_2]_u \big) \big) \leq \diam^*\big( \J_u\big( [K_1]_u \big) \big) + \diam^*\big( \J_u\big( [K_2]_u \big) \big) + O(1),
\end{equation}
so $f(t) \le 2f(t-1) + O(1)$, as claimed. 

Now, recall that $\K^0 = \big\{ \{x\} : x \in K \big\}$ and observe that, by Lemma~\ref{lem:extremal:basecase} and our choice of $\lambda$, we have $\diam^*\big( \J_u\big( [\{x\}]_u \big) \big) < \lambda$ for every $x \in K$. Since $\< K \>_u = \{ \J \}$ and $\lambda \le m \le \diam^*(\J)$, it follows that there exists a minimum $t \in [T]$ such that $f(t) \ge m$, that is, such that there exists a set  $K' \in \K^t$ with 
$$\diam^*\big( \J_u\big( [K']_u \big) \ge m.$$
Now, since $m \ge \lambda$, and hence $f(t) \le 3f(t-1)$, the minimality of $t$ implies that 
$$m \le \diam^*\big( \J_u\big( [K']_u \big) \big) \le 3m,$$
as required. Finally, note that we have $\< K' \>_u = \big\{ \J_u\big( [K']_u \big) \big\}$ by Lemma~\ref{lem:output-of-span}, since the set $[K']_u$ is strongly connected.
\end{proof}

Next we prove our extremal lemma, which will be used in Section~\ref{sec:ind} to prove the base case of our induction hypothesis (see Definition~\ref{def:ih} and Lemma~\ref{lem:basecase}). This lemma, which generalizes~\cite[Lemma~6.17]{BDMS} to the iceberg setting, follows easily from the iceberg spanning algorithm, together with Lemmas~\ref{lem:subadd} and~\ref{lem:extremal:basecase}.


\begin{lemma}[Extremal lemma for icebergs]\label{lem:extremal}
Let $u$ be a non-leaf vertex of\/ $T$ and let $\J = (J,v)$ be a $u$-iceberg spanned iceberg. Then $|J \cap A| = \Omega\big( \diam^*(\J) \big)$.
\end{lemma}

\begin{proof}
Since $\J$ is $u$-iceberg spanned, by Definition~\ref{def:iceberg spanned} there exists a set $K \subset J \cap A$ such that $[K]_u$ is strongly connected and $\J = \J_u\big( [K]_u \big)$. If we apply the $u$-iceberg spanning algorithm to the set $K$, then the algorithm begins with $|K|$ sets containing the individual elements of $K$, and finishes with the single set $K$. At each step of the algorithm the number of sets in the collection decreases by 1, and we claim that the quantity
\[
\sum_{K'\in \K^t} \diam^*\big( \J_u\big( [K']_u \big) \big)
\]
increases (additively) by at most $O(1)$, by Lemma~\ref{lem:subadd}. To see this, recall (cf.~the proof of Lemma~\ref{lem:AL}) that if $K_1$ and $K_2$ are united at time $t$ then the sets $[K_1]_u$, $[K_2]_u$ and $[K_1]_u \cup [K_2]_u$ are all strongly connected, and $K_1 \cup K_2 \subset K$. Moreover, we have $K \subset J \subset \Lambda \setminus H_T(u)$, by Lemma~\ref{lem:iceberg-cont-closed}, because $\J$ is the $u$-iceberg container of $[K]_u$. The inequality~\eqref{eq:subadd} therefore follows from Lemma~\ref{lem:subadd}, and this proves the claim.

Thus, comparing these quantities at times $t = 0$ and $t = |K| - 1$, it follows that
\[
\diam^*(\J) \le \sum_{x \in K} \diam^*\big( \J_u\big( [\{x\}]_u \big) \big) + O(|K|) = O(|K|),
\]
since $\diam^*\big( \J_u\big( [\{x\}]_u \big) \big) = O(1)$ for every $x \in \Lambda \setminus H_T(u)$ by Lemma~\ref{lem:extremal:basecase}. This implies that $|J \cap A| \ge |K| = \Omega\big( \diam^*(\J) \big)$, as required.
\end{proof}

We finish this section with two lemmas about the size and number of icebergs. The first, which provides an upper bound on $|J|$ in terms of $\diam^*(\J)$, where $\J = (J,u)$, will allow us 
to bound the number of choices for the set $J \cap A$.

\begin{lemma}\label{lem:iceberg-volume}
There exists $C > 0$ depending only on $T$ such that the following holds. Let $u$ be a non-leaf vertex of\/ $T$ and let\/ $\J = (J,u)$ be an iceberg. Then 
$$|J| \le \big( C \cdot \diam^*(\J) + 1 \big)^d.$$
\end{lemma}

\begin{proof}
Observe that there exists $x \in \R^d$ such that $J$ is contained in $x + [0,\diam(J)]^d$. Hence, since $\diam(J) \le C \cdot \wdiam_u(J)$ for some $C$ depending only on $T$, by Lemma~\ref{lem:wd-prop-diam}, and $\wdiam_u(J) \le \delta^{-d} \diam^*(\J)$, by 
Definition~\ref{def:diamstar}, the lemma follows.
\end{proof}

Our final lemma of the section allows us to count icebergs themselves, given limited information about their position and size.

\begin{lemma}\label{lem:iceberg-count}
There exists $C > 0$ depending only on $T$ such that the following holds. Let\/ $u$ be a non-leaf vertex of\/ $T$, and let\/ $\J = (J,u)$ be an iceberg with\/ $\diam^*(\J) \geq 2$. Then the number of icebergs $\J' = (J',v)$, with $v$ a non-leaf vertex of $T$, such that 
$$J \cap J' \neq \emptyset \qquad \text{and} \qquad \diam^*(\J') \leq \diam^*(\J),$$ 
is at most $\diam^*(\J)^C$.
\end{lemma}

\begin{proof}
Recall from Definition~\ref{def:iceberg} that the iceberg $\J'$ is uniquely defined by the choices of $v \in T$ and the set $\{ d_w(J') : w \in N_T^\to(v) \}$, where $d_w(J') := \max_{x \in J'} \< x,w \>$. Let us fix $y \in J \cap J'$ and observe that, by Lemma~\ref{lem:iceberg-volume}, there are at most $\big( C \cdot \diam^*(\J) + 1 \big)^d$ choices for $y$. We claim that there exists $C'$, depending only on $T$, such that 
\begin{equation}\label{eq:containment:inacube}
J' \subset y + [-C'\diam^*(\J),C'\diam^*(\J)]^d.
\end{equation}
Indeed, this follows as in the proof of Lemma~\ref{lem:iceberg-volume}, since we have
$$\diam(J') \leq C \cdot \wdiam_v(J') \leq C \delta^{-d} \cdot \diam^*(\J'),$$
where the first inequality holds by Lemma~\ref{lem:wd-prop-diam}, and the second by Definition~\ref{def:diamstar}, and since $J'\subset G_v(J')$ for any $v$-iceberg droplet $J'$. 

It follows from~\eqref{eq:containment:inacube} that, given $y \in J \cap J'$, for each $w \in N_T^\to(v)$ we have at most $\big( 2C'\diam^*(\J) + 1 \big)^d$ choices for the $x \in J'$ that maximizes $\< x,w \>$. Putting the pieces together, and recalling that $\diam^*(\J) \geq 2$, the number of choices for $\J'$ is at most
$$\sum_{v \in V(T)} \big( C \cdot \diam^*(\J) + 1 \big)^d \big( 2C'\diam^*(\J) + 1 \big)^{d \,\cdot\, |N_T^\to(v)|} \le \diam^*(\J)^{C''},$$
for some $C''$ depending only on $T$, as required.
\end{proof}

\section{One-step hierarchies}\label{one:step:sec}

Hierarchies, originally introduced by Holroyd~\cite{Hol} in the context of the 2-neighbour model in two dimensions, and subsequently developed further in~\cite{BBM3d,BBDM,DE,DH,BDMS,GHM,HM}, are a tool for encoding the formation of an internally spanned droplet $D$ by a relatively small number of internally spanned sub-droplets and `crossing' events (see Definition~\ref{def:delta-event}). Roughly speaking, this is achieved by running the spanning algorithm (for droplets) with initial set $D \cap A$, and recording only some of the sets that appear during the evolution of the algorithm. One significant problem with using hierarchies in our setting is that we must deal with droplets (and icebergs of various types) on very different scales: from polynomial in $p$, up to towers of exponentials of height $r - 2$. 

We overcome this issue by using a different (and, in fact, much simpler) notion of hierarchy, which we call `one-step hierarchies'. The reader who is familiar with the various notions of hierarchies used previously in the literature should think of these as consisting of a root vertex together with its neighbour(s) (see Lemma~\ref{lem:goodsat}). There will be enough leeway in our induction hypothesis in Section~\ref{sec:ind} that we can afford to apply Lemma~\ref{lem:goodsat} at each step rather than building up larger hierarchies as in previous works~\cite{Hol,BBDM,BDMS}.

We remark that for update families with resistance $r = 2$, we do not need the lemmas of this section. In those cases, Lemmas~\ref{lem:AL} and~\ref{lem:extremal} suffice to prove Theorem~\ref{thm:lower} (see Lemma~\ref{lem:basecase}). Moreover, for such models the only non-leaf vertex of $T$ is the root, in which case both lemmas are statements about $N_T^\to(\0)$-droplets. Thus, to prove the case $r = 2$ of Theorem~\ref{thm:lower} we do not require any of the icebergs framework that we developed in the previous two sections, and in fact the proof follows in much the same way as the lower bound for critical two-dimensional models in~\cite{BSU}. 

We begin by defining precisely the event that one iceberg `grows' into a larger iceberg. We shall only be interested in this event when the two icebergs have the same type.

\begin{definition}\label{def:delta-event}
Let $u$ be a non-leaf vertex of $T$, and let $\J = (J,u)$ and $\J' = (J',u)$ be icebergs such that $J' \subset J$. Then we define
\[
\Delta_u(\J',\J) := \big\{ \< K \>_u = \{\J\} \text{ for some } K \subset J' \cup (J \cap A) \big\}.
\]
\end{definition}


We are now ready to state and prove our only hierarchies-type lemma. It says, roughly speaking, that if $\J = (J,u)$ is a $u$-iceberg spanned iceberg (see Definition~\ref{def:iceberg spanned}), then one of two possibilities occurs. The first is that there exist two disjointly\footnote{Two increasing events $E$ and $F$ (depending on the set $A$) are said to \emph{occur disjointly} if there exist disjoint `witness' sets $X,Y \subset A$ such that $\{ X \subset A \} \Rightarrow E$ and $\{ Y \subset A \} \Rightarrow F$.} $u$-iceberg spanned icebergs, $\J_1'$ and $\J_2'$, whose iceberg droplets intersect $J$, and such that $\diam^*(\J_1')$ and $\diam^*(\J_2')$ are both reasonably large, and $\diam^*(\J_1') + \diam^*(\J_2')$ is not much smaller than $\diam^*(\J)$. The second is that there exists a $u$-iceberg spanned iceberg $\J' = (J',v)$ whose iceberg droplet $J'$ intersects $J$, with $\diam^*(\J')$ not much smaller than $\diam^*(\J)$, such that either $v \ne u$, or the event $\Delta_u(\J',\J)$ occurs. 

The precise (and slightly technical) statement of this `one-step hierarchies' lemma is as follows. Recall from Lemmas~\ref{lem:AL} and~\ref{lem:extremal} that $\lambda$ is a sufficiently large constant. 


\begin{lemma}\label{lem:goodsat}
Let $u$ be a non-leaf vertex of $T$, let $\J = (J,u)$ be a $u$-iceberg spanned $u$-iceberg, and let $\lambda \le y \leq \diam^*(\J)/4$. Then one of the following two events occurs:
\begin{itemize}
\item[$(a)$] There exist icebergs $\J_1' = (J_1',u_1)$ and $\J_2' = (J_2',u_2)$ such that $J_1' \cap J \neq \emptyset$ and $J_2' \cap J \neq \emptyset$, satisfying
\begin{equation}\label{eq:goodsat-diams}
y \leq \diam^*(\J_2') \leq \diam^*(\J_1') \leq \diam^*(\J) - 3y
\end{equation}
and
\begin{equation}\label{eq:goodsatsubadd}
\diam^*(\J) \leq \diam^*(\J_1') + \diam^*(\J_2') + 2y,
\end{equation}
and such that $\J_1'$ and $\J_2'$ are disjointly $u$-iceberg spanned.\smallskip
\item[$(b)$] There exists an iceberg $\J' = (J',v)$ with $J' \cap J \neq \emptyset$, satisfying 
$$y \leq \diam^*(\J) - \diam^*(\J') \leq 3y,$$
such that $\ice_u(\J')$ occurs. If $v = u$, then also $J' \subset J$ and $\Delta_u(\J',\J)$ occurs.
\end{itemize}
\end{lemma}

Observe that in~$(b)$ the two events $\ice_u(\J')$ and $\Delta_u(\J',\J)$ automatically occur disjointly, since for each fixed $\J'$ the two events depend on disjoint subsets of $A$.

\begin{proof}[Proof of Lemma~\ref{lem:goodsat}]
We claim that there exists a sequence
\[
J \cap A \supset K_0 \supset K_1 \supset \dots \supset K_m 
\]
such that setting
\[
\J_i := \J_u\big( [K_i]_u \big)
\]
for $0 \leq i \leq m$, we have the following: $[K_0]_u$ is strongly connected and $\J_0 = \J$, $|K_m| = 1$, and for every $1 \leq i \leq m$, $[K_i]_u$, $[K_{i-1} \setminus K_i]_u$ and $[K_i]_u \cup [K_{i-1} \setminus K_i]_u$ are all strongly connected, and
\begin{equation}\label{eq:goodsat-diam}
\diam^*\big( \J_u\big( [K_i]_u\big) \big) \geq \diam^*\big( \J_u\big( [K_{i-1} \setminus K_i]_u \big) \big).
\end{equation}
To construct the sequence, first we choose $K_0$ such that $[K_0]_u$ is strongly connected and $\J = \J_u\big( [K_0]_u \big)$, which we can do, by Definition~\ref{def:iceberg spanned}, because $\J$ is $u$-iceberg spanned. Now, supposing $K_{i-1} \subset J \subset \Lambda \setminus H_T(u)$ has already been chosen for some $i \geq 1$ and is such that $|K_{i-1}| \geq 2$ and $[K_{i-1}]_u$ is strongly connected, we apply Lemma~\ref{lem:penultimate} to $K_{i-1}$. This gives a subset $K_i$ of $K_{i-1}$ such that $[K_i]_u$, $[K_{i-1} \setminus K_i]_u$ and $[K_i]_u \cup [K_{i-1} \setminus K_i]_u$ are all strongly connected. Moreover, without loss of generality this subset satisfies~\eqref{eq:goodsat-diam}. 
Now, let $k \geq 1$ be minimal such that
\[
\diam^*(\J) - \diam^*(\J_k) \geq y.
\]
Note that such a $k$ exists because $\diam^*(\J) \ge 4y$, and since $|K_m| = 1$ and $K_m \subset J \subset \Lambda \setminus H_T(u)$, so $\diam^*(\J_m) \le \lambda \le y$, by Lemma~\ref{lem:extremal:basecase}. Suppose first that
\begin{equation}\label{eq:goodsat-3a}
\diam^*(\J) - \diam^*(\J_k) \geq 3y.
\end{equation}
We claim that in this case the conditions of~$(a)$ hold with
\begin{equation}\label{eq:goodsat-J1}
\J_1' = (J_1',u_1) := \J_k = \J_u\big( [K_k]_u \big)
\end{equation}
and
\begin{equation}\label{eq:goodsat-J2}
\J_2' = (J_2',u_2) := \J_u\big( [K_{k-1} \setminus K_k ]_u \big).
\end{equation}
Indeed, $\J_1'$ and $\J_2'$ are disjointly $u$-iceberg spanned by their definitions in~\eqref{eq:goodsat-J1} and~\eqref{eq:goodsat-J2}, since $[K_k]_u$ and $[K_{k-1} \setminus K_k]_u$ are both strongly connected by construction, and since $K_k \subset K_{k-1} \subset A$. To prove~\eqref{eq:goodsatsubadd}, observe that
\begin{equation}\label{eq:goodsat-parta-subadd}
\diam^*(\J_{k-1}) \le \diam^*(\J_1') + \diam^*(\J_2') + O(1)
\end{equation}
by Lemma~\ref{lem:subadd}, since the sets $[K_k]_u$, $[K_{k-1} \setminus K_k]_u$ and $[K_k]_u \cup [K_{k-1} \setminus K_k]_u$ are strongly connected, and $K_{k-1} \subset J \subset \Lambda \setminus H_T(u)$. It follows that
\[
\diam^*(\J) \leq \diam^*(\J_{k-1}) + y \le \diam^*(\J_1') + \diam^*(\J_2') + 2y
\]
where the first step follows from by the minimality of $k$, and the second from~\eqref{eq:goodsat-parta-subadd} and the bound $y \ge \lambda$, and therefore~\eqref{eq:goodsatsubadd} holds as claimed. To prove~\eqref{eq:goodsat-diams}, firstly observe that
\[
\diam^*(\J_1') \geq \diam^*(\J_2') \geq \diam^*(\J_{k-1}) - \diam^*(\J_k) - O(1) \geq y,
\]
where the first inequality holds by~\eqref{eq:goodsat-diam}, the second by~\eqref{eq:goodsat-parta-subadd} and since $\J_1' = \J_k$, and the third since $\diam^*(\J_{k-1}) - \diam^*(\J_k) \ge 2y$, which follows from~\eqref{eq:goodsat-3a} and the minimality of~$k$. Together with~\eqref{eq:goodsat-3a}, this gives~\eqref{eq:goodsat-diams}. Finally, both $J_1'$ and $J_2'$ intersect $J$ because both $K_k$ and $K_{k-1} \setminus K_k$ are contained in $J$. Thus, $\J_1'$ and $\J_2'$ satisfy the conditions of~$(a)$.

Now suppose that~\eqref{eq:goodsat-3a} does not hold, so we instead have
\begin{equation}\label{eq:goodsat-a}
y \leq \diam^*(\J) - \diam^*(\J_k) \leq 3y.
\end{equation}
We claim that the conditions of~$(b)$ are satisfied with
\begin{equation}\label{eq:goodsat-J'}
\J' = (J',v) := \J_k = \J_u\big( [K_k]_u \big).
\end{equation}
Indeed, $\J'$ is $u$-iceberg spanned by its definition in~\eqref{eq:goodsat-J'}, since $[K_k]_u$ is strongly connected by construction, and $K_k \subset A$. The bounds on $\diam^*(\J) - \diam^*(\J')$ are the same as those in~\eqref{eq:goodsat-a}, so they hold automatically, and we have $J' \cap J \neq \emptyset$ for the same reason as in part~$(a)$: the set $K_k$ is a subset of $J$. These are all the conditions that apply when $v \neq u$, so to deal with the remaining two conditions we assume in addition that $v = u$. These conditions are in fact straightforward: the inclusion $J' \subset J$ follows from $K_k \subset K_0$ and from $J'$ and $J$ being the smallest $u$-iceberg droplets containing $[K_k]_u$ and $[K_0]_u$ respectively (see Definition~\ref{def:iceberg-container}), and the event $\Delta_u(\J',\J)$ is an immediate consequence of the event $\ice_u(\J)$ and Lemma~\ref{lem:icebergspanned}, since we may take the set $K$ in Definition~\ref{def:delta-event} to be equal to any $K \subset J \cap A$ such that $\< K \>_u = \{ \J \}$.
\end{proof}

The `crossing' event $\Delta_u(\J',\J)$ appeared naturally in the proof of Lemma~\ref{lem:goodsat}, but it is not immediately clear how we can bound its probability. We will do so using our final deterministic property of icebergs, Lemma~\ref{lem:sideways-iceberg}, which will allow us to deduce from the event $\Delta_u(\J',\J)$ the occurrence of another event that we know better how to control: namely, the existence of a certain iceberg spanned iceberg. In the proof of Lemma~\ref{lem:sideways-iceberg} we shall use the following simple property of a pair of nested $u$-icebergs. 


\begin{lemma}\label{lem:nested-icebergs-diam}
Let $u$ be a non-leaf vertex of $T$, let $k := d_T(\0,u)$, and let $\J = (J,u)$ and $\J' = (J',u)$ be $u$-icebergs with $J' \subset J$. Then
$$\diam^*(\J) - \diam^*(\J') \leq \delta^k \cdot \max\big\{ d_v(J) - d_v(J') : v \in N_T^\to(u) \big\}.$$
\end{lemma}

\begin{proof}
By Definitions~\ref{def:weighted-diam} and~\ref{def:diamstar}, we have 
$$\diam^*(\J) - \diam^*(\J') = \delta^k \sum_{v \in M_T(u)} \lambda_v \Big( d_v\big( G_u(J) \big) - d_v\big( G_u(J') \big) \Big)$$
where $\Xi_u = \big\{ \lambda_v \,:\, v \in M_T(u) \big\}$. Now, by Definition~\ref{def:FuGu}, we have
$$d_v(G_u(J)) = d_v(J) \qquad \text{and} \qquad d_v(G_u(J')) = d_v(J')$$
for each $v \in N_T^\to(u)$. Moreover, if $u \neq \0$, then $d_{-u}(G_u(J)) = d_{-u}(G_u(J')) = 0$, since $N_T^\to(u)$ is not an $\SS^{d-1}$-bounding set, by Lemma~\ref{lem:icebergs-finite}. It follows that
\begin{align*}
\diam^*(\J) - \diam^*(\J') & \, = \, \delta^k \sum_{v \in N_T^\to(u)} \lambda_v \big( d_v(J) - d_v(J') \big)\\
& \, \le \, \delta^k \cdot \max\big\{ d_v(J) - d_v(J') : v \in N_T^\to(u) \big\},
\end{align*}
as required, since $\lambda_v \geq 0 $ for every $v \in M_T(u)$, and $\sum_{v \in M_T(u)} \lambda_v = 1$.
\end{proof}

Up to this point, we have defined icebergs (and their related structures and events) with respect to unions of half-spaces passing through the origin. However, in our final deterministic lemma we need to allow the role of the origin to be played by other elements of $\R^d$. As noted at the beginning of Section~\ref{sec:icebergs}, this is why we have worked in the shifted lattice $\Lambda$, instead of in $\Z^d$. 

In order to state and prove Lemma~\ref{lem:sideways-iceberg}, let us generalize some of the definitions of Sections~\ref{sec:icebergs} and~\ref{sec:subadd} by replacing the origin by some arbitrary element $x \in \R^d$. First, if $\0 \to u_1 \to \cdots \to u_k = u$ is a path in $T$, then we define
\begin{equation}\label{def:HTux}
H_T(u,x) := \bigcup_{i=1}^k \big( \HH_{u_i}^d + x \big) \cap \Lambda.
\end{equation}
Let $u$ be a non-leaf vertex of $T$, and let $K \subset \Lambda \setminus H_T(u,x)$ be a finite, non-empty set.  
\begin{itemize}
\item The \emph{$(u,x)$-iceberg closure} of $K$ is 
\begin{equation}\label{def:ux:closure}
[K]_{u,x} := \big[ K \cup H_T(u,x) \big]_\U \setminus H_T(u,x).
\end{equation}
\item An \emph{$x$-shifted $u$-iceberg} is a triple $\J = (J,u,x)$, where $J$ is a \emph{$(u,x)$-iceberg droplet}, i.e., a non-empty set of the form $D \setminus H_T(u,x)$, where $D$ is an $N_T^\to(u)$-droplet.\smallskip
\item The $(u,x)$-iceberg spanning algorithm is defined as in Definition~\ref{def:icespanalg}, except replacing the $u$-iceberg closure by the $(u,x)$-iceberg closure. The output of the algorithm is the \emph{$(u,x)$-iceberg span} of $K$, 
$$\< K \>_{u,x} := \big\{ \J_{u,x}\big( [K']_u \big) : K' \in \K^T \big\},$$
where $\J_{u,x}(K)$, the \emph{$(u,x)$-iceberg container} of $K$, is defined as in Definition~\ref{def:iceberg-container}, except with $J_i$ defined to be the minimal $(u_i,x)$-iceberg droplet containing $K$, and with each half-space $\HH_{u_i}^d$ replaced by $\HH_{u_i}^d + x$.\smallskip
\item The event $\ice_{u,x}(\J)$ holds (and we say that $\J$ is \emph{$(u,x)$-iceberg spanned}) if there exists $K \subset J \cap A$ such that $[K]_{u,x}$ is strongly connected and $\J = \J_{u,x}\big( [K]_{u,x} \big)$.\footnote{By Lemma~\ref{lem:icebergspanned}, this is equivalent to the event that $\< K \>_{u,x} = \{ \J \}$ for some $K \subset J \cap A$.}\smallskip
\item If $\J = (J,u,x)$ is an $x$-shifted $u$-iceberg, and $d_T(\0,u) = k$, then 
\begin{equation}\label{def:diamstar:x}
\diam_x^*(\J) := \delta^k \wdiam_u\big( G_{u,x}(J) \big),
\end{equation}
where $G_{u,x}(J) = \big\{ z \in F^*_u(J) : \< z,u\> \ge \<x,u\> \big\}$. 
\end{itemize}
Observe that, because we were working in the lattice $\Lambda = y + \Z^d$, where $y \in \R^d$ was arbitrary, all of the results proved in Sections~\ref{sec:icebergs} and~\ref{sec:subadd} also hold in this setting. 

\begin{lemma}\label{lem:sideways-iceberg}
Let $u$ be a vertex of $T$ with $d_T(\0,u) \le r - 3$, 
and let $\J = (J,u)$ and $\J' = (J',u)$ be icebergs with $J' \subset J$ and $\diam^*(\J) - \diam^*(\J') \ge \lambda$.  

If the event\/ $\Delta_u(\J',\J)$ occurs, then there exist $x \in \R^d$ and an $x$-shifted iceberg $\J''$, with
\begin{equation}\label{eq:sideways-diam}
\diam_x^*(\J'') \ge \delta^3 \cdot \big( \diam^*(\J) - \diam^*(\J') \big),
\end{equation}
such that $\< K \>_{v,x} = \{\J''\}$ for some $K \subset (J \cap A) \setminus J'$ and $v \in \{u\} \cup N_T^\to(u)$.
\end{lemma}

\begin{proof}
By Definition~\ref{def:delta-event}, the event $\Delta_u(\J',\J)$ implies that there exists a set $K' \subset J' \cup (J \cap A)$ such that $\< K' \>_u = \{\J\}$. Suppose first $K' \subset (J \cap A) \setminus J'$, and observe that the iceberg $\J'' = \J$ has the required properties. Indeed,~\eqref{eq:sideways-diam} holds trivially with $x = \0$, and $\< K' \>_{u,\0} = \< K' \>_u = \{\J''\}$, as required. We may therefore assume that $K' \cap J' \ne \emptyset$. 

We begin by choosing $v \in N_T^\to(u)$ and $x \in \R^d$. To do so, note first that, by Lemma~\ref{lem:nested-icebergs-diam}, there exists $v \in N_T^\to(u)$ such that
\begin{equation}\label{eq:sideways-diam-v}
\diam^*(\J) - \diam^*(\J') \le \delta^k \cdot \big( d_v(J) - d_v(J') \big),
\end{equation}
where $k := d_T(\0,u)$. Moreover, $v$ is not a leaf of $T$, since $k \le r - 3$ and $T$ has depth $r - 1$. Note that the set\footnote{We think of $S$ as the `strip' between $J'$ and $J$ in direction $v$. Note that $S$ is non-empty, by~\eqref{eq:sideways-diam-v} and since $\diam^*(\J) - \diam^*(\J') \ge \lambda > 0$.}
$$S := \big\{ z \in J : \< z,v \> > d_v(J') \big\}$$  
is finite, by Lemma~\ref{cor:icebergs-finite}, and choose $a \in S$ with $\< a,v \>$ minimal.

Now, let $\0 \to u_1 \to \cdots \to u_{k-1} \to u \to v$ be the path in $T$ from the root to $v$, and recall from Definition~\ref{def:tree} that the directions $\{u_1,\ldots, u_{k-1},u, v\}$ are linearly independent. Since there are at most $r - 2 \le d - 2$ linear constraints, it follows that there exists an element $x \in \R^d$ such that $\< x - a, v \> = 0$ and $\< x, w \> = 0$ for every $w \in \{u_1,\ldots, u_{k-1},u \}$. 

Having chosen $v \in N_T^\to(u)$ and $x \in \R^d$, we now define an $x$-shifted iceberg 
\begin{equation}\label{eq:sideways-J''}
\J'' := \J_{v,x}(L) = (J'',w,x),
\end{equation}
where $L$ is any strongly connected component of $[K' \cap S]_{v,x}$ such that $d_v(L) \ge d_v(J)$. We~claim that such a set $L$ exists, and that $\J''$ has the required properties. To prove this, we need the following simple claim, which follows from our choice of $x$.


\begin{claim}\label{clm:sideways1}
$[K']_u \setminus H_T(v,x) \subset [K' \cap S]_{v,x}$.
\end{claim}

\begin{clmproof}{clm:sideways1}
Note first that $H_T(v,x) = H_T(u) \cup \big( (\HH_v^d + a) \cap \Lambda \big)$, by~\eqref{def:HTux} and our choice of~$x$. Noting that $K' \subset (K' \cap S) \cup (\HH_v^d + a)$, since $K' \subset J$, we obtain
$$[K']_u \subset \big[ K' \cup H_T(u) \big]_\U \subset \big[ (K' \cap S) \cup H_T(v,x) \big]_\U.$$
Removing the set $H_T(v,x)$ from both sides, and recalling~\eqref{def:ux:closure}, the claim follows.
\end{clmproof}

Now, recall that $\< K' \>_u = \{\J\}$ and $\J = (J,u)$, and note that therefore $\J = \J_u\big( [K']_u \big)$, by Lemma~\ref{lem:output-of-span}. By Definition~\ref{def:iceberg-container}, it follows that $J$ is the minimal $u$-iceberg droplet containing $[K']_u$. Since $v \in N_T^\to(u)$, this implies that 
\begin{equation}\label{eq:yv:bigger:than:xv}
\<y,v\> = d_v([K']_u) = d_v(J) \ge \< x,v\>
\end{equation}
for some $y \in [K']_u$, the inequality holding because $\< x,v \> = \< a,v \>$ and $a \in J$. 

To deduce that $L$ exists, observe that, by Claim~\ref{clm:sideways1}, we have
$$[K']_u \subset [K' \cap S]_{v,x} \cup \big( \HH_v^d + x \big),$$
since $[K']_u \cap H_T(u) = \emptyset$. It follows that $y \in [K' \cap S]_{v,x}$ since, by~\eqref{eq:yv:bigger:than:xv}, we have $\<y,v\> \ge \<x,v\> > \<z,v\>$ for all $z \in \HH_v^d + x$. In particular, if $L$ is the strongly connected component of $[K' \cap S]_{v,x}$ containing $y$, then $d_v(L) \ge \<y,v\> \ge d_v(J)$. 

Next, observe that $L = [L \cap K' \cap S]_{v,x}$, since $L$ is a strongly connected component, so must contain every site of $K' \cap S$ that is used to infect it. Set 
$$K := L \cap K' \cap S \subset (J \cap A) \setminus J',$$ 
so $L = [K]_{v,x}$ is strongly connected. Recalling from~\eqref{eq:sideways-J''} that $\J'' = \J_{v,x}(L)$, it follows by Lemma~\ref{lem:output-of-span} that $\< K \>_{v,x} = \{\J''\}$, as required. 

It remains to verify the lower bound~\eqref{eq:sideways-diam} on $\diam_x^*(\J'')$. To do so, recall~\eqref{def:diamstar:x}, and observe first that
\begin{equation}\label{eq:sideways-diam:step}
\diam_x^*(\J'') \geq \delta^{k+1} \cdot \wdiam_w(J'') \geq \delta^{k+2} \cdot \diam(J'') \end{equation}
where the first inequality holds because $J'' \subset G_{w,x}(J'')$ and $d_T(\0,w) \leq k+1$, and the second by Lemma~\ref{lem:wd-prop-diam}. To deduce~\eqref{eq:sideways-diam}, we shall use~\eqref{eq:sideways-diam-v} and the following claim. 

\begin{claim}\label{clm:sideways2}
$L$ is strongly connected to $\HH_v^d + x$. 
\end{claim}

\begin{clmproof}{clm:sideways2}
Observe first that, by Claim~\ref{clm:sideways1}, we have
\begin{equation}\label{clm:sideways:second:application}
[K' \cap S]_{v,x} = \big[ [K']_u \cap S \big]_{v,x}.
\end{equation}
We shall show that $\big( [K']_u \cap S \big) \cup \big( \HH_v^d + x \big)$  is strongly connected, which will in turn imply that every strongly connected component of $[K']_u \cap S$ is strongly connected to $\HH_v^d + x$. 
By~\eqref{clm:sideways:second:application}, it will then follow that every strongly connected component of $[K' \cap S]_{v,x}$ is strongly connected to $\HH_v^d + x$, as required.

To show that $\big( [K']_u \cap S \big) \cup \big( \HH_v^d + x \big)$  is strongly connected, recall that $\< K' \>_u = \{\J\}$, which implies that $[K']_u$ is strongly connected, by Lemma~\ref{lem:output-of-span}. It also follows from $\< K' \>_u = \{\J\}$ that $[K']_u \subset J$, and therefore that $[K']_u \setminus \big( \HH_v^d + x \big) \subset S$, by the definitions of $S$ and $x$. Since $K' \cap J'$ is non-empty and $J' \subset \HH_v^d + x$, it follows that the set
\[
[K']_u \cup \big( \HH_v^d + x \big) = \big( [K']_u \cap S \big) \cup \big( \HH_v^d + x \big)
\]
is strongly connected, as required.
\end{clmproof}

By Claim~~\ref{clm:sideways2}, and since $d_v(L) \ge d_v(J)$, it follows that
$$\diam(J'') \ge d_v(J) - \< x,v\> - R,$$
and hence, by~\eqref{eq:sideways-diam:step}, and noting that $\< x,v\> \le d_v(J') + R$, we obtain
$$\diam_x^*(\J'') \ge \delta^{k+2} \big( d_v(J) - d_v(J') - 2R \big).$$
Combining this with~\eqref{eq:sideways-diam-v}, it follows that
\[
\diam_x^*(\J'') \ge \delta^3 \cdot \big( \diam^*(\J) - \diam^*(\J') \big),
\]
as required, since $\diam^*(\J) - \diam^*(\J') \ge \lambda$. 
\end{proof}

\begin{remark}\label{remark:choices}
Observe that, in the proof above, we took $x$ equal either to $\0$ (when $v = u$), or to a deterministic function of the pair $(J,J')$ and the vertex $v \in N_T^\to(u)$. This fact will be used in the next section in order to bound the number of choices for $\J''$. 
\end{remark}

\section{The proof of Theorem~\ref{thm:lower}}\label{sec:ind}

Our toolbox of deterministic lemmas is now complete, and we are ready to begin assembling the parts of the proof of Theorem~\ref{thm:lower}. The proof will be by induction on $r$, and also by induction on the diameter of the iceberg.

To begin, we need to fix some constants. In view of Lemma~\ref{lem:extremal}, choose $c > 0$ such that, whenever $u \in V(T)$ is not a leaf\footnote{Recall that the tree $T$ was fixed in Section~\ref{subsec:tree}.} and $\J = (J,v)$ is $u$-iceberg spanned, we have
\begin{equation}\label{eq:little-c}
| J \cap A | \geq c \cdot \diam^*(\J).
\end{equation}
Next, select an arbitrary sequence $\eps(2), \dots, \eps(r)$ 
such that
\begin{equation}\label{eq:constants1}
\frac{1}{d-1} > \eps(2) > \cdots > \eps(r) > 0.
\end{equation}
The function
\begin{equation}\label{eq:xs}
x(s) := \exp_{(s-2)}\big( p^{-\eps(s)} \big),
\end{equation}
defined for each $2 \leq s \leq r$, will be (up to a constant factor) the maximum diameter of the icebergs of type $u$ that we consider when $d_T(\0,u) = r - s$. 

Recall from Section~\ref{sec:subadd} that 
$\lambda$ is a sufficiently large constant, 
which is allowed to depend on all other constants (in particular, it will depend on $\eps(2)$ and $c$).\footnote{We shall need to assume that $\lambda > \exp_{(r-1)}\big( c^{-1} (1 - \eps(2)(d-1))^{-1} \big)$.} We emphasize that $p$ is not considered to be a constant, and will always be assumed to be sufficiently small. 
Now, let $q \colon \R \to \R$ be the function defined by $q(x) = 1$ for all $x < 1$, and
\begin{equation}\label{eq:qx}
q(x) := \exp\Bigg( -x \cdot \frac{\log 1/p}{\log_{(r-1)}\big( \lambda \cdot x \big)} \Bigg)
\end{equation}
for every $x \ge 1$. Thus, up to constant factors in the exponents, $q(x)$ interpolates between $p^x$ when $x = O(1)$, and $e^{-x}$ when $x$ is close to $x(r)$. This will be our bound on the probability that an iceberg of diameter $x$ is $u$-iceberg spanned (for appropriate $u$).

When we introduced the iceberg diameter $\diam^*(\J)$ in Definition~\ref{def:diamstar}, we `gave away' a factor of $\delta^k$, which allowed us to prove Lemma~\ref{lem:subadd} even when the icebergs involved were of different types. We now give away another factor of $\delta^k$ (though in this case a factor of $2^k$ would be more than sufficient), which will allow us to use the trivial bound on the probability of the event $\Delta_u(\J',\J)$ whenever $\J$ and $\J'$ are icebergs of different types. To that end, if $\J = (J,u)$ is an iceberg and $k := d_T(\0,u)$, then we define
\begin{equation}\label{def:diamstarstar}
\diam^{**}(\J) := \delta^k \diam^*(\J).
\end{equation}
We can now state our induction hypothesis. 

\begin{definition}[The induction hypothesis]\label{def:ih}
For each $2 \leq s \leq r$ and each $x \leq x(s)$, let $\ih(s,x)$ denote the following statement:
\begin{itemize}
\item[$(\ast)$] For each $u \in V(T)$ with $d_T(\0,u) = r - s$, and every iceberg $\J$ with $\diam^*(\J) \leq x$,
\begin{equation}\label{eq:indbound}
\P_p \big( \ice_u(\J) \big) \leq q\big( \diam^{**}(\J) \big).
\end{equation}
\end{itemize}
\end{definition}

We shall prove the induction hypothesis in two steps, as follows.

\begin{lemma}\label{lem:basecase}
The statement $\ih\big(s,x(2)\big)$ holds for each $2 \leq s \leq r$.
\end{lemma}

\begin{lemma}\label{lem:indstep}
For each $3 \leq s\leq r$ and each $x(2) < x \leq x(s)$, we have
\[
\ih\big(s-1,x(t-1)\big) \, \cap \, \ih(s,x-1) \, \cap \, \Big( \ih(s+1,x) \, \cap \, \dots \cap \, \ih\big(r,x\big) \Big) \, \Rightarrow \, \ih(s,x),
\]
where $3 \leq t \leq s$ is such that $x(t-1) < x \leq x(t)$.
\end{lemma}

Before proving Lemmas~\ref{lem:basecase} and~\ref{lem:indstep}, let us quickly observe that together they imply that $\ih(s,x)$ holds for all $2 \leq s \leq r$ and all $x \leq x(s)$, and show that this is sufficient to prove Proposition~\ref{prop:critdrop}. In fact, we shall prove the following slightly stronger proposition, which moreover specifies the set $\T$ explicitly (this will turn out to be useful later).

\begin{prop}\label{prop:critdrop-v2}
Let $\J = (J,\0)$ be a $\0$-iceberg such that $\diam^*(\J) \leq x(r)$. Then 
$$\P_p \big( \ice_\0(\J) \big) \le \exp\big( - \diam^*(\J)\big).$$
\end{prop}

\begin{proof}
Suppose Lemmas~\ref{lem:basecase} and~\ref{lem:indstep} hold. We shall prove, by induction on $t$, that $\ih(s,x)$ holds for every $t \leq s \leq r$ and all $x \le x(t)$. The case $t = 2$ holds by Lemma~\ref{lem:basecase}, so let $3 \le t \le r$ and suppose that we have already proved that $\ih(s,x)$ holds for every $t - 1 \le s \le r$ and all $x \le x(t-1)$. By Lemma~\ref{lem:indstep} and induction on $x$, it follows that $\ih(r,x)$ holds for all $x \le x(t)$. Now, again by Lemma~\ref{lem:indstep} and induction on $x$, it follows that $\ih(r-1,x)$ holds for all $x \le x(t)$. Repeating this argument $r - t + 1$ times, we deduce that $\ih(s,x)$ holds for every $t \leq s \leq r$ and all $x \le x(t)$, as claimed. 

It follows that $\ih(s,x)$ holds for every $2 \leq s \leq r$ and $x \le x(s)$, and so, in particular, $\ih(r,x(r))$ holds. Therefore, if $\J$ is a $\0$-iceberg with $\diam^*(\J) \le x(r)$, then 
\[
\P_p \big( \ice_\0(\J) \big) \leq q\big( \diam^{**}(\J) \big) = q\big( \diam^*(\J) \big) \leq \exp\big( -\diam^*(\J) \big),
\]
where the first step holds by $\ih(r,x(r))$, applied with $u = \0$, the second because $\J$ is a $\0$-iceberg, so $\diam^{**}(\J) = \diam^*(\J)$, and the third by~\eqref{eq:xs} and~\eqref{eq:qx}, since $\eps(r) < 1$. 
\end{proof}

Proposition~\ref{prop:critdrop} follows easily from Proposition~\ref{prop:critdrop-v2}, but let us spell out the details for completeness. 

\begin{proof}[Proof of Proposition~\ref{prop:critdrop}]
Set $\T := N_T^\to(\0)$ and $\eps := \eps(r)$, 
and let $D$ be a $N_T^\to(\0)$-droplet with $\diam(D) \le x(r) = \exp_{(r-2)}(p^{-\eps})$. Note that $D$ is a $\0$-iceberg droplet, by Definition~\ref{def:iceberg}, set $\J := (D,\0)$, and observe that
\begin{equation}\label{eq:diam:0iceberg:xr}
\diam^*(\J) = \wdiam_\0\big( G_\0(D) \big) = \wdiam_\0(D) \le \diam(D) \le x(r),
\end{equation}
where the first step is Definition~\ref{def:diamstar}, the second follows from Definitions~\ref{def:weighted-diam} and~\ref{def:FuGu}, the third holds by Lemma~\ref{lem:wd-prop-diam}, and the fourth holds by assumption. Hence, by Proposition~\ref{prop:critdrop-v2}, we have  
$$\P_p \big( \ice_\0(\J) \big) \le \exp\big( - \diam^*(\J) \big) \le \exp\big( - \delta' \cdot \diam(D) \big)$$
for some sufficiently small constant $\delta' > 0$, where the final inequality holds because $\diam^*(\J) = \wdiam_\0(D)$,  by~\eqref{eq:diam:0iceberg:xr}, and by Lemma~\ref{lem:wd-prop-diam}. 

To complete the proof, we shall show that the event $\ice_\0(\J)$ is the same as the event that the droplet $D$ is internally $N_T^\to(\0)$-spanned. To see this, recall from Definition~\ref{def:iceberg spanned} that, since $\J = (D,\0)$, the event $\ice_\0(\J)$ holds if and only if there exists $K \subset D \cap A$ such that $[K]_\0$ is strongly connected and $\J_\0( [K]_\0 ) = (D,\0)$. Note that $[K]_\0 = [K]_\U$, by Definition~\ref{def:u-closure}, and that $\J_\0( [K]_\0 ) = (J,\0)$, where $J$ is the minimal $N_T^\to(\0)$-droplet containing $[K]_\0$, by Definitions~\ref{def:iceberg} and~\ref{def:iceberg-container}. It follows that  $\ice_\0(\J)$ holds if and only if there exists $K \subset D \cap A$ such that $[K]_\U$ is strongly connected and $D$ is the minimal $N_T^\to(\0)$-droplet containing $[K]_\U$. By Definition~\ref{def:intspan}, this is the same as the event that $D$ is internally $N_T^\to(\0)$-spanned, as required. 
\end{proof}

In order to prove Proposition~\ref{prop:critdrop-v2}, it remains to prove Lemmas~\ref{lem:basecase} and~\ref{lem:indstep}. We begin with Lemma~\ref{lem:basecase}, which follows easily from our extremal lemma, Lemma~\ref{lem:extremal}, and the simple bound on the volume of an iceberg given by Lemma~\ref{lem:iceberg-volume}. 

\begin{proof}[Proof of Lemma~\ref{lem:basecase}]
Let $2 \leq s \leq r$, and let $u \in V(T)$ with $d_T(\0,u) = r - s$, so in particular $u$ is not a leaf. Let $\J = (J,v)$ be an iceberg with $\diam^*(\J) \le x(2) = p^{-\eps(2)}$. We are required to show that
\begin{equation}\label{eq:ih:basecase:need}
\P_p \big( \ice_u(\J) \big) \le q\big( \diam^{**}(\J) \big).
\end{equation}
The bound holds trivially if $\diam^{**}(\J) < 1$, since $q(x) = 1$ for every $0 \le x < 1$. We may therefore assume that $\diam^{**}(\J) \ge 1$, and hence 
\begin{equation}\label{eq:volume:app1}
|J| \le \big( C \cdot \diam^*(\J) + 1 \big)^d \le \lambda \cdot \diam^*(\J)^d,
\end{equation}
by Lemma~\ref{lem:iceberg-volume}, and since $\lambda$ is sufficiently large. Now, recall from~\eqref{eq:little-c} (which followed from Lemma~\ref{lem:extremal}) that if $\J$ is $u$-iceberg spanned, then $|J \cap A| \ge c \cdot \diam^*(\J)$. Hence, by Markov's inequality, $\J$ is $u$-iceberg spanned with probability at most
\begin{equation}\label{eq:basecase}
{\lambda x^d \choose c x} p^{c x} = O\big( x^{d-1} p \big)^{c x} \leq O\big( p^{1 - \eps(2)(d-1)} \big)^{c x} = p^{\Omega(x)} \le q(x),
\end{equation} 
where $x = \diam^*(\J)$, the second step holds since $x \le p^{-\eps(2)}$, then third holds because $\eps(2) < 1 / (d - 1)$, and the last follows from~\eqref{eq:qx}, since $\lambda$ is sufficiently large.

Finally, note that $x = \diam^*(\J) \ge \diam^{**}(\J)$, and that $q$ is decreasing. It follows that $q(x) \le q\big( \diam^{**}(\J) \big)$, and hence~\eqref{eq:ih:basecase:need} follows from~\eqref{eq:basecase}.
\end{proof}

The proof of Lemma~\ref{lem:indstep} (the induction step) is significantly more challenging. The idea is to apply Lemma~\ref{lem:goodsat} to the iceberg $\J$, and then bound the probabilities of the various events using Lemma~\ref{lem:sideways-iceberg} and the induction hypothesis. We shall also make use of the following fundamental lemma of van den Berg and Kesten~\cite{BK}. If $E$ and $F$ are increasing events, then we write $E \circ F$ for the event that $E$ and $F$ occur disjointly.  

\begin{lemma}[The van den Berg--Kesten inequality]\label{lem:vBK}
Let $E$ and $F$ be any two increasing events and let $p \in (0,1)$. Then
$$\Pr_p(E \circ F) \, \le \, \Pr_p(E)\,\Pr_p(F).$$
\end{lemma}

Three straightforward inequalities will be needed in order to show that the bounds we obtain are sufficiently strong. Since the proofs are unenlightening, they are deferred to Appendix~\ref{app:calcs}. 
Let $C > 0$ be a sufficiently large constant (depending on $T$), and fix a sufficiently small constant $\eta > 0$ with $(1+\eta)\eps(3) < \eps(2)$. The first inequality will be used when event~$(a)$ of Lemma~\ref{lem:goodsat} occurs. We write $\beta \gg 1$ to denote that $\beta \to \infty$ as $p \to 0$. 

\begin{lemma}\label{lem:splitter}
Let $r \geq 3$, and let $\alpha \ge \beta \gg 1$ and $\gamma \in \R$ satisfy 
$$\gamma \le C (\log \alpha)^{1+\eta} \le C^2 \beta.$$ 
Then
\[
q(\alpha) \cdot q(\beta) \cdot (\alpha + \beta + \gamma)^{4C} \le q(\alpha + \beta + \gamma).
\]
\end{lemma}

The second inequality will be used when event~$(b)$ of Lemma~\ref{lem:goodsat} occurs, but only when the icebergs $\J$ and $\J'$ of that lemma are of different types.

\begin{lemma}\label{lem:giveaway}
Let $\alpha \ge \gamma \gg 1$. Then
\[
q(\alpha) \cdot (\alpha + \gamma)^{3C} \le q\big( \delta(\alpha + \gamma) \big).
\]
\end{lemma}

The third inequality will also be used when event~$(b)$ of Lemma~\ref{lem:goodsat} occurs and the icebergs of that lemma are of the same type.

\begin{lemma}\label{lem:sideways}
Let $r \ge 3$ and $\alpha,\gamma \gg 1$ be such that
\begin{equation}\label{eq:sideways-z}
C^{-1} \big( \log(\alpha + \gamma) \big)^{1 + \eta} \le \gamma \le C \big( \log(\alpha + \gamma) \big)^{1 + \eta}.
\end{equation}
Then
\[
q(\alpha) \cdot q(\delta^4 \gamma) \cdot (\alpha + \gamma)^{5C} \leq q(\alpha + \gamma).
\]
\end{lemma}

\pagebreak

We are now ready to prove the induction step, Lemma~\ref{lem:indstep}.

\begin{proof}[Proof of Lemma~\ref{lem:indstep}]
Let $3 \leq s \leq r$, let $x(2) < x \leq x(s)$, and let $3 \leq t \leq s$ be such that $x(t-1) < x \leq x(t)$. Our aim is to show that $\ih(s,x)$ holds, so let $u \in V(T)$ be a vertex with $d_T(\0,u) = k := r - s$, and let $\J = (J,w)$ be an iceberg with $\diam^*(\J) \leq x$. Note in particular that $u$ is not a leaf of $T$. Under the assumption that
\begin{equation}\label{eq:indstep-hyp}
\ih\big(s-1,x(t-1)\big) \, \cap \, \ih(s,x-1) \, \cap \, \Big( \ih(s+1,x) \, \cap \, \dots \cap \, \ih\big(r,x\big) \Big)
\end{equation}
holds, we must prove
\begin{equation}\label{eq:indstep-aim}
\P_p\big( \ice_u(\J) \big) \leq q\big( \diam^{**}(\J) \big).
\end{equation}

To begin, let us eliminate cases that are either trivial or dealt with by~\eqref{eq:indstep-hyp}. First, by Lemma~\ref{lem:basecase}, we may assume that $\diam^*(\J) > x(2)$. Next, we may assume that $w \in [\0,u]_T$, because if not then $\J$ cannot be $u$-iceberg spanned, so $\P_p\big( \ice_u(\J) \big) = 0$. Third, note that the event $\ice_u(\J)$ implies $\ice_w(\J)$, by Lemma~\ref{lem:iceberg span-reduce}, and therefore if $d_T(\0,w) = r - s'$ for some $s + 1 \le s' \le r$, then 
\[
\P_p\big( \ice_u(\J) \big) \leq \P_p\big( \ice_w(\J) \big) \leq q\big( \diam^{**}(\J) \big)
\]
by $\ih(s',x)$, which holds by~\eqref{eq:indstep-hyp}. We may therefore assume that $w = u$. 

Suppose then that $\J = (J,u)$ and that $\ice_u(\J)$ occurs. Apply Lemma~\ref{lem:goodsat} to $\J$, with
\begin{equation}\label{eq:indstep-a}
y := \big( \log \diam^*(\J) \big)^{1+\eta},
\end{equation}
observing that $\lambda \le y \leq \diam^*(\J)/4$, since $\diam^*(\J) > x(2) \gg 1$ as $p \to 0$. The remainder of the proof is divided into three claims, according to whether event~$(a)$ or event~$(b)$ in Lemma~\ref{lem:goodsat} occurs, and in the case of event~$(b)$, whether or not the new iceberg $\J'$ is also of type $u$. We refer to these three events as $E_1$, $E_2$ and $E_3$ respectively. That is:
\begin{itemize}
\item $E_1$ is the event that there exist icebergs $\J_1 = (J_1,u_1)$ and $\J_2 = (J_2,u_2)$, with $J_1 \cap J \neq \emptyset$ and $J_2 \cap J \neq \emptyset$, satisfying
\begin{equation}\label{eq:indstep-case-a-diams1}
y \le \diam^*(\J_2) \le \diam^*(\J_1) \le \diam^*(\J) - 3y
\end{equation}
and
\begin{equation}\label{eq:indstep-case-a-diams2}
\diam^*(\J) \le \diam^*(\J_1) + \diam^*(\J_2)+ 2y,
\end{equation}
and such that $\J_1$ and $\J_2$ are disjointly $u$-iceberg spanned.\smallskip
\item $E_2$ is the event that there exists an iceberg $\J' = (J',v)$ with $v \neq u$, $J' \cap J \neq \emptyset$,
\begin{equation}\label{eq:indstep-case-b-diams}
y \le \diam^*(\J) - \diam^*(\J') \le 3y,
\end{equation}
and such that $\ice_u(\J')$ occurs.\smallskip
\item $E_3$ is the event that there exists an iceberg $\J' = (J',u)$ such that $J' \subset J$, the bounds~\eqref{eq:indstep-case-b-diams} hold, and the events $\ice_u(\J')$ and $\Delta_u(\J',\J)$ both occur.
\end{itemize}

We shall show in the next three claims that
\begin{equation}\label{eq:indstep-aim-cond}
\P_p(E_i) \leq q\big( \diam^{**}(\J) \big) / 3
\end{equation}
for each $i \in \{1,2,3 \}$. Since, by Lemma~\ref{lem:goodsat}, the event $\ice_u(\J)$ implies that the event $E_1 \cup E_2 \cup E_3$ holds, these three bounds together imply~\eqref{eq:indstep-aim} via the union bound. The three claims below will therefore complete the proof of the lemma.

\begin{claim}\label{clm:indstep-a}
$\P_p(E_1) \leq q\big( \diam^{**}(\J) \big) / 3$.
\end{claim}

To prove the claim we shall use the union bound over pairs $(\J_1,\J_2)$ of icebergs, bounding the number of choices for $\J_1$ and $\J_2$ using Lemma~\ref{lem:iceberg-count}, and the probability that the icebergs $\J_1$ and $\J_2$ are disjointly $u$-iceberg spanned using the van den Berg--Kesten inequality (Lemma~\ref{lem:vBK}) and $\ih(s,x-1)$ (which is one of our assumptions in~\eqref{eq:indstep-hyp}). 

\begin{clmproof}{clm:indstep-a}
First, set
$$\alpha := \diam^*(\J_1), \qquad \beta := \diam^*(\J_2) \qquad \text{and} \qquad \gamma := \diam^*(\J) - \alpha - \beta.$$
By Lemma~\ref{lem:iceberg-count}, it follows that we have at most 
$$\diam^*(\J)^{2C} \le \big( \delta^k (\alpha + \beta + \gamma) \big)^{3C}$$
choices for the icebergs $\J_1$ and $\J_2$. Indeed, the inequality follows since $\diam^*(\J) = \alpha + \beta + \gamma$ and since $\diam^*(\J) > x(2) \gg 1$. Next, observe that
$$\P_p \big( \ice_u(\J_i) \big) \le q\big( \diam^{**}(\J_i) \big)$$
for each $i \in \{1,2\}$, by~$\ih(s,x-1)$, since $\diam^*(\J_i) \le \diam^*(\J) - 3y$ and $\diam^*(\J) \le x$. Moreover, recalling~\eqref{def:diamstarstar} and that $q$ is a decreasing function, we have 
$$q\big( \diam^{**}(\J_1) \big) \leq q(\delta^k \alpha) \qquad \text{and} \qquad q\big( \diam^{**}(\J_2) \big) \leq q(\delta^k \beta),$$
since if $\J_i$ is $u$-iceberg spanned then $u_i \in [\0,u]_T$.   

Thus, by the van den Berg--Kesten inequality (Lemma~\ref{lem:vBK}), and applying the union bound over all choices of $\J_1$ and $\J_2$, we obtain
\begin{equation}\label{eq:indstep-case-a-sup}
\P_p(E_1) \leq \sup_{\alpha,\beta,\gamma} \; \big( \delta^k (\alpha + \beta + \gamma) \big)^{3C} \cdot q(\delta^k \alpha) \cdot q(\delta^k \beta),
\end{equation}
where the supremum is taken over the valid choices of $\alpha$, $\beta$ and $\gamma$.

We shall bound~\eqref{eq:indstep-case-a-sup} using Lemma~\ref{lem:splitter}, applied to $\delta^k \alpha$, $\delta^k \beta$ and $\delta^k \gamma$. We may do so, since 
$$\alpha \ge \beta \ge y = \big( \log \diam^*(\J) \big)^{1+\eta} \gg 1,$$
by~\eqref{eq:indstep-case-a-diams1} and since $\diam^*(\J) > x(2)$, and 
$$\gamma \le 2y = 2\big( \log(\alpha+ \beta + \gamma) \big)^{1+\eta} \le 3(\log \alpha)^{1+\eta} \le 3y \le 3 \beta,$$ 
by~\eqref{eq:indstep-case-a-diams1} and~\eqref{eq:indstep-case-a-diams2}, and noting that $\beta + \gamma \le 3\alpha$. It follows, by Lemma~\ref{lem:splitter} and~\eqref{eq:indstep-case-a-sup}, that
\begin{equation}\label{eq:indstep-bound1}
\P_p(E_1) \le \sup_{\alpha,\beta,\gamma} \; \frac{q\big( \delta^k (\alpha + \beta + \gamma) \big)}{\big( \delta^k (\alpha + \beta + \gamma) \big)^C} \le \frac{q\big( \diam^{**}(\J) \big)}{3},
\end{equation}
as claimed, since $\diam^{**}(\J) = \delta^k (\alpha + \beta + \gamma) \gg 1$. 
\end{clmproof}

We next turn to the events $E_2$ and $E_3$. In both events we have an iceberg $\J' = (J',v)$, with $J' \cap J \neq \emptyset$, such that~\eqref{eq:indstep-case-b-diams} holds and the event $\ice_u(\J')$ occurs. We shall take a union bound over choices for $\J'$, bounding the number of choices using Lemma~\ref{lem:iceberg-count}, and the probability that $\J'$ is $u$-iceberg spanned using the $\ih(s,x-1)$. For the event $E_3$, an extra argument using Lemma~\ref{lem:sideways-iceberg} will be needed to bound the probability of $\Delta_u(\J',\J)$. 

To save repetition, let us perform the common part of the argument here, before separating into cases. Set
\begin{equation}\label{eq:alpha:gamma}
\alpha := \diam^*(\J') \qquad \text{and} \qquad \gamma := \diam^*(\J) - \alpha,
\end{equation}
and observe that, by Lemma~\ref{lem:iceberg-count}, we have at most 
\begin{equation}\label{eq:counting:icebergs:claims}
\diam^*(\J)^{C} \le \big( \delta^k (\alpha + \gamma) \big)^{2C}
\end{equation}
choices for the iceberg $\J'$, where the inequality follows since $\diam^*(\J) = \alpha + \gamma > x(2)$. Next, observe that
\begin{equation}\label{eq:applying:induction:claims}
\P_p \big( \ice_u(\J') \big) \le q\big( \diam^{**}(\J') \big)
\end{equation}
by~$\ih(s,x-1)$, since $\diam^*(\J') \le \diam^*(\J) - y$, by~\eqref{eq:indstep-case-b-diams}, and $\diam^*(\J) \le x$. 

We are now ready to bound the probability of the event $E_2$. The reader may be concerned that we ought to lose here, because $\diam^*(\J') < \diam^*(\J)$. We are saved, however, by the extra factor of $\delta^k$ in $\diam^{**}(\J)$ (compared with $\diam^*(\J)$), and the fact that $v$ is strictly closer to $\0$ than $u$, which imply that $\diam^{**}(\J') > \diam^{**}(\J)$. 

\begin{claim}\label{clm:indstep-b}
$\P_p(E_2) \leq q\big( \diam^{**}(\J) \big) / 3$.
\end{claim}

\begin{clmproof}{clm:indstep-b}
Observe first that $d_T(\0,v) \leq k - 1$, since if $\J' = (J',v)$ is $u$-iceberg spanned then $v \in [\0,u]_T$, and $v \ne u$. We therefore have $\diam^{**}(\J') \ge \delta^{k-1} \alpha$, by~\eqref{def:diamstarstar}. Since $q$ is a decreasing function, it follows from~\eqref{eq:applying:induction:claims} that
$$\P_p \big( \ice_u(\J') \big) \le q\big( \diam^{**}(\J') \big) \le q(\delta^{k-1} \alpha).$$
Hence, taking a union bound over the choices of $\J'$, and recalling~\eqref{eq:counting:icebergs:claims}, we obtain
\begin{equation}\label{eq:indstep-case-b1-sup}
\P_p(E_2) \leq \sup_{\alpha,\gamma} \big( \delta^{k} (\alpha + \gamma) \big)^{2C} \cdot q( \delta^{k-1} \alpha),
\end{equation}
where the supremum is taken over valid values of $\alpha$ and $\gamma$.

We shall bound~\eqref{eq:indstep-case-b1-sup} using Lemma~\ref{lem:giveaway}, applied to $\delta^{k-1} \alpha$ and $\delta^{k-1} \gamma$. The conditions of the lemma are satisfied because $y \le \gamma \le 3y$ and $\alpha \ge \diam^*(\J) - 3y$, by~\eqref{eq:indstep-case-b-diams}, and therefore $\alpha \ge \gamma \gg 1$, since $y = ( \log \diam^*(\J))^{1+\eta}$ and $\diam^*(\J) > x(2) \gg 1$. Thus, by~\eqref{eq:indstep-case-b1-sup} and Lemma~\ref{lem:giveaway}, we have 
\begin{equation}\label{eq:indstep-bound2}
\P_p(E_2) \le \sup_{\alpha,\gamma} \frac{q\big( \delta^k (\alpha + \gamma) \big)}{\big( \delta^{k} (\alpha + \gamma) \big)^C} \le \frac{q\big( \diam^{**}(\J) \big)}{3},
\end{equation}
since $\diam^{**}(\J) = \delta^k (\alpha + \gamma) \gg 1$. 
\end{clmproof}

The proof of the final claim is similar, but the details are somewhat more complicated because we also need to bound the probability of the event $\Delta_u(\J',\J)$. 
This will be done by first applying Lemma~\ref{lem:sideways-iceberg} to deduce that there exists a $(v,x)$-spanned $x$-shifted iceberg $\J''$, for some $x \in \R^d$ and $v \in \{u\} \cup N_T^\to(u)$, and then by bounding the probability that such an iceberg exists by taking a union bound over choices of $\J''$. We shall use Lemma~\ref{lem:iceberg-count} and Remark~\ref{remark:choices} to bound the number of choices, and $\ih\big( s', x(t-1) \big)$ for some $s' \in \{s-1,s\}$ 
(which holds by~\eqref{eq:indstep-hyp}) to bound the probability that $\J''$ is $(v,x)$-spanned. 

\smallskip
\pagebreak

\begin{claim}\label{clm:indstep-c}
$\P_p(E_3) \leq q\big( \diam^{**}(\J) \big) / 3$.
\end{claim}

\begin{clmproof}{clm:indstep-c}
We shall again use~\eqref{eq:counting:icebergs:claims} and~\eqref{eq:applying:induction:claims} to count the possible icebergs $\J'$ and to bound the probability of the event $\ice_u(\J')$ for each of them. Our task will therefore be to bound the probability of the event $\Delta_u(\J',\J)$ for fixed $\J$ and $\J' = (J',u)$ such that $J' \subset J$ and the bounds~\eqref{eq:indstep-case-b-diams} hold.

Suppose that the event $\Delta_u(\J',\J)$ holds, and apply Lemma~\ref{lem:sideways-iceberg}. Observe that the conditions of the lemma are satisfied, since $d_T(\0,u) = r - s$ and $s \ge 3$, by assumption, and since 
$J' \subset J$ and $\gamma \ge y \ge \lambda$, by~\eqref{eq:indstep-case-b-diams} and~\eqref{eq:alpha:gamma}, and since $\diam^*(\J) > x(2) \gg 1$. By Lemma~\ref{lem:sideways-iceberg}, it follows that there exists $x \in \R^d$, an $x$-shifted iceberg $\J''$, and a set $K \subset (J \cap A) \setminus J'$, such that
$$\diam_x^*(\J'') \ge \delta^3 \gamma \qquad \text{and} \qquad \< K \>_{v,x} = \{\J''\}$$
for some $v \in \{u\} \cup N_T^\to(u)$. Recall that this implies that the event $\ice_{v,x}(\J'')$ holds. 


To bound the probability that there exists such an iceberg, observe first that, by Lemma~\ref{lem:AL}, we may assume that
\begin{equation}\label{eq:indstep-case-b2-diams2}
\delta^3 \gamma \le \diam_x^*(\J'') \le \gamma.
\end{equation}
Indeed, if $\diam_x^*(\J'') > \gamma$ then, by Lemma~\ref{lem:AL}, there exists a set $K' \subset K$ and an $x$-shifted iceberg $\J'''$, with $\gamma/3 \le \diam_x^*(\J''') \le \gamma$, such that $\< K' \>_{v,x} = \{ \J''' \}$. 

Now, 
to count the number of choices for the iceberg $\J'' = (J'',w)$ with $J'' \cap J \ne \emptyset$ and $\diam^*(\J'') \le \gamma \le \diam^*(\J)$, we would like to apply Lemma~\ref{lem:iceberg-count}. The lemma gives a bound on the number of icebergs, not $x$-shifted icebergs; however, by Remark~\ref{remark:choices}, 
$x$ is a deterministic function of the icebergs $\J$ and $\J'$ (which were fixed above) and the vertex $v \in \{u\} \cup N_T^\to(u)$. There are therefore only $|N_T^\to(u)| + 1$ possible choices for $x$ (given $\J$ and $\J'$), and so we have at most 
\begin{equation}\label{eq:number:of:xshifted:icebergs}
\big( |N_T^\to(u)| + 1 \big) \cdot \diam^*(\J)^{C} \le \big( \delta^k (\alpha + \gamma) \big)^{2C}
\end{equation}
choices for the iceberg $\J''$, since $\diam^*(\J) = \alpha + \gamma \gg 1$. 

Next we use the induction hypothesis to bound the probability of the event $\ice_{v,x}(\J'')$. Set $s' := r - d_T(\0,v)$, and observe that $s' \in \{s-1,s\}$, 
since $v \in \{u\} \cup N_T^\to(u)$. 
Now, recalling from~\eqref{eq:indstep-a},~\eqref{eq:indstep-case-b-diams},~\eqref{eq:alpha:gamma} and~\eqref{eq:indstep-case-b2-diams2} that 
$$\diam_x^*(\J'') \leq \gamma \leq 3y = 3 \big( \log \diam^*(\J) \big)^{1+\eta},$$
and that $\diam^*(\J) \le x(t) = \exp_{(t-2)}(p^{-\eps(t)})$, it follows that
\begin{equation}\label{eq:J''-clm3}
\diam_x^*(\J'') \leq 3 \big( \exp_{(t-3)}(p^{-\eps(t)}) \big)^{1+\eta} \leq \exp_{(t-3)}(p^{-\eps(t-1)}) - 1 = x(t - 1) - 1,
\end{equation}
since $\eps(t) < \eps(t-1)$, and by our choice of $\eta$. 

Now, recall from~\eqref{eq:indstep-hyp} that $\ih\big( s-1, x(t-1) \big)$ and $\ih(s,x-1)$ both hold, and also that $x > x(t-1)$, by our choice of $t$. Hence $\ih\big( s', x(t-1) - 1 \big)$ also holds, which together with~\eqref{eq:indstep-case-b2-diams2} and~\eqref{eq:J''-clm3} implies that
$$\P_p\big( \ice_{v,x}(\J'') \big) \le q\big( \diam_x^{**}(\J'') \big) \le q\big( \delta^{k+4} \gamma \big),$$ 
since $q$ is decreasing and $d_T(\0,v) \le k+1$, so $\diam_x^{**}(\J'') \ge \delta^{k+1} \diam_x^{*}(\J'') \ge \delta^{k+4} \gamma$. 

By the union bound, and recalling~\eqref{eq:number:of:xshifted:icebergs}, it follows that
\begin{equation}\label{eq:indstep-delta}
\P_p\big( \Delta_u(\J',\J) \big) \leq \big( \delta^k (\alpha + \gamma) \big)^{2C} \cdot q\big( \delta^{k+4} \gamma \big).
\end{equation}
Moreover, since $\diam^{**}(\J') = \delta^k \alpha$ by~\eqref{eq:alpha:gamma} and because $\J'$ is an iceberg of type $u$, it follows from~\eqref{eq:applying:induction:claims} that
\begin{equation}\label{eq:recalling:induction:claims}
\P_p \big( \ice_u(\J') \big) \le q\big( \diam^{**}(\J') \big) = q(\delta^k \alpha).
\end{equation}
Now, for each fixed $\J'$, the events $\ice_u(\J')$ and $\Delta_u(\J',\J)$ are independent, because they depend on disjoint subsets of $A$. Hence, by the union bound, and using~\eqref{eq:counting:icebergs:claims},~\eqref{eq:indstep-delta} and~\eqref{eq:recalling:induction:claims}, we obtain
\begin{align}
\P_p(E_3) & \le \sum_{\J'} \P_p \big( \ice_u(\J') \big) \cdot \P_p\big( \Delta_u(\J',\J) \big) \notag \\
&\leq \sup_{\alpha,\gamma} \; \big( \delta^k(\alpha + \gamma) \big)^{4C} \cdot q\big( \delta^k \alpha \big) \cdot q\big( \delta^{k+4} \gamma \big), \label{eq:indstep-case-b2-sup}
\end{align}
where the sum is over all valid choices of $\J'$.

Finally, we bound the right-hand side of~\eqref{eq:indstep-case-b2-sup} using Lemma~\ref{lem:sideways}, applied to $\delta^k \alpha$ and $\delta^k \gamma$. To check that the conditions of the lemma are satisfied, recall from~\eqref{eq:indstep-a},~\eqref{eq:indstep-case-b-diams} and~\eqref{eq:alpha:gamma} that $y = \big( \log( \alpha + \gamma ) \big)^{1+\eta}$ and $y \le \gamma \le 3y$. Moreover, $\alpha \gg \gamma \ge y \gg 1$, so by ~\eqref{eq:indstep-case-b2-sup} and Lemma~\ref{lem:sideways}, we have 
\begin{equation}\label{eq:indstep-bound3}
\P_p(E_3) \le \sup_{\alpha,\gamma} \; \frac{q\big( \delta^k (\alpha + \gamma) \big)}{\big( \delta^k (\alpha + \gamma) \big)^C} \le \frac{q\big( \diam^{**}(\J) \big)}{3}
\end{equation}
as required, since $\diam^{**}(\J) = \delta^k (\alpha + \gamma) \gg 1$. 
\end{clmproof}

As observed earlier, the event $\ice_u(\J)$ implies that the event $E_1 \cup E_2 \cup E_3$ holds, by Lemma~\ref{lem:goodsat}, and therefore~\eqref{eq:indstep-aim} follows from Claims~\ref{clm:indstep-a},~\ref{clm:indstep-b} and~\ref{clm:indstep-c}. This completes the proof of the lemma.
\end{proof}

Now that Proposition~\ref{prop:critdrop-v2} has been proved, all that remains before we can complete the proof of Theorem~\ref{thm:lower} is to verify that if $[A]_\U = \Z_n^d$  then $\Z_n^d$ contains a $\0$-iceberg spanned $\0$-iceberg of diameter roughly $x(r)$. We therefore now switch our setting to the torus, and for the rest of the proof we assume that $A \subset \Z_n^d$. 

The following Aizenman--Lebowitz-type lemma is not, strictly speaking, a consequence of Lemma~\ref{lem:AL}, because of the wrap-around effect of the torus. 
However, since we only need the lemma for droplets, the proof is particularly simple. 

\begin{lemma}\label{lem:AL-for-torus}
Suppose that $[A]_\U = \Z_n^d$, for some $n \gg x(r)$. Then there exists a $\0$-iceberg spanned $\0$-iceberg\/ $\J$ with 
$$\frac{x(r)}{3} \le \diam^*(\J) \le x(r).$$
\end{lemma}

\begin{proof}
We run the $\0$-iceberg spanning algorithm on $\Z_n^d$ with initial set $A$. Observe that, since $[A]_\0 = [A]_\U = \Z_n^d$, the algorithm cannot terminate before
$$\diam^*\big( \J_\0([K']_\U) \big) \ge x(r)$$
for some $t \in \N$ and $K' \in \K^t$. Now, by Lemma~\ref{lem:subadd}, the quantity
$$\max\Big\{ \diam^*\big( \J_\0\big( [K']_\U \big) \big) : K' \in \K^t \Big\}$$
at most triples at each step, provided that this maximum is at least an absolute constant (cf.~the proof of Lemma~\ref{lem:AL}). It follows that for some $t \le T$, there exists $K' \in \K^t$ such that $x(r)/3 \le \diam^*(\J) \le x(r)$, where $\J_\0([K']_\U) = \{\J\}$, as required. 
\end{proof}


We are finally ready to complete the proof of Theorem~\ref{thm:lower}.

\begin{proof}[Proof of Theorem~\ref{thm:lower}]
We are required to prove that for any $0 < \eps < \eps(r)$, if
\begin{equation}\label{eq:proof:p}
p \le \big( \log_{(r-1)} n \big)^{-1/\eps},
\end{equation}
then $\Pr_p\big( [A]_\U = \Z_n^d \big) \to 0$ as $n \to \infty$. Since, by~\eqref{eq:constants1}, the constant $\eps(r)$ may be chosen arbitrarily close to $1/(d-1)$, this suffices to prove the theorem. 

Let $A$ be a $p$-random subset of $\Z_n^d$, and suppose that $[A]_\U = \Z_n^d$. By Lemma~\ref{lem:AL-for-torus}, there~exists a $\0$-iceberg spanned $\0$-iceberg $\J$ such that
\begin{equation}\label{eq:proofdiam}
\frac{1}{3} \cdot \exp_{(r-2)} \big( p^{-\eps(r)} \big) \le \diam^*(\J) \leq \exp_{(r-2)} \big( p^{-\eps(r)} \big).
\end{equation}
Now, the number of $\0$-icebergs in $\Z_n^d$ 
is $n^{O(1)}$, since we have at most $n^d$ choices for each face, 
and each is $\0$-iceberg spanned with probability at most
$$\exp\big( -\diam^*(\J) \big) \leq \exp\Big( -\tfrac{1}{3} \cdot \exp_{(r-2)} \big( p^{-\eps(r)} \big) \Big),$$
by Proposition~\ref{prop:critdrop-v2}. 
Thus, by~\eqref{eq:proof:p} and since $\eps < \eps(r)$, we have finally that
$$\P_p\big( [A] = \Z_n^d \big) \le n^{O(1)} \cdot \exp\Big( -\tfrac{1}{3} \cdot \exp_{(r-2)} \big( p^{-\eps(r)} \big) \Big) \to 0$$
as $n \to \infty$, and this completes the proof of Theorem~\ref{thm:lower}.
\end{proof}

\section*{Acknowledgements}

The authors would like to thank Hugo Duminil-Copin for a very helpful discussion at the outset of this project, and for many other interesting conversations over the years.

\bibliographystyle{amsplain}
\bibliography{bprefs}

\appendix

\section{Simple calculations}\label{app:calcs}

This appendix contains the proofs of three simple bounds that were used in Section~\ref{sec:ind}, namely Lemmas~\ref{lem:splitter},~\ref{lem:giveaway} and~\ref{lem:sideways}. For the reader's convenience, we restate each of these here, as Lemmas~\ref{lem:splitter-app},~\ref{lem:giveaway-app} and~\ref{lem:sideways-app}, respectively. All constants used in this appendix should be interpreted as they were used in Section~\ref{sec:ind}. Recall that $p$ is assumed to be sufficiently small, so in particular $\log 1/p \gg C$. 

\begin{lemma}\label{lem:splitter-app}
Let $r \geq 3$, and let $\alpha \ge \beta \gg 1$ and $\gamma \in \R$ satisfy 
$$\gamma \le C (\log \alpha)^{1+\eta} \le C^2 \beta.$$ 
Then
\begin{equation}\label{eq:splitter}
q(\alpha) \cdot q(\beta) \cdot (\alpha + \beta + \gamma)^{4C} \le q(\alpha + \beta + \gamma).
\end{equation}
\end{lemma}

\begin{proof}
The left-hand side of~\eqref{eq:splitter} is increasing in $\gamma$ and the right-hand side is decreasing in $\gamma$ (since $q$ is decreasing), so we may assume that $\gamma = C(\log \alpha)^{1+\eta} > 0$. Recalling the definition~\eqref{eq:qx} of $q(x)$, we are required to show that
\begin{equation}\label{eq:splitterstp}
\frac{\alpha}{\log_{(r-1)}( \lambda \alpha )} + \frac{\beta}{\log_{(r-1)}( \lambda \beta )}
\geq \frac{\alpha+\beta+\gamma}{\log_{(r-1)}\big( \lambda (\alpha+\beta+\gamma) \big)} + \frac{4C \log(\alpha+\beta+\gamma)}{\log 1/p}.
\end{equation}
To prove~\eqref{eq:splitterstp}, observe first that if
\begin{equation}\label{eq:splitter:case1}
\log_{(r-1)}\big( \lambda (\alpha+\beta+\gamma) \big) \ge C^3 \cdot \log_{(r-1)}( \lambda \beta ),
\end{equation}
then, since $\gamma \le C^2\beta$, we have 
$$\frac{\beta}{\log_{(r-1)}( \lambda \beta )} \,-\, \frac{\beta + \gamma}{\log_{(r-1)}\big( \lambda (\alpha+\beta+\gamma) \big)} \ge  \frac{\beta}{2 \log_{(r-1)}( \lambda \beta )}.$$
Now, since $\beta \ge C^{-1} (\log \alpha)^{1+\eta} \gg 1$ and $r \ge 2$, we have
$$\frac{\beta}{\log_{(r-1)}( \lambda \beta )} \ge \log \alpha \gg \frac{\log(\alpha+\beta+\gamma)}{\log 1/p},$$
so in this case~\eqref{eq:splitterstp} holds. On the other hand, if~\eqref{eq:splitter:case1} fails to hold, then it follows 
that $\beta \ge (\log \alpha)^{2C}$, and therefore
\begin{equation}\label{eq:splitterbeats}
\frac{\beta}{(\log \alpha)^C} \ge (\log \alpha)^C \gg \frac{\gamma}{\log_{(r-1)}\big( \lambda (\alpha+\beta+\gamma) \big)} + \frac{4C \log(\alpha+\beta+\gamma)}{\log 1/p},
\end{equation}
since $\gamma \le C(\log \alpha)^{1+\eta}$. It will therefore suffice to  show that 
\begin{equation}\label{eq:splitterstp1}
\frac{\alpha}{\log_{(r-1)}( \lambda \alpha )} - \frac{\alpha}{\log_{(r-1)}\big( \lambda (\alpha+\beta+\gamma) \big)} \geq \frac{\beta}{( \log \alpha )^C}.
\end{equation}
To see this, note first that, by the mean value theorem, we have
$$\log_{(r-1)}\big( \lambda (\alpha+\beta+\gamma) \big) - \log_{(r-1)}( \lambda \alpha ) = \frac{\beta + \gamma}{(\alpha + \xi) \log\big(\lambda (\alpha+\xi) \big) \cdots \log_{(r-2)}\big(\lambda (\alpha+\xi)\big)}$$
for some $\xi \in (0,\beta+\gamma)$. Now, observe that the right-hand side is at least
$$\frac{\beta}{C^4 \alpha \cdot (\log \alpha)^{r-2}} \ge \frac{\beta}{\alpha \cdot (\log \alpha)^{C-3}},$$
since $\alpha + \xi < \alpha + \beta + \gamma \le C^3 \alpha$ and $\gamma \ge 0$, and recalling that $\alpha \gg 1$ and $C$ is a large constant. Rearranging, we obtain~\eqref{eq:splitterstp1}, and as noted above this proves~\eqref{eq:splitterstp}.\end{proof}

\begin{lemma}\label{lem:giveaway-app}
Let $\alpha \ge \gamma \gg 1$. Then
\[
q(\alpha) \cdot (\alpha + \gamma)^{3C} \le q\big( \delta(\alpha + \gamma) \big).
\]
\end{lemma}

\begin{proof}
Using~\eqref{eq:qx}, the conclusion of the lemma may be rewritten as
\begin{equation}\label{eq:giveaway}
\frac{\alpha}{\log_{(r-1)}( \lambda \alpha )} \geq \frac{\delta(\alpha + \gamma)}{ \log_{(r-1)}\big( \lambda \delta(\alpha + \gamma) \big)} + \frac{3C \log(\alpha + \gamma)}{\log 1/p}.
\end{equation}
Observe that $\log_{(r-1)}( \delta \lambda \alpha ) \ge 4\delta \cdot \log_{(r-1)}( \lambda \alpha )$, since $\alpha \gg 1$ and $\delta < 1/4$, and therefore
\[
\frac{\alpha}{\log_{(r-1)}( \lambda \alpha )} - \frac{\delta(\alpha + \gamma)}{ \log_{(r-1)}\big( \lambda \delta(\alpha + \gamma) \big)} \geq \frac{\alpha}{2 \log_{(r-1)}( \lambda \alpha )},
\]
since $\alpha \ge \gamma$. The bound~\eqref{eq:giveaway} now follows easily.
\end{proof}

\begin{lemma}\label{lem:sideways-app}
Let $r \geq 3 $ and $\alpha,\gamma \gg 1$ be such that
\begin{equation}\label{eq:sideways-z2}
C^{-1} \big( \log(\alpha + \gamma) \big)^{1 + \eta} \le \gamma \le C \big( \log(\alpha + \gamma) \big)^{1 + \eta}.
\end{equation}
Then
\[
q(\alpha) \cdot q(\delta^4 \gamma) \cdot (\alpha + \gamma)^{5C} \leq q(\alpha + \gamma).
\]
\end{lemma}

\begin{proof}
Using~\eqref{eq:qx}, the conclusion of the lemma may be rewritten as
$$\frac{\alpha}{\log_{(r-1)}( \lambda \alpha )} + \frac{\delta^4 \gamma}{\log_{(r-1)}( \lambda \delta^4 \gamma )} \ge \frac{\alpha+\gamma}{\log_{(r-1)}\big( \lambda (\alpha+\gamma) \big)} + \frac{5C \log(\alpha+\gamma)}{\log 1/p}.$$
Since $\gamma > 0$ and $5C < \log(1/p)$, this will follow if we can prove that
\begin{equation}\label{eq:sidewaysstp}
\frac{\delta^4 \gamma}{\log_{(r-1)}( \lambda \delta^4 \gamma )} - \frac{\gamma}{\log_{(r-1)}\big( \lambda (\alpha+\gamma) \big)} \ge \log(\alpha+\gamma).
\end{equation}
To prove~\eqref{eq:sidewaysstp}, note that, since $\alpha,\gamma \gg 1$, it follows from~\eqref{eq:sideways-z2} that  
$$\delta^4 \cdot \log_{(r-1)}\big(\lambda(\alpha+\gamma)\big) \ge 2 \cdot \log_{(r-1)}(\lambda \delta^4 \gamma),$$
and therefore
\[
\frac{\delta^4 \gamma}{\log_{(r-1)}(\lambda \delta^4 \gamma)} - \frac{\gamma}{\log_{(r-1)}(\lambda(\alpha+\gamma))} \ge \frac{\delta^4 \gamma}{2 \cdot \log_{(r-1)}(\lambda \delta^4 \gamma)}.
\]
Since $\gamma \ge C^{-1} \big( \log(\alpha + \gamma) \big)^{1+\eta}$, by~\eqref{eq:sideways-z2}, we obtain~\eqref{eq:sidewaysstp}, as required.
\end{proof}

\end{document}